\documentclass[10pt]{amsart}

\usepackage{float}
\usepackage{amscd,amsmath,amssymb,amsfonts,amsthm,ascmac}
\usepackage{enumerate,mathrsfs,stmaryrd,latexsym, comment, mathtools, mathdots} 
\usepackage{colonequals}

\usepackage{amsthm, amssymb}
\usepackage{bm}
\usepackage[leqno]{amsmath}
\usepackage{caption}
\usepackage{subcaption}
\usepackage{hyperref}
\usepackage{relsize}
\usepackage{booktabs}
\usepackage{a4wide}
\usepackage{ifthen}
\usepackage[mathscr]{euscript}
\usepackage{color}

\usepackage[all]{xy}

\usepackage{tikz-cd}
\usepackage{tikz}
\usetikzlibrary{snakes,3d, matrix,decorations.pathreplacing,calc,decorations.pathmorphing, fit, patterns, trees}
\usetikzlibrary {positioning}

\numberwithin{equation}{section}
\allowdisplaybreaks[1]

\allowdisplaybreaks

\newcommand{\C}{{\mathbb{C}}}
\newcommand{\R}{{\mathbb{R}}}
\newcommand{\Z}{{\mathbb{Z}}}

\newcommand{\Paths}{\mathcal{GP}}
\newcommand{\RCn}{R(w_0^{(C_n)})}

\newcommand{\Rcn}[1]{R(w_0^{(C_{#1})})}
\newcommand{\Ran}[1]{R(w_0^{(A_{#1})})}

\newcommand{\chamber}{\mathscr{C}}
\newcommand{\F}{\mathscr{F}}
\newcommand{\mcal}{\mathcal}
\newcommand{\vs}{\vspace}

\def\node{\textrm{node}}
\def\symp{\rm{symp}}

\newtheorem{theorem}{Theorem}[section]
\newtheorem{lemma}[theorem]{Lemma}
\newtheorem{proposition}[theorem]{Proposition}
\newtheorem{corollary}[theorem]{Corollary}

\theoremstyle{definition}
\newtheorem{example}[theorem]{Example}
\newtheorem{definition}[theorem]{Definition}
\newtheorem{remark}[theorem]{Remark}

\hypersetup{colorlinks}
\hypersetup{
colorlinks=true,       
linkcolor=blue,          
citecolor=red,        
filecolor=violet,      
urlcolor=violet       
}

\begin{document}

\author{Yunhyung Cho}
\address{Department of Mathematics Education, Sungkyunkwan University, Seoul, Republic of Korea}
\email{yunhyung@skku.edu}

\author{Naoki Fujita}
\address{Faculty of Advanced Science and Technology, Kumamoto University, 2-39-1 Kurokami, Chuo-ku, Kumamoto 860-8555, Japan}
\email{fnaoki@kumamoto-u.ac.jp}

\author{Eunjeong Lee}
\address{Department of Mathematics, Chungbuk National University, 28644, Republic of Korea}
\email{eunjeong.lee@chungbuk.ac.kr}

\thanks{Cho was supported by the National Research Foundation of Korea(NRF) grant funded by the Korea government(MSIP; Ministry of Science, ICT \& Future Planning) (No.2020R1C1C1A01010972) and (No.2020R1A5A1016126). Fujita was supported by Grant-in-Aid for JSPS Fellows (No.\ 19J00123), by JSPS Grant-in-Aid for Early-Career Scientists (No.\ 20K14281), and by MEXT Japan Leading Initiative for Excellent Young Researchers (LEADER) Project. Lee was supported by the National Research Foundation of Korea(NRF) grant funded by the Korea government(MSIT) (No. RS-2022-00165641 and No. RS-2023-00239947). }

\title{On combinatorics of string polytopes in types $B$ and $C$}

\date{\today}



\begin{abstract} 
A string polytope is a rational convex polytope whose lattice points parametrize a highest weight crystal basis, which is obtained from a string cone by explicit affine inequalities depending on a highest weight. 
It also inherits geometric information of a flag variety such as toric degenerations, Newton--Okounkov bodies, mirror symmetry, Schubert calculus, and so on.
In this paper, we study combinatorial properties of string polytopes in types $B$ and $C$ by giving an explicit description of string cones in these types which is analogous to Gleizer--Postnikov's description of string cones in type $A$. 
As an application, we characterize string polytopes in type~$C$ which are unimodularly equivalent to the Gelfand--Tsetlin polytope in type $C$ for a specific highest weight.
\end{abstract}

\maketitle
\setcounter{tocdepth}{1} \tableofcontents

\section{Introduction}
A crystal basis is a combinatorial skeleton of a representation of a semisimple Lie algebra $\mathfrak{g}$, which was introduced by Kashiwara \cite{Kas90, Kas91} via the quantized enveloping algebra $U_q(\mathfrak{g})$ associated with $\mathfrak{g}$.
To study crystal bases, it is important to give their concrete parameterizations. 
In the present paper, we focus on a specific polyhedral parametrization, called Berenstein--Littelmann--Zelevinsky's string polytope (see \cite{Littelmann98, BerensteinZelevinsky01}). 
Let $W$ be the Weyl group of $\mathfrak{g}$, $w_0 \in W$ the longest element, $R(w_0)$ the set of reduced words for $w_0$, and $P_+$ the set of dominant integral weights. 
A string polytope $\Delta_{\bm i}(\lambda)$ is a rational convex polytope defined from ${\bm i} \in R(w_0)$ and $\lambda \in P_+$, which is obtained from a string cone $\mathcal{C}_{\bm i}$ by explicit affine inequalities depending on~$\lambda$.
In type $X_n$, where $X$ is $A$ or $C$, Littelmann \cite{Littelmann98} proved that there exists a specific reduced word ${\bm i}_X \in R(w_0)$ such that the string polytope $\Delta_{{\bm i}_X}(\lambda)$ is unimodularly equivalent to the Gelfand--Tsetlin polytope $GT_{X_n}(\lambda)$ for each $\lambda \in P_+$.
String polytopes also inherit geometric information of the full flag variety associated with $\mathfrak{g}$ such as toric degenerations \cite{GL96, Cal02, KM05}, Newton--Okounkov bodies \cite{Kav15, FO17, FO_cluster}, mirror symmetry \cite{BCFKvS00, AB04, Rusinko08}, Schubert calculus \cite{KST12, Fujita_Schubert}, and so on. 
Hence it is interesting to study combinatorics of string polytopes. 
Since combinatorial properties of $\Delta_{\bm i}(\lambda)$ heavily depend on the choice of a reduced word~${\bm i} \in R(w_0)$, it is a fundamental problem to classify the string polytopes $\{\Delta_{\bm i}(\lambda)\}_{\bm i}$ up to unimodular equivalence, where $\lambda$ is fixed. 
In a joint work \cite{CKLP21_GC} with Kim and Park, the first and third named authors addressed this problem in type $A_n$.
More precisely, using Gleizer--Postnikov's description~\cite{GleizerPostnikov} of string cones in type $A_n$, they classified reduced words ${\bm i} \in R(w_0)$ such that $\Delta_{\bm i}(\lambda)$ is unimodularly equivalent to the Gelfand--Tsetlin polytope $GT_{A_n}(\lambda)$ for all $\lambda \in P_+$.
In the present paper, we address the problem in types $B_n$ and $C_n$ by giving an explicit description of string cones which is analogous to the description by Gleizer--Postnikov.

String cones in types $A_{2n-1}$, $B_n$, and $C_n$ are closely related as follows. 
Crystal bases in type~$B_n$ can be realized as specific subsets of crystal bases in type $A_{2n-1}$ (see \cite{Kas96, NS05}).
From this, we know that string cones in type $B_n$ are identified with slices of string cones in type $A_{2n-1}$ (see Theorem \ref{t:relation_between_A_and_B} for more details). 
In addition, Kashiwara \cite{Kas96} gave a similarity between crystal bases in types $B_n$ and $C_n$, which induces a \emph{similarity} between their string cones (see Theorem~\ref{t:similarity_between_B_and_C}). 
The second named author \cite{Fujita18} also proved that string cones in type $C_n$ can be obtained as quotients of string cones in type $A_{2n-1}$ (see Theorem~\ref{t:relation_between_A_and_C} for more details).
Using these relations, we show that Gleizer--Postnikov's description \cite{GleizerPostnikov} in type $A_{2n-1}$ induces an explicit description of string cones in types $B_n$ and $C_n$. 
As an application of the description, we classify simplicial string cones in types $B_n$ and $C_n$ as follows; this is the main result of the present paper. 

\begin{theorem}[{see Theorem~\ref{thm_simplicial_string_cones_BC} and Remark~\ref{r:second_main_type_B}}]\label{thm_intro_2}
Let $\mathfrak{g}$ be a simple Lie algebra of type $B_n$ or $C_n$ with $n \geq 2$. 
Then, for ${\bm i} \in R(w_0)$, the following are equivalent.
\begin{enumerate}
\item The number of facets\footnote{A \emph{facet} of an $N$-dimensional polytope is a face that has dimension $N-1$.} of $\Delta_{\bm i}(\lambda)$ is $2N$ for every $\lambda \in P_{++}$, where $P_{++} \subseteq P_+$ denotes the set of regular dominant integral weights. 
\item The string cone $\mathcal{C}_{\bm i}$ is simplicial.
\item The reduced word ${\bm i}$ is either 
\[
\begin{split}
{\bm i}_{C}^{(n)} &\coloneqq (n, \underbrace{n-1, n, n-1}_{3}, \underbrace{n-2, n-1, n, n-1, n-2}_{5}, \ldots, \underbrace{2, \ldots, n, \ldots, 2}_{2n-3}, \underbrace{1, 2, \ldots, n, \ldots, 2, 1}_{2n-1}); \text{ or} \\
{\bm j}_{C}^{(n)} &\coloneqq (n-1, n, n-1, n, \underbrace{n-2, n-1, n, n-1, n-2}_{5}, \ldots, \underbrace{2, \ldots, n, \ldots, 2}_{2n-3}, \underbrace{1, 2, \ldots, n, \ldots, 2, 1}_{2n-1}).
\end{split}
\]
\end{enumerate}
\end{theorem}

Let $\rho \in P_{++}$ be the sum of fundamental weights. 
In type $C_n$, we also classify string polytopes~$\Delta_{\bm i}(\rho)$ which are unimodularly equivalent to $GT_{C_n}(\rho)$. 

\begin{theorem}[{Theorem~\ref{thm_GT_string_polytopes_C}}]\label{thm_intro_3}
Let $\mathfrak{g}$ be the simple Lie algebra of type $C_n$ with $n \geq 2$, and ${\bm i} \in R(w_0)$.
Then the string polytope $\Delta_{\bm i}(\rho)$ is unimodularly equivalent to the Gelfand--Tsetlin polytope $GT_{C_n}(\rho)$ if and only if ${\bm i} = {\bm i}_{C}^{(n)}$.
\end{theorem}

The present paper is organized as follows. 
In Section~\ref{sec_GP_description}, we recall some basic definitions on Gleizer--Postnikov paths, and review their description of string cones in type $A$. 
Section~\ref{sec_folding} is devoted to recalling some relations among string cones in types $A_{2n-1}$, $B_n$, and $C_n$, which induce systems of explicit linear inequalities defining string cones in types $B_n$ and $C_n$.
In Section~\ref{sec_symmetric_GP_paths}, we study non-redundancy of the inequalities.  
Section~\ref{sec_GCtype} is devoted to proving Theorems~\ref{thm_intro_2} and~\ref{thm_intro_3} above.

\section{Gleizer--Postnikov description of string cones}
\label{sec_GP_description}
Let $\mathfrak{g} \coloneqq \mathfrak{sl}_m(\mathbb{C})$ be the special linear Lie algebra over $\mathbb{C}$, and write $[k] \coloneqq \{1, 2, \ldots, k\}$ for $k \in \mathbb{Z}_{>0}$.
We identify the set $I$ of vertices of the Dynkin diagram of $\mathfrak{g}$ with $[m-1]$ as follows:
\begin{align*}
&A_{m-1}\ \begin{xy}
\ar@{-} (50,0) *++!D{1} *\cir<3pt>{};
(60,0) *++!D{2} *\cir<3pt>{}="C"
\ar@{-} "C";(70,0) *++!D{3} *\cir<3pt>{}="D"
\ar@{-} "D";(75,0) \ar@{.} (75,0);(80,0)^*!U{}
\ar@{-} (80,0);(85,0) *++!D{m-1} *\cir<3pt>{}="E"
\end{xy}.
\end{align*}
Let $\mathfrak{S}_{m}$ be the symmetric group on $[m]$, which we identify with the Weyl group of $\mathfrak{sl}_{m}(\C)$.
Denote by $w_0^{(A_{m-1})} \in \mathfrak{S}_{m}$ the longest element, and by $\ell$ the length of $w_0^{(A_{m-1})}$, that is, we have $\ell = \frac{m(m-1)}{2}$. 
Let 
\[
R(w_0^{(A_{m-1})}) \colonequals \{{\bm i} = (i_1, \dots, i_{\ell}) \in [m-1]^{\ell} \mid 
w_0^{(A_{m-1})} = s_{i_1} s_{i_2} \cdots s_{i_{\ell}}\}
\]
be the set of reduced words for $w_0^{(A_{m-1})}$, where $s_i$ is the simple transposition $(i, i+1)$ for $i \in [m-1]$. 
For each reduced word ${\bm i} \in R(w_0^{(A_{m-1})})$, 
one can associate a rational convex polyhedral cone $\mathcal{C}_{\bm i}^{(A_{m-1})} \subseteq \R^{\ell}$ as in Littelmann \cite{Littelmann98} and Berenstein--Zelevinsky \cite{BerensteinZelevinsky01}, which is called a {\em string cone}.

In this section, we shall review Gleizer--Postnikov's description of string cones in type $A_{m-1}$ using wiring diagrams 
(see \cite{GleizerPostnikov}). 
Note that there is another description of string cone inequalities using \textit{rhombic tilings of a regular $(2m)$-gon} (see~\cite[Remarks 3.12 and 5.5]{GKS21} and references therein).	

\subsection{Wiring diagrams and rigorous paths }
\label{ssecWiringDiagrams}\label{ssecRigorousPaths}
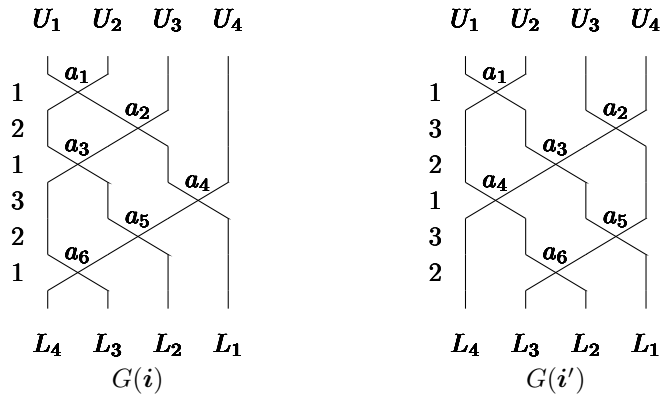
\begin{figure}[h]
\begin{tikzpicture}[scale = 0.8, yscale=0.6]
\def\nwires{4} 
\def\length{6} 

\foreach \i\y in {1/1,2/2,3/1,4/3,5/2,6/1}{
	\pgfmathsetmacro\ycoordinate{\length-\i+1}
	\pgfmathsetmacro\yplusone{\y+1}
	\foreach \x in {1,...,\nwires}{
		\ifthenelse{\x = \y}
			{\coordinate (\i SW) at (\x, \ycoordinate+0);
			\coordinate (\i SE) at (\x+1, \ycoordinate+0);
			\coordinate (\i NW) at (\x, \ycoordinate+1);
			\coordinate (\i NE) at (\x+1,\ycoordinate+1);

			\draw (\x,\ycoordinate) -- (\i SW)-- (\i NE) -- (\x+1, \ycoordinate+1);
			\draw (\x+1,\ycoordinate) --(\i SE) --(\i NW) --(\x, \ycoordinate+1);
			}{}
		\ifthenelse{\x=\y \OR \x = \yplusone}{}{\draw (\x,\ycoordinate)--(\x,\ycoordinate+1);}

	\node at (0.5,\ycoordinate+0.5) {$\y$};
	\foreach \x in {1,...,\nwires}{
		\node at (\x, \length+2.5) {$U_{\x}$};
		\node at (\nwires-\x+1,-0.5) {$L_{\x}$};

	\node[above] at (\y+0.5,\ycoordinate+0.5) {$a_{\i}$};
}}}

\foreach \x in{1,2,3,4}{
\draw (\x,\length+1)--(\x,\length+1.5);
\draw (\x,1)--(\x,0.5);
}

\node at (2.5,-1.5) {$G(\bm i)$};
\end{tikzpicture}
\hspace{2cm}
\begin{tikzpicture}[scale = 0.8, yscale=0.6]
\def\nwires{4} 
\def\length{6} 

\foreach \i\y in {1/1,2/3,3/2,4/1,5/3,6/2}{
	\pgfmathsetmacro\ycoordinate{\length-\i+1}
	\pgfmathsetmacro\yplusone{\y+1}
	\foreach \x in {1,...,\nwires}{
		\ifthenelse{\x = \y}
			{\coordinate (\i SW) at (\x, \ycoordinate+0);
			\coordinate (\i SE) at (\x+1, \ycoordinate+0);
			\coordinate (\i NW) at (\x, \ycoordinate+1);
			\coordinate (\i NE) at (\x+1,\ycoordinate+1);

			\draw (\x,\ycoordinate) -- (\i SW)-- (\i NE) -- (\x+1, \ycoordinate+1);
			\draw (\x+1,\ycoordinate) --(\i SE) --(\i NW) --(\x, \ycoordinate+1);
			}{}
		\ifthenelse{\x=\y \OR \x = \yplusone}{}{\draw (\x,\ycoordinate)--(\x,\ycoordinate+1);}

	\node at (0.5,\ycoordinate+0.5) {$\y$};
	\foreach \x in {1,...,\nwires}{
		\node at (\x, \length+2.5) {$U_{\x}$};
		\node at (\nwires-\x+1,-0.5) {$L_{\x}$};

	\node[above] at (\y+0.5,\ycoordinate+0.5) {$a_{\i}$};
}}}

\foreach \x in{1,2,3,4}{
\draw (\x,\length+1)--(\x,\length+1.5);
\draw (\x,1)--(\x,0.5);
}

\node at (2.5,-1.5) {$G(\bm i')$};

\end{tikzpicture}

		\caption{\label{figure_wiring_diagram_GC} Wiring diagrams for ${\bm i} = (1,2,1,3,2,1)$ and ${\bm i}' = (1,3,2,1,3,2)$.}	
	\end{figure}

Each reduced word ${\bm i} = (i_1, \dots, i_{\ell}) \in R(w_0^{(A_{m-1})})$ can be represented by a {\em wiring diagram} (that is also called a {\em pseudoline arrangement} or a {\em string diagram}). 
The diagram corresponding to a reduced word ${\bm i}$ is denoted by~$G({\bm i})$.
As depicted in Figure \ref{figure_wiring_diagram_GC}, 
the wiring diagram $G({\bm i})$ consists of a family of~$m$ vertical piecewise straight lines labeled by $\ell_1, \ell_2, \ldots, \ell_{m}$. 
For each $k \in [m]$, the upper end and the lower end of the line~$\ell_k$ are labeled by $U_k$ and $L_k$, respectively.
We call $\ell_k$ the $k$th \emph{wire}.
The crossing patterns of pairs of wires are determined by ${\bm i}$.
Specifically, the position of the $j$th crossing (from the top) should be located on the $i_j$th column of $G({\bm i})$ (see Figure \ref{figure_wiring_diagram_GC}). 
We call each crossing a {\em node} and name them as $a_1, a_2, \ldots, a_{\ell}$ from the top to the bottom. 

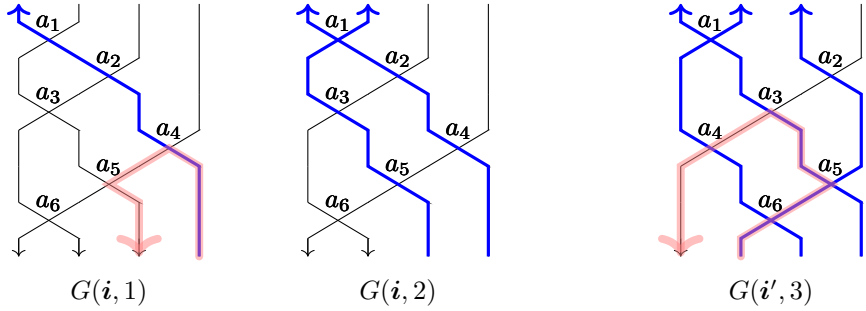
\begin{figure}[h]
\begin{tikzpicture}[scale = 0.8, yscale=0.6]
\def\nwires{4} 
\def\length{6} 

\foreach \i\y in {1/1,2/2,3/1,4/3,5/2,6/1}{
	\pgfmathsetmacro\ycoordinate{\length-\i+1}
	\pgfmathsetmacro\yplusone{\y+1}
	\foreach \x in {1,...,\nwires}{
		\ifthenelse{\x = \y}
			{\coordinate (\i SW) at (\x, \ycoordinate+0);
			\coordinate (\i SE) at (\x+1, \ycoordinate+0);
			\coordinate (\i NW) at (\x, \ycoordinate+1);
			\coordinate (\i NE) at (\x+1,\ycoordinate+1);
			
			\coordinate (\i C) at (\x+0.5, \ycoordinate + 0.5);

			\draw (\x,\ycoordinate) -- (\i SW)-- (\i NE) -- (\x+1, \ycoordinate+1);
			\draw (\x+1,\ycoordinate) --(\i SE) --(\i NW) --(\x, \ycoordinate+1);
			}{}
		\ifthenelse{\x=\y \OR \x = \yplusone}{}{\draw (\x,\ycoordinate)--(\x,\ycoordinate+1);}

\node[above] at (\y+0.5,\ycoordinate+0.5) {$a_{\i}$};
}}

\foreach \x in{1,...,\nwires}{
\draw (\x,\length+1)--(\x,\length+1.5);
\draw (\x,1)--(\x,0.5);
}

\draw[color=blue, very thick, ->] (4,0.5)--(4SE)--(4NW)--(2SE)--(2NW)--(1SE)--(1NW)--(1,7.5);

\foreach \x in {1,2,3}{
	\draw[->] (\x,1)--(\x,0.5);
}

\draw[color=red!50, line cap=round, line width=3, opacity=0.5, ->] (4,0.5)--(4SE)--(4C)--(4SW)--(5NE)--(5C)--(5SE)--(3,0.5);

\node at (2.5,-0.5) {$G(\bm i,1)$};
\end{tikzpicture}
\hspace{1cm} 
\begin{tikzpicture}[scale = 0.8, yscale=0.6]
\def\nwires{4} 
\def\length{6} 

\foreach \i\y in {1/1,2/2,3/1,4/3,5/2,6/1}{
	\pgfmathsetmacro\ycoordinate{\length-\i+1}
	\pgfmathsetmacro\yplusone{\y+1}
	\foreach \x in {1,...,\nwires}{
		\ifthenelse{\x = \y}
			{\coordinate (\i SW) at (\x, \ycoordinate);
			\coordinate (\i SE) at (\x+1, \ycoordinate);
			\coordinate (\i NW) at (\x, \ycoordinate+1);
			\coordinate (\i NE) at (\x+1,\ycoordinate+1);

			\coordinate (\i C) at (\x+0.5, \ycoordinate + 0.5);

			\draw (\x,\ycoordinate) -- (\i SW)-- (\i NE) -- (\x+1, \ycoordinate+1);
			\draw (\x+1,\ycoordinate) --(\i SE) --(\i NW) --(\x, \ycoordinate+1);
			}{}
		\ifthenelse{\x=\y \OR \x = \yplusone}{}{\draw (\x,\ycoordinate)--(\x,\ycoordinate+1);}

\node[above] at (\y+0.5,\ycoordinate+0.5) {$a_{\i}$};

}}

\foreach \x in{1,...,\nwires}{
\draw (\x,\length+1)--(\x,\length+1.5);
\draw (\x,1)--(\x,0.5);
}

\draw[color=blue, very thick, ->] (4,0.5)--(4SE)--(4NW)--(2SE)--(2NW)--(1SE)--(1NW)--(1,7.5);
\draw[color=blue, very thick, ->] 
(3,0.5)--(5SE)--(5NW)--(3SE)--(3NW)--(1SW)--(1NE)--(2,7.5);

\foreach \x in {1,2}{
	\draw[->] (\x,1)--(\x,0.5);
}

\node at (2.5,-0.5) {$G(\bm i,2)$};
\end{tikzpicture}
\hspace{2cm} 
\begin{tikzpicture}[scale = 0.8, yscale=0.6]
\def\nwires{4} 
\def\length{6} 

\foreach \i\y in {1/1,2/3,3/2,4/1,5/3,6/2}{
	\pgfmathsetmacro\ycoordinate{\length-\i+1}
	\pgfmathsetmacro\yplusone{\y+1}
	\foreach \x in {1,...,\nwires}{
		\ifthenelse{\x = \y}
			{\coordinate (\i SW) at (\x, \ycoordinate);
			\coordinate (\i SE) at (\x+1, \ycoordinate);
			\coordinate (\i NW) at (\x, \ycoordinate+1);
			\coordinate (\i NE) at (\x+1,\ycoordinate+1);

			\coordinate (\i C) at (\x+0.5, \ycoordinate + 0.5);

			\draw (\x,\ycoordinate) -- (\i SW)-- (\i NE) -- (\x+1, \ycoordinate+1);
			\draw (\x+1,\ycoordinate) --(\i SE) --(\i NW) --(\x, \ycoordinate+1);
			}{}
		\ifthenelse{\x=\y \OR \x = \yplusone}{}{\draw (\x,\ycoordinate)--(\x,\ycoordinate+1);}

\node[above] at (\y+0.5,\ycoordinate+0.5) {$a_{\i}$};

}}

\foreach \x in{1,...,\nwires}{
\draw (\x,\length+1)--(\x,\length+1.5);
\draw (\x,1)--(\x,0.5);
}

\draw[color=blue, very thick, ->] 
(4,0.5)--(5SE)--(5NW)--(3SE)--(3NW)--(1SE)--(1NW)--(1,7.5);
\draw[color=blue, very thick, ->]
(3,0.5)--(6SE)--(6NW)--(4SE)--(4NW)--(1SW)--(1NE)--(2,7.5); 
\draw[color=blue, very thick, ->]
(2,0.5)--(6SW)--(6NE)--(5SW)--(5NE)--(2SE)--(2NW)--(3,7.5);

\foreach \x in {1}{
	\draw[->] (\x,1)--(\x,0.5);
}

\draw[color=red!50, line cap=round, line width=3, opacity=0.5, ->] 
(2,0.5)--(6SW)--(6NE)--(5SW)--(5C)--(5NW)--(3SE)--(3C)--(3SW)--(4NE)--(4SW)--(1,0.5);

\node at (2.5,-0.5) {$G({\bm i}',3)$};
\end{tikzpicture}
\caption{\label{figure_wd_oriented} Oriented wiring diagrams for ${\bm i} = (1,2,1,3,2,1)$ and ${\bm i}' = (1,3,2,1,3,2)$.}
\end{figure}

A {\em rigorous path} is an oriented path on $G({\bm i})$ defined as follows. 
For each $k \in [m-1]$, let $G({\bm i}, k)$ be the wiring diagram $G({\bm i})$ together with the orientation on the wires, where the first $k$ wires $\ell_1, \ldots, \ell_k$ are oriented 
upward and the other wires $\ell_{k+1}, \ldots, \ell_{m}$ are oriented downward (see Figure~\ref{figure_wd_oriented}). 
	
\begin{definition}[{\cite[Section~5.1]{GleizerPostnikov}}]\label{definition_rigorous_path}
	For each $k \in [m-1]$, a {\em rigorous path} (or a {\em Gleizer--Postnikov path}) is an oriented path on $G({\bm i}, k)$ obeying the following properties:
\begin{itemize}
	\item it starts at $L_k$ and ends at $L_{k+1}$, 
	\item it respects the orientation of $G({\bm i}, k)$,
	\item it passes through each node at most once, and
	\item it does \emph{not} include a {\em forbidden fragment} given in Figure~\ref{figure_avoiding}.
\end{itemize}
\begin{figure}[h]
\begin{tikzpicture}
\tikzset{red line/.style = {line width=0.5ex, red, semitransparent}}
\draw[->] (0,1)--(1,0);
\draw[red line, ->] (1,1)--(0,0);
\begin{scope}[xshift = 2.5cm]
\draw[->, red line] (1,0)--(0,1);
\draw[->] (0,0)--(1,1);
\end{scope}					
\end{tikzpicture}	
\caption{\label{figure_avoiding} Forbidden fragments.}
\end{figure}
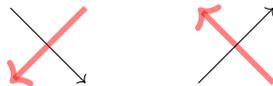
\vspace{-0.2cm}
\end{definition}
We denote by $\Paths({\bm i}, k)$ the set of rigorous paths in the oriented diagram $G(\bm i, k)$ for $k \in [m-1]$, and set
\[
\Paths(\bm i) \colonequals \bigsqcup_{k=1}^{m-1} \Paths({\bm i}, k).
\]
	
\begin{remark}\label{remark_avoiding}
A forbidden fragment in Figure~\ref{figure_avoiding} occurs when $\ell_i$ (red one) crosses over $\ell_j$ (black one) such that 
\begin{itemize}
	\item $i > j$, where the orientation of both wires is downward, or 
	\item $i < j$, where the orientation of both wires is upward.  
\end{itemize}	
\end{remark}

A node $t$ is called a {\em peak} of a rigorous path $P \in \Paths({\bm i})$ if $t$ is a local maximum of the path~$P$ with respect to the height of the diagram $G({\bm i})$.
Note that $P$ may have many peaks. 
We denote by $\Lambda(P)$ the set of peaks of $P$.
	
We use two expressions for a rigorous path $P$: a node-expression and a wire-expression.
Following the notation of \cite{GleizerPostnikov}, we express a rigorous path $P$ in $G({\bm i}, k)$ by 
\begin{equation}\label{equation_node_expression}
	P = (L_k \rightarrow a_{j_1} \rightarrow \cdots \rightarrow a_{j_s} \rightarrow L_{k+1}),
\end{equation}
where $a_{j_1},  \dots, a_{j_s}$ are the nodes at which the path $P$ crosses in order from one wire to another wire. 
We call \eqref{equation_node_expression} the {\em node-expression} of $P$ and write
\[
\node(P) \colonequals \{a_{j_1}, \dots, a_{j_s}\}.
\]
For instance, the path in the first diagram $G(\bm i,1)$ in Figure~\ref{figure_wd_oriented} is expressed as $L_1 \to a_4 \to a_5 \to L_2$. Similarly, the path in the third diagram $G(\bm i',3)$ in the same figure is expressed as $L_3 \to a_5 \to a_3 \to L_4$.
	
Also, a rigorous path can be expressed by recording the wires in the order of travel through. 
In other words, the rigorous path $P$ given in \eqref{equation_node_expression} can be written as
\begin{equation}\label{equation_wire_expression}
	\ell_{r_1} \rightarrow \cdots \rightarrow \ell_{r_{s+1}} \quad \quad \quad (r_1 = k, \quad r_{s+1} = k+1),
\end{equation}
where the node $a_{j_t}$ is at the intersection of $\ell_{r_t}$ and $\ell_{r_{t+1}}$ for each $t=1,\dots,s$.
The expression \eqref{equation_wire_expression} is called a {\em wire-expression}. 
For instance, the path in the first diagram in Figure~\ref{figure_wd_oriented} is expressed as $\ell_1 \to \ell_4 \to \ell_2$. 
Also, the third diagram in the same figure is expressed as $\ell_3 \to \ell_1 \to \ell_4$.
	
\subsection{String cone inequalities}\label{ssecStringInequalities}

We now introduce defining inequalities of the string cone $\mathcal{C}_{\bm i}^{(A_{m-1})}$ for~${\bm i} \in R(w_0^{(A_{m-1})})$.
	
\begin{definition}\label{definition_string_inequality}
Let ${\bm i} \in R(w_0^{(A_{m-1})})$ and 
$P$ a rigorous path in~$G({\bm i}, k)$ for some $k \in [m-1]$. 
The {\em string inequality associated with $P$} is defined by 
\begin{equation}\label{equation_stringcone}
\hat{\F}_P({\bm a}) \colonequals \sum_{j=1}^{\ell} A_j a_j \geq 0, \quad \quad \text{ where }A_j \colonequals \begin{cases}
	1 & \text{if $P$ travels from $\ell_r$ to $\ell_s$ at $a_j$ and $r < s$}, \\ 
	-1 & \text{if $P$ travels from $\ell_r$ to $\ell_s$ at $a_j$ and $r > s$}, \\ 
	 0 & \text{otherwise}.
\end{cases}
\end{equation}
Here, we regard ${\bm a} = (a_1,\dots,a_{\ell})$ as a coordinate system on $\R^\ell$.
\end{definition}

\begin{theorem}[{\cite[Corollary~5.8]{GleizerPostnikov}}]\label{theorem_GP_BZ}
Let ${\bm i} \in R(w_0^{(A_{m-1})})$.
The string cone $\mathcal{C}_{\bm i}^{(A_{m-1})}$ coincides with the set of ${\bm a} \in \R^{\ell}$ satisfying the string inequalities in Definition \ref{definition_string_inequality}, that is,
\begin{equation}\label{eq:explicit_description_string_cone_type_A}
\begin{aligned}
\mathcal C_{\bm i}^{(A_{m-1})} = \{{\bm a} \in \R^{\ell} \mid \hat{\F}_{P}({\bm a}) \geq 0 \quad \text{ for all } P \in \Paths(\bm i)\}.
\end{aligned}
\end{equation}
\end{theorem}

\begin{example}\label{example_string_cone_GC}
Let ${\bm i} = (1,2,1,3,2,1)$ and ${\bm i}' = (1,3,2,1,3,2)$ where the corresponding wiring diagrams are given in Figure \ref{figure_wiring_diagram_GC}.
\noindent 
Then the string cones $\mathcal{C}_{\bm i}^{(A_3)}$ and $\mathcal{C}_{\bm i'}^{(A_3)}$ are defined as follows:
\[
\begin{cases}
G({\bm i}, 1): a_1 \geq 0, a_2 - a_3 \geq 0, a_4 - a_5 \geq 0, \\
G({\bm i}, 2): a_3 \geq 0, a_5 - a_6 \geq 0, \\
G({\bm i}, 3): a_6 \geq 0, 
\end{cases}
\, \text{and} \,\,
\begin{cases}
G({\bm i'},1): a_1 \geq 0, a_3 - a_4 \geq 0, a_5 - a_6 \geq 0,\\
G({\bm i'},2): a_6 \geq 0, \\
G({\bm i'},3): a_2 \geq 0, a_3 - a_5 \geq 0, a_4 - a_6 \geq 0.
\end{cases}
\]
\end{example}

Notice that in Example~\ref{example_string_cone_GC}, the number of inequalities for $\mathcal{C}_{\bm i}^{(A_3)}$ is six, while that for $\mathcal{C}_{\bm i'}^{(A_3)}$ is seven. 
Depending on the choice of a reduced word, the number of inequalities for a string cone may vary.  

To study the non-redundancy of inequalities, \emph{chamber variables} $u_1,\dots,u_{\ell}$ are useful. 
For later use, we recall from~\cite{CKLP21_GC} a specific coordinate change $\Phi$ from $(a_1,\dots,a_{\ell})$ to $(u_1,\dots,u_{\ell})$.
The wiring diagram $G(\bm i)$ divides the plane into bounded or unbounded regions, and we call each region a \emph{chamber}. 
For each node $a_j$, we label a chamber having the top node $a_j$ by $\chamber_j$ (see Figure~\ref{figure_chamber}). 
We denote by $I_j$ the set of nodes contained in the boundary of $\chamber_j$. 
We divide the set $I_j$ into two sets $I_j^+$ and $I_j^-$ such that $I_j^+$ consists of nodes in the same column as $a_j$ and $I_j^- \coloneqq I_j \setminus I_j^+$. 
Define an $\mathbb{R}$-linear transformation $\Phi \colon \mathbb{R}^\ell \rightarrow \mathbb{R}^\ell$, $(a_1,\dots,a_{\ell}) \mapsto (u_1,\dots,u_{\ell})$, by 
\begin{equation}\label{eq_def_of_uj}
u_j \colonequals \sum_{a_k \in I_j^+} a_k - \sum_{a_k \in I_j^-} a_k
\end{equation}
for each $j \in [\ell]$. 
For each rigorous path $P$, denote by $\chamber(P)$ the region enclosed by $P$.
\begin{figure}[b]
\begin{tikzpicture}[scale = 0.8, yscale=0.6]
\def\nwires{4} 
\def\length{6} 

\foreach \i\y in {1/1,2/3,3/2,4/1,5/3,6/2}{
	\pgfmathsetmacro\ycoordinate{\length-\i+1}
	\pgfmathsetmacro\yplusone{\y+1}
	\foreach \x in {1,...,\nwires}{
		\ifthenelse{\x = \y}
			{\coordinate (\i SW) at (\x, \ycoordinate);
			\coordinate (\i SE) at (\x+1, \ycoordinate);
			\coordinate (\i NW) at (\x, \ycoordinate+1);
			\coordinate (\i NE) at (\x+1,\ycoordinate+1);

			\coordinate (\i C) at (\x+0.5, \ycoordinate + 0.5);

			\draw (\x,\ycoordinate) -- (\i SW)-- (\i NE) -- (\x+1, \ycoordinate+1);
			\draw (\x+1,\ycoordinate) --(\i SE) --(\i NW) --(\x, \ycoordinate+1);
			}{}
		\ifthenelse{\x=\y \OR \x = \yplusone}{}{\draw (\x,\ycoordinate)--(\x,\ycoordinate+1);}

}}

\foreach \x in{1,...,\nwires}{
\draw (\x,\length+1)--(\x,\length+1.5);
\draw (\x,1)--(\x,-0.5);
}

\filldraw[draw=none,pattern=north east lines,pattern color=black!30] (1,-0.5)--(2,-0.5)--(6SW)--(6C)--(6NW)--(4SE)--(4C)--(4SW)--cycle;

\filldraw[draw=none, pattern=north west lines, pattern color = blue!30]
(2,-0.5)--(3,-0.5)--(6SE)--(6C)--(6SW)--cycle;

\filldraw[draw=none, pattern= dots, pattern color = purple!20]
(3,-0.5)--(4,-0.5)--(5SE)--(5C)--(6NE)--(6C)--(6SE)--cycle;

\filldraw[draw=none, pattern=fivepointed stars, pattern color=magenta!20] 
(4C)--(3C)--(3NW)--(1SE)--(1C)--(1SW)--(4NW)--cycle;

\filldraw[draw=none, pattern=horizontal lines, pattern color = green!30] (6C)--(5C)--(5NW)--(3SE)--(3C)--(4C)--(4SE)--(6NW)--cycle;

\filldraw[draw=none, pattern = vertical lines, pattern color = black!30] 
(5C)--(5NE)--(2SE)--(2C)--(3C)--(3SE)--(5NW)--cycle;

\foreach \x in {1,...,6}{
\node[above] at (\x C) {$a_{\x}$};

\node[color=red, above =-1cm of \x C]
 {$\mathscr{C}_{\x}$};

}
\end{tikzpicture}
\caption{Chambers $\mathscr C_j$ for $\bm i = (1,3,2,1,3,2)$.}\label{figure_chamber}
\end{figure}
Then we have the following lemma.

\begin{lemma}[{\cite[Section 4.1]{CKLP21_GC}}]\label{lemma_chamber_varialbe}
The map $\Phi$ is a unimodular transformation. 
In particular, by regarding ${\bm u} = (u_1,\dots,u_{\ell})$ as a coordinate system on $\mathbb{R}^\ell$, the coordinate change $(a_1,\dots,a_{\ell}) \mapsto (u_1,\dots,u_{\ell})$ is a change of bases of a $\mathbb{Z}$-lattice in $(\mathbb{R}^\ell)^\ast$. 
Moreover, it holds that
\[
\hat{\F}_P ({\bm a}) = \sum_{\chamber_j \subseteq \chamber(P)} \Phi^\ast(u_j),
\]
where $\Phi^\ast \colon (\mathbb{R}^\ell)^\ast \rightarrow (\mathbb{R}^\ell)^\ast$ denotes the dual map, and $\Phi^\ast(u_j)$ is given by the formula \eqref{eq_def_of_uj}.
\end{lemma}

We call $u_j$ the \emph{$j$th chamber variable} for $j \in [\ell]$.

\begin{example}
For $\bm i = (1,3,2,1,3,2)$, consider a rigorous path $P = (\ell_3 \to \ell_1 \to \ell_4)$. This path encloses chambers $\chamber_3$ and $\chamber_4$ (see Figure~\ref{figure_chamber}). Moreover, we have $\Phi^\ast(u_3) = a_3 + a_6 - (a_4 + a_5)$ and $\Phi^\ast(u_4) = a_4 - a_6$. 
Accordingly, it holds that
\[
\Phi^\ast(u_3) + \Phi^\ast(u_4) = \left(a_3 + a_6 - (a_4 + a_5)\right) + (a_4 - a_6) 
= a_3 - a_5 = \hat{\F}_P ({\bm a}).
\]
\end{example}

\begin{proposition}[{see \cite[Proposition~4.5]{CKLP21_GC}}]\label{prop_non_redundancy_type_A}
Let ${\bm i} \in R(w_0^{(A_{m-1})})$. 
Then the expression~\eqref{eq:explicit_description_string_cone_type_A} of $\mathcal{C}_{\bm i}^{(A_{m-1})}$ is non-redundant, and the number of facets of the string cone~$\mathcal{C}_{\bm i}^{(A_{m-1})}$ is 
\[
\# \Paths(\bm i)
= \sum_{k=1}^{m-1} \# \Paths(\bm i,k).
\]
\end{proposition}

Now we define a {\em string polytope} in type $A$.

\begin{definition}[{\cite{Littelmann98}}]\label{def_string_polytope}
Let $\bm i = (i_1,\dots,i_{\ell}) \in R(w_0^{(A_{m-1})})$. 
For a dominant integral weight $\lambda$, the string polytope~$\Delta_{\bm i}^{(A)}(\lambda)$ is the rational convex polytope defined as the intersection of the string cone $\mathcal C_{\bm i}^{(A_{m-1})}$ and the cone $\mathcal C^{\lambda}_{\bm i}$ (called a {\em $\lambda$-cone}) defined by the following inequalities:
\[
\begin{split}
a_1 &\leq \langle \lambda - a_{\ell} \alpha_{i_{\ell}} - \cdots - a_2 \alpha_{i_2}, h_{i_1} \rangle, \\
\vdots & \\
a_{\ell-1} & \leq \langle \lambda- a_{\ell} \alpha_{i_{\ell}}, h_{i_{\ell-1}} \rangle,\\
a_{\ell} &\leq \langle \lambda, h_{i_{\ell}} \rangle. 
\end{split}
\]
Here, $\{\alpha_i \mid i \in [m-1]\}$ is the set of simple roots, and $\{h_i \mid i \in [m-1]\}$ the set of simple coroots. Moreover, $\langle \cdot, \cdot \rangle$ denotes the canonical pairing.
\end{definition}

Let $\{\varpi_i \mid i \in [m-1]\}$ be the set of fundamental weights. 
We note that for Lie type $A$, we have 
	\[
	\begin{split}
	\langle \varpi_i, h_j \rangle &= \delta_{i,j}; \\
	\langle \alpha_i, h_j \rangle &= \begin{cases}
		2 & \text{ if } i = j; \\
		-1 & \text{ if } |i-j| = 1;  \\
		0 & \text{ otherwise}.
	\end{cases}
	\end{split}
	\]
\begin{example}
Let $\bm i = (1,3,2,1,3,2)$ and $\lambda = \lambda_1 \varpi_1 + \lambda_2 \varpi_2 + \lambda_3 \varpi_3$ a dominant integral weight. 
We have $\langle \lambda, h_{i_6} \rangle = 
\langle \lambda, h_2 \rangle = \lambda_2$. Similarly, we obtain 
\[
\langle \lambda - a_{6} \alpha_{i_6}, h_{i_5} \rangle 
= \langle \lambda, h_3 \rangle - a_6 \langle \alpha_2, h_3 \rangle 
= \lambda_3 + a_6. 
\]
Accordingly, the $\lambda$-cone $\mathcal C^{\lambda}_{\bm i}$ is given by 
\[
\begin{split}
a_1 &\leq \lambda_1 + a_3 - 2 a_4 + a_6, \\
a_2 &\leq \lambda_3 + a_3 - 2 a_5 + a_6, \\
a_3 &\leq \lambda_2 + a_4 + a_5 -2a_6, \\
a_4 &\leq \lambda_1 + a_6, \\
a_5 &\leq \lambda_3 + a_6, \\
a_6 &\leq \lambda_2.
\end{split}
\]

\end{example}

\section{Folding procedure for string cones}
\label{sec_folding}
In this section, we review relations among string cones in types $A_{2n-1}$, $B_n$, and $C_n$, following~\cite{Kas96, NS05, Fujita18}.
Let $m = 2n$ and consider the special linear Lie algebra $\mathfrak{g} \coloneqq \mathfrak{sl}_{2n}(\mathbb{C})$. 
Recall that the set $I$ of vertices of the Dynkin diagram of $\mathfrak{g}$ is identified with $[2n-1]$. 
We set $\overline{i} = 2n-i$ for $1 \leq i \leq 2n-1$, and define a bijection $\omega \colon I \rightarrow I$ by $i \mapsto \overline{i}$.
Denoting the Cartan matrix of $\mathfrak{g}$ by $C = (c_{i, j})_{i, j \in I}$, it holds that $c_{\omega(i), \omega(j)} = c_{i, j}$ for all $i, j \in I$, which implies that $\omega$ corresponds to an automorphism of the Dynkin diagram. 
Let us write $\breve{I} \coloneqq \{1, 2, \ldots, n\} \subseteq I$, which is a complete set of representatives for the $\omega$-orbits in $I$. 
We set $m_i \coloneqq \min\{k \in \mathbb{Z}_{>0} \mid \omega^k(i) = i\}$ for $i \in \breve{I}$, that is, 
\begin{equation}\label{equation_m}
	m_1 = m_2 = \cdots = m_{n-1} = 2 \quad \text{and} \quad m_n = 1.  
\end{equation}
Define an integer matrix $\breve{C} \coloneqq (\breve{c}_{i, j})_{i, j \in \breve{I}}$ by 
\[\breve{c}_{i, j} \coloneqq \sum_{0 \le k < m_j} c_{i, \omega^k(j)}\]
for $i, j \in \breve{I}$. 
Let $\breve{\mathfrak{g}}$ denote the finite-dimensional semisimple Lie algebra over $\mathbb{C}$ with Cartan matrix~$\breve{C}$, which is called the \emph{orbit Lie algebra} associated with $\omega$.
Since $\breve{C}$ is the indecomposable Cartan matrix of type $B_n$, the orbit Lie algebra $\breve{\mathfrak{g}}$ is isomorphic to $\mathfrak{so}_{2n+1}(\mathbb{C})$, where $\breve{I}$ is identified with the set of vertices of the Dynkin diagram of type $B_n$ as follows:
\begin{align*}
&B_n\ \begin{xy}
\ar@{-} (50,0) *++!D{1} *\cir<3pt>{};
(60,0) *++!D{2} *\cir<3pt>{}="C"
\ar@{-} "C";(65,0) \ar@{.} (65,0);(70,0)^*!U{}
\ar@{-} (70,0);(75,0) *++!D{n-1} *\cir<3pt>{}="D"
\ar@{=>} "D";(85,0) *++!D{n} *\cir<3pt>{}="E"
\end{xy}.
\end{align*}
Let $W^{B_n}$ denote the Weyl group of $\breve{\mathfrak{g}}$, $w_0^{(B_n)} \in W^{B_n}$ the longest element, and $R(w_0^{(B_n)})$ the set of reduced words for $w_0^{(B_n)}$. 
We write $N \coloneqq n^2$, which is the length of $w_0^{(B_n)}$. 
The Weyl group~$W^{B_n}$ is naturally regarded as a subgroup of the Weyl group $\mathfrak{S}_{2n}$ of type $A_{2n-1}$ (see, for instance, \cite[Section 3]{FRS97}).
We denote by 
\begin{equation}\label{equation_Theta_inclusion}
\Theta \colon W^{B_n} \hookrightarrow \mathfrak{S}_{2n}
\end{equation} 
the inclusion map. 
Then a reduced word ${\bm i} = (i_1, i_2, \ldots, i_N) \in R(w_0^{(B_n)})$ for $w_0^{(B_n)}$ induces a reduced word $\hat{\bm i} \coloneqq (i_{1, 1}, \ldots, i_{1, m_{i_1}}, \ldots, i_{N, 1}, \ldots, i_{N, m_{i_N}})$ for $\Theta(w_0^{(B_n)}) = w_0^{(A_{2n-1})}$, where 
\[(i_{k, 1}, \ldots, i_{k, m_{i_k}}) \coloneqq 
\begin{cases}
(i_k, \overline{i_k})\quad &(i_k = 1, 2, \ldots, n-1),\\
(n)\quad &(i_k = n).
\end{cases}
\]
Here, we notice that $m_{i_k} = 2$ if $i_k = 1,2,\dots,n-1$ and $m_{i_k} = 1$ for $i_k = n$ as in \eqref{equation_m}.
We call $\hat{\bm i}$ the \emph{lift} of ${\bm i}$. 
Let $\mathcal{C}_{\bm i}^{(B_n)} \subseteq \mathbb{R}^N$ denote the string cone associated with ${\bm i} = (i_1, i_2, \ldots, i_N) \in R(w_0^{(B_n)})$.
We define an injective $\mathbb{R}$-linear map $\Upsilon_{\bm i}^{B, A} \colon \mathbb{R}^N \hookrightarrow \mathbb{R}^{m_{i_1} + \cdots + m_{i_N}}$ by 
\[\Upsilon_{\bm i}^{B, A} (a_1, \ldots, a_N) \coloneqq (\underbrace{a_1, \ldots, a_1}_{m_{i_1}}, \ldots, \underbrace{a_N, \ldots, a_N}_{m_{i_N}})\]
for $(a_1, \ldots, a_N) \in \mathbb{R}^N$. 
Then the string cone $\mathcal{C}_{\bm i}^{(B_n)}$ is identified with a slice of the string cone~$\mathcal C_{\hat{\bm i}}^{(A_{2n-1})}$ as follows. 

\begin{theorem}[{see \cite[Theorem 1]{NS05} and \cite[Corollary 4.8]{Fujita18}}]\label{t:relation_between_A_and_B}
Let ${\bm i} = (i_1, i_2, \ldots, i_N) \in R(w_0^{(B_n)})$. 
Then the following equality holds: 
\[\Upsilon_{\bm i}^{B, A} (\mathcal{C}_{\bm i}^{(B_n)}) = \{(a_{k, \ell})_{1 \leq k \leq N, 1 \leq \ell \leq m_{i_k}} \in \mathcal C_{\hat{\bm i}}^{(A_{2n-1})} \mid a_{k, 1} = a_{k, 2} = \cdots = a_{k, m_{i_k}},\ 1 \leq k \leq N\}.\]
\end{theorem}

Define a Lie algebra automorphism $\hat{\omega} \colon \mathfrak{g} \xrightarrow{\sim} \mathfrak{g}$ of $\mathfrak{g} = \mathfrak{sl}_{2n}(\mathbb{C})$ by 
\[\hat{\omega}(X) \coloneqq (\overline{w}_0)^{-1} \cdot (-X^T) \cdot \overline{w}_0\] 
for $X \in \mathfrak{g}$, where $X^T$ denotes the transpose of $X$, and $\overline{w}_0$ is an integer $2n \times 2n$ matrix given by
\[\overline{w}_0 \coloneqq \begin{pmatrix}
0 & 0 & 0 & \cdots & -1 \\
\vdots & \vdots & \vdots & \iddots & \vdots \\
0 & 0 & 1 & \cdots & 0 \\
0 & -1 & 0 & \cdots & 0 \\
1 & 0 & 0 & \cdots & 0
\end{pmatrix}.\]
Then the fixed point Lie subalgebra 
\[\mathfrak{g}^{\hat{\omega}} \coloneqq \{X \in \mathfrak{g} \mid \hat{\omega}(X) = X\}\]
of $\mathfrak{g}$ coincides with the symplectic Lie algebra  
\[\mathfrak{sp}_{2n}(\mathbb{C}) \coloneqq \{X \in \mathfrak{g} \mid X^T \overline{w}_0 + \overline{w}_0 X = 0\}\]
with respect to the skew-symmetric matrix $\overline{w}_0$. 
This is the simple Lie algebra of type $C_n$. 
For $1 \leq i, j \leq 2n$, we denote by $E_{i, j}$ the $2n \times 2n$ matrix whose $(i, j)$-entry is $1$ and other entries are all $0$. 
The automorphism $\hat{\omega}$ coincides with the Lie algebra automorphism $\mathfrak{g} \xrightarrow{\sim} \mathfrak{g}$ induced by $\omega$, that is, 
\[\hat{\omega}(e_i) = e_{\omega(i)},\ \hat{\omega}(f_i) = f_{\omega(i)},\ \hat{\omega}(h_i) = h_{\omega(i)}\] 
for all $i \in I$, where $e_i, f_i, h_i \in \mathfrak{g}$, $i \in I$, are Chevalley generators of $\mathfrak{g}$ given by $e_i \coloneqq E_{i, i+1}$, $f_i \coloneqq E_{i+1, i}$, and $h_i \coloneqq E_{i, i} - E_{i+1, i+1}$. 
We identify $\breve{I}$ with the set of vertices of the Dynkin diagram of type $C_n$ as follows:
\begin{align*}
&C_n\ \begin{xy}
\ar@{-} (50,0) *++!D{1} *\cir<3pt>{};
(60,0) *++!D{2} *\cir<3pt>{}="C"
\ar@{-} "C";(65,0) \ar@{.} (65,0);(70,0)^*!U{}
\ar@{-} (70,0);(75,0) *++!D{n-1} *\cir<3pt>{}="D"
\ar@{<=} "D";(85,0) *++!D{n} *\cir<3pt>{}="E"
\end{xy}.
\end{align*}
Then the Cartan matrix $C^\prime = (c^\prime_{i, j})_{i, j \in \breve{I}}$ of type $C_n$ coincides with the transpose of $\breve{C}$. 
In other words, the orbit Lie algebra $\breve{\mathfrak{g}}$ is the Langlands dual Lie algebra of the fixed point Lie subalgebra~$\mathfrak{g}^{\hat{\omega}}$.  
Summarizing, we obtain the following diagram: 
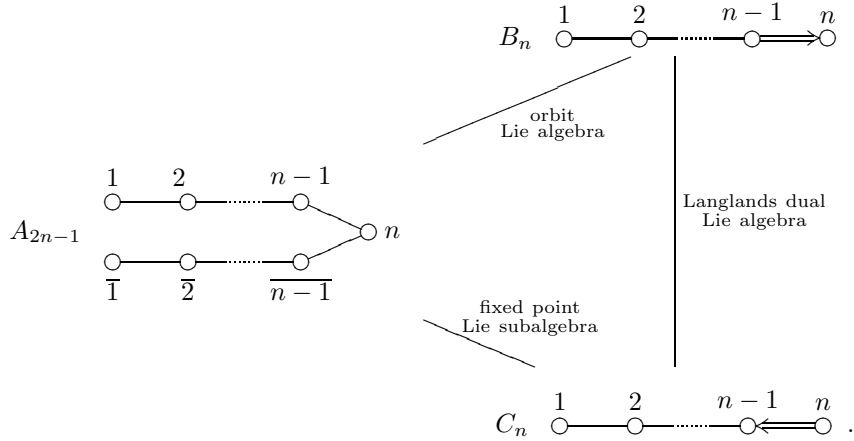
\begin{align*}
\xymatrix{                   
      & \ar@{-}[ld]^-{\substack{{\rm orbit}\\{\rm Lie\ algebra}}} B_n\ \begin{xy}
\ar@{-} (50,0) *++!D{1} *\cir<3pt>{};
(60,0) *++!D{2} *\cir<3pt>{}="C"
\ar@{-} "C";(65,0) \ar@{.} (65,0);(70,0)^*!U{}
\ar@{-} (70,0);(75,0) *++!D{n-1} *\cir<3pt>{}="D"
\ar@{=>} "D";(85,0) *++!D{n} *\cir<3pt>{}="E"
\end{xy} \ar@{-}[dd]^-{\substack{{\rm Langlands\ dual}\\{\rm Lie\ algebra}}}   \\
A_{2n-1}\ \begin{xy}
\ar@{-} (20,4) *++!D{1} *\cir<3pt>{};
(30,4) *++!D!R(0.4){2} *\cir<3pt>{}="B",
\ar@{-} "B";(35,4) \ar@{.} (35,4);(40,4)^*!U{}
\ar@{-} (40,4);(45,4) *++!D{n-1} *\cir<3pt>{}="C"
\ar@{-} "C";(54,0) *++!L{n} *\cir<3pt>{}="D"
\ar@{-} "D";(45,-4) *++!U{\overline{n-1}} *\cir<3pt>{}="E"
\ar@{-} "E";(40,-4) \ar@{.} (40,-4);(35,-4)^*!U{}
\ar@{-} (35,-4);(30,-4) *++!U{\overline{2}} *\cir<3pt>{}="F"
\ar@{-} "F";(20,-4) *++!U{\overline{1}} *\cir<3pt>{}
\end{xy} & \\
                         & \ar@{-}[lu]_-{\substack{\quad{\rm fixed\ point}\\\quad{\rm Lie\ subalgebra}}} C_n\ \begin{xy}
\ar@{-} (50,0) *++!D{1} *\cir<3pt>{};
(60,0) *++!D{2} *\cir<3pt>{}="C"
\ar@{-} "C";(65,0) \ar@{.} (65,0);(70,0)^*!U{}
\ar@{-} (70,0);(75,0) *++!D{n-1} *\cir<3pt>{}="D"
\ar@{<=} "D";(85,0) *++!D{n} *\cir<3pt>{}="E"
\end{xy}.
} 
\end{align*}
Note that the Weyl group $W^{C_n}$ of type $C_n$ is naturally isomorphic to $W^{B_n}$.
Under the isomorphism, the longest element $w_0^{(C_n)}$ of $W^{C_n}$ corresponds to $w_0^{(B_n)} \in W^{B_n}$, and the set $R(w_0^{(B_n)})$ coincides with the set $R(w_0^{(C_n)})$ of reduced words for $w_0^{(C_n)}$ as a subset of $\breve{I}^N$.
We write $m_1^\prime = m_2^\prime = \cdots = m_{n-1}^\prime = 1$ and $m_n^\prime = 2$. 
For ${\bm i} = (i_1, i_2, \ldots, i_N) \in R(w_0^{(B_n)}) = R(w_0^{(C_n)})$, define $\mathbb{R}$-linear automorphisms $\Gamma_{\bm i}^{B, C} \colon \mathbb{R}^N \rightarrow \mathbb{R}^N$ and $\Gamma_{\bm i}^{C, B} \colon \mathbb{R}^N \rightarrow \mathbb{R}^N$ by 
\begin{align*}
&\Gamma_{\bm i}^{B, C} (a_1, \ldots, a_N) \coloneqq (m_{i_1} a_1, \ldots, m_{i_N} a_N),\\
&\Gamma_{\bm i}^{C, B} (a_1, \ldots, a_N) \coloneqq (m_{i_1}^\prime a_1, \ldots, m_{i_N}^\prime a_N)
\end{align*}
for $(a_1, \ldots, a_N) \in \mathbb{R}^N$, respectively. 
Then we have $\Gamma_{\bm i}^{B, C} \circ \Gamma_{\bm i}^{C, B} = \Gamma_{\bm i}^{C, B} \circ \Gamma_{\bm i}^{B, C} = 2 \cdot {\rm id}_{\mathbb{R}^N}$.
Kashiwara~\cite[Section~5]{Kas96} gave a similarity of crystal bases between types $B_n$ and $C_n$. 
The following theorem is a straightforward consequence of this similarity. 

\begin{theorem}[{see \cite[Section~5]{Kas96} and \cite[Proposition~6.3]{Fujita18}}]\label{t:similarity_between_B_and_C}
Let ${\bm i} \in R(w_0^{(B_n)}) = R(w_0^{(C_n)})$. 
Then the following equalities hold: 
\begin{align*}
&\Gamma_{\bm i}^{B, C} (\mathcal{C}_{\bm i}^{(B_n)}) = \mathcal{C}_{\bm i}^{(C_n)}, & &\Gamma_{\bm i}^{C, B} (\mathcal{C}_{\bm i}^{(C_n)}) = \mathcal{C}_{\bm i}^{(B_n)}.
\end{align*}
\end{theorem}

For ${\bm i} = (i_1, i_2, \ldots, i_N) \in R(w_0^{(C_n)})$, we define a surjective $\mathbb{R}$-linear map $\Omega_{\bm i}^{A, C} \colon \mathbb{R}^{m_{i_1} + \cdots + m_{i_N}} \twoheadrightarrow \mathbb{R}^N$ by 
\begin{equation}\label{equation_Omega_AC}
\Omega_{\bm i}^{A, C} (a_{1, 1}, \ldots, a_{1, m_{i_1}}, \ldots, a_{N, 1}, \ldots, a_{N, m_{i_N}}) \coloneqq (a_{1, 1} + \cdots + a_{1, m_{i_1}}, \ldots, a_{N, 1} + \cdots + a_{N, m_{i_N}})
\end{equation}
for $(a_{1, 1}, \ldots, a_{1, m_{i_1}}, \ldots, a_{N, 1}, \ldots, a_{N, m_{i_N}}) \in \mathbb{R}^{m_{i_1} + \cdots + m_{i_N}}$.
Then we have $\Gamma_{\bm i}^{B, C} = \Omega_{\bm i}^{A, C} \circ \Upsilon_{\bm i}^{B, A}$.
Combining Theorem \ref{t:relation_between_A_and_B} with Theorem \ref{t:similarity_between_B_and_C}, it follows that 
\begin{align*}
\mathcal{C}_{\bm i}^{(C_n)} &= \Gamma_{\bm i}^{B, C} (\mathcal{C}_{\bm i}^{(B_n)})\\
&= \Omega_{\bm i}^{A, C} \circ \Upsilon_{\bm i}^{B, A} (\mathcal{C}_{\bm i}^{(B_n)})\\ 
&\subseteq \Omega_{\bm i}^{A, C} (\mathcal C_{\hat{\bm i}}^{(A_{2n-1})}). 
\end{align*}
More strongly, we know the following. 

\begin{theorem}[{see \cite[Theorem 5.7]{Fujita18}}]\label{t:relation_between_A_and_C}
Let ${\bm i} \in R(w_0^{(C_n)})$. 
Then the following equality holds: 
\begin{align*}
&\Omega_{\bm i}^{A, C} (\mathcal C_{\hat{\bm i}}^{(A_{2n-1})}) = \mathcal{C}_{\bm i}^{(C_n)}.
\end{align*}
\end{theorem}

\begin{definition}\label{def_string_inequalities_BC}
Let ${\bm i} \in R(w_0^{(B_n)}) = R(w_0^{(C_n)})$ and $P$ a rigorous path in $G(\hat{\bm i}, k)$ for some $k \in [2n-1]$. 
Then the \emph{string inequality in type $B_n$ associated with $P$} is defined by 
\begin{equation}\label{equation_string_cone_type_B}
\hat{\F}_P^{(B)}({\bm a}) \coloneqq \hat{\F}_P (\Upsilon_{\bm i}^{B, A} ({\bm a})) \geq 0
\end{equation}
for ${\bm a} \in \mathbb{R}^N$. 
In addition, the \emph{string inequality in type $C_n$ associated with $P$} is defined as 
\begin{equation}\label{equation_string_cone_type_C}
\hat{\F}_P^{(C)}({\bm a}) \coloneqq \hat{\F}_P^{(B)} (\Gamma_{\bm i}^{C, B} ({\bm a})) \geq 0
\end{equation}
for ${\bm a} \in \mathbb{R}^N$. 
\end{definition}

The following theorem is a straightforward consequence of Theorems \ref{t:relation_between_A_and_B} and \ref{t:similarity_between_B_and_C}. 

\begin{theorem}\label{thm_string_cones_BC}
For ${\bm i} \in R(w_0^{(B_n)}) = R(w_0^{(C_n)})$, the following equalities hold: 
\begin{align*}
&\mathcal{C}_{\bm i}^{(B_n)} = \{{\bm a} \in \R^N \mid \hat{\F}_P^{(B)}({\bm a}) \geq 0 \quad \text{ for all } P \in \Paths(\hat{\bm i})\},\\
&\mathcal{C}_{\bm i}^{(C_n)} = \{{\bm a} \in \R^N \mid \hat{\F}_P^{(C)}({\bm a}) \geq 0 \quad \text{ for all } P \in \Paths(\hat{\bm i})\}.
\end{align*}
\end{theorem}

\section{Folding procedure for rigorous paths}
\label{sec_symmetric_GP_paths}
\begin{figure}[hbt]
\begin{tikzpicture}[scale = 0.8, yscale=0.6]
\def\nwires{3} 
\def\length{9} 

\foreach \i\y in {1/1,2/2,3/3,4/1,5/2,6/3,7/1,8/2,9/3}{
	\pgfmathsetmacro\ycoordinate{\length-\i+1}
	\pgfmathsetmacro\yplusone{\y+1}
	\foreach \x in {1,...,\nwires}{
		\ifthenelse{\x = \y}
			{\coordinate (\i LSW) at (-\nwires+\x-0.5, \ycoordinate+0);
			\coordinate (\i LSE) at (-\nwires+\x +0.5, \ycoordinate+0);
			\coordinate (\i LNW) at (-\nwires+\x-0.5, \ycoordinate+1);
			\coordinate (\i LNE) at (-\nwires+\x +0.5,\ycoordinate+1);

			\coordinate (\i RSW) at (\nwires - \x -0.5, \ycoordinate );
			\coordinate (\i RSE) at (\nwires - \x + 0.5, \ycoordinate );
			\coordinate (\i RNW) at (\nwires - \x - 0.5, \ycoordinate +1);
			\coordinate (\i RNE) at (\nwires - \x + 0.5, \ycoordinate +1);

			\coordinate (\i LC) at (-\nwires + \x, \ycoordinate+ 0.5);
			\coordinate (\i RC) at (\nwires - \x, \ycoordinate+0.5);
		
			\draw (-\nwires+\x-0.5,\ycoordinate) -- (\i LSW)-- (\i LNE) -- (-\nwires+\x+0.5, \ycoordinate+1);
			\draw (-\nwires+\x+0.5,\ycoordinate) --(\i LSE) --(\i LNW) --(-\nwires+\x-0.5, \ycoordinate+1);

			\draw (\nwires - \x -0.5, \ycoordinate) -- (\i RSW)-- (\i RNE) -- (\nwires - \x + 0.5, \ycoordinate + 1);
			\draw (\nwires - \x +0.5, \ycoordinate) -- (\i RSE)-- (\i RNW) -- (\nwires - \x -0.5, \ycoordinate + 1);
			}{}

		\ifthenelse{\x=\y \OR \x = \yplusone}{}{
			\draw
			(-\nwires+\x-0.5,\ycoordinate)--(-\nwires+\x-0.5,\ycoordinate+1);
			\draw
			(\nwires-\x+0.5,\ycoordinate)--(\nwires-\x+0.5,\ycoordinate+1);
		}}

	\ifthenelse{\y=\nwires}{
		\node[above] at (\i LC) {$t_{\i}$};}
	{\node[above] at (\i LC) {$\bar{t}_{\i}$};
	\node[above] at (\i RC) {$t_{\i}$};}

	\node at (-\nwires,\ycoordinate+0.5) {$\y$};
	\foreach \x in {1,...,\nwires}{
		\node at (-\nwires+\x-0.5, \length+2.5) {$U_{\x}$};
		\node at (\nwires-\x+0.5, \length+2.5) {$U_{\bar{\x}}$};
		\node at (-\nwires + \x - 0.5, -0.5) {$L_{\bar{\x}}$};
		\node at (\nwires - \x + 0.5, -0.5) {$L_{\x}$};
}}

\foreach \x in{1,...,\nwires}{
	\draw (\x-0.5,\length+1)--(\x-0.5,\length+1.5);
	\draw (\x-0.5,1)--(\x-0.5,0.5);

	\draw (-\x+0.5,\length+1)--(-\x+0.5,\length+1.5);
	\draw (-\x+0.5,1)--(-\x+0.5,0.5);

}

\draw[teal!50!green, dashed] (0,0.5)--(0,\length+1.5);

\end{tikzpicture}
\caption{Symplectic wiring diagram for $\bm i = (1,2,3,1,2,3,1,2,3)$.}\label{figure_symplectic_wiring_diagram_123123123}
\end{figure}
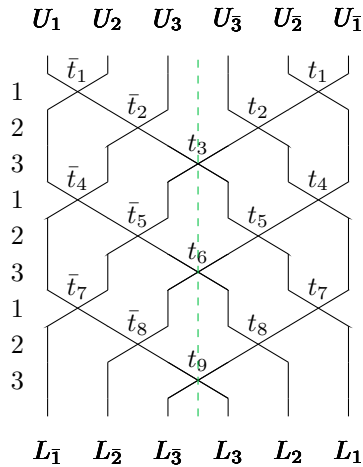

In this section, we introduce \emph{symplectic} wiring diagrams to provide a combinatorial way of describing string inequalities in types $B_n$ and $C_n$. 
Let $\bm i$ be a reduced word for $w_0^{(C_n)}$.
As before, we draw a symplectic wiring diagram to represent $\bm i$. 
As presented in Figure~\ref{figure_symplectic_wiring_diagram_123123123}, the symplectic wiring diagram $G^{\symp}({\bm i})$ consists of $2n$ wires labeled by $\ell_{1},\dots,\ell_{n},\ell_{\overline{n}},\dots,\ell_{\bar{1}}$. 
For each $k = 1,\dots,n,\bar{n},\dots,\bar{1}$, the upper and lower ends of the line $\ell_k$ are labeled by $U_k$ and $L_k$, respectively. 
Moreover, the $j$th crossing (from the top) is located on the $i_j$th and $(2n-i_j)$th column of $G^{\symp}(\bm i)$. 
We call each crossing a \emph{node} as before. 
If $i_j \neq n$, then we label the $j$th crossing (from the top) on the $i_j$th column with $\bar{t}_j$, and the $j$th crossing on the $(2n-i_j)$th column with $t_j$. If $i_j = n$, we label the $j$th crossing on the $n$th column with $t_j$.
This can be regarded as the wiring diagram of $\hat{\bm i}$ in type $A_{2n-1}$ via the identification $\bar{k} = 2n+1-k$ for $1 \le k \le n$. We denote by ${\bm t}$ the coordinate system on $\R^{n(2n-1)}$ in which the string cone 
$\mathcal{C}_{\bm i}^{(A_{2n-1})}$ lives where the components of ${\bm t}$ are {\em pairwisely ordered}, for instance, 
\[
	{\bm t} = (t_1, \bar{t}_1, t_2, \bar{t}_2, t_3, t_4, \bar{t}_4, \dots, t_8, \bar{t}_8, t_9)
\]
in Figure \ref{figure_symplectic_wiring_diagram_123123123}.

Note that the wiring diagram $G^{\symp}(\bm i)$ is \emph{symmetric} with respect to the central vertical line, called the \emph{wall}, between $L_{\bar{n}}$ and~$L_{n}$. For instance, in Figure~\ref{figure_symplectic_wiring_diagram_123123123}, 
we represent the wall using the green-dashed line.

For $k \in [n]$, we denote by $G^{\symp}(\bm i,k)$ the oriented symplectic wiring diagram such that the wires $\ell_1,\dots,\ell_k$ are oriented upward and the other wires $\ell_{k+1},\dots,\ell_n, \ell_{\bar{n}},\dots,\ell_{\bar{1}}$ are oriented downward. 
Similarly, let $G^{\symp}(\bm i,\bar{k})$ be the oriented symplectic wiring diagram such that the wires $\ell_1,\dots,\ell_n,\ell_{\bar{n}},\dots,\ell_{\bar{k}}$ are oriented upward and the other wires $\ell_{\overline{k-1}},\dots,\ell_{\bar{1}}$ are oriented downward.
Denote by $\Paths^{\symp}(\bm i,k)$ the set of rigorous paths in the oriented diagram $G^{\symp}(\bm i,k)$ that are defined in a way similar to Section \ref{ssecWiringDiagrams}.

For each $P \in \Paths^{\symp}(\bm i,k)$, we define a linear function $\hat{\F}_P(\bm t)$ on $\R^{n(2n-1)}$ in an exactly the same way as in~\eqref{equation_stringcone}. 
In addition, we use a wire-expression for $P \in G^{\symp}(\bm i,k)$ as in Section~\ref{ssecWiringDiagrams}, and denote by $\Lambda(P)$ the set of peaks of $P$.
There is a natural bijective correspondence between $\Paths^{\symp}(\bm i, k)$ and $\Paths(\hat{\bm i},k)$. 
For $P \in \Paths^{\symp}(\bm i,k)$, we denote by $\hat{P} \in \Paths(\hat{\bm i},k)$ the corresponding path, called the \emph{lift} of $P$. 
We note that a path $P$ and its lift $\hat{P}$ look same but have different labelings.
Set 
\[
\Paths^{\symp}({\bm i}) \colonequals \bigsqcup_{k=1}^{n} \Paths^{\symp}({\bm i}, k),
\]
which corresponds to $\Paths(\hat{\bm i})$.

\begin{definition}
Let $\bm i \in \RCn$, and $k \in [n]$.
For a rigorous path $P = (\ell_{r_1} \to \cdots \to \ell_{r_{s+1}}) \in \Paths^{\symp}({\bm i},k)$, its \emph{mirror}~$P^{\vee}$ is defined by
\[
P^{\vee} \colonequals (\ell_{\overline{r_{s+1}}} \to \cdots \to \ell_{\overline{r_1}}) \in \Paths^{\symp}({\bm i},\overline{k+1}),
\]
where $\overline{n+1} \coloneqq n$ and $\bar{\bar{r}} \coloneqq r$ for $r \in [n]$.
A path $P \in \Paths^{\symp}(\bm i,n)$ is said to be \emph{symmetric} if $P = P^{\vee}$.
\end{definition}

\begin{example}\label{example_symmetric_paths}
Let $\bm i = (2,1,2,1) \in R(w_0^{(C_2)})$. There are five rigorous paths in $\Paths^{\symp}(\bm i, 2)$ as follows.

\smallskip
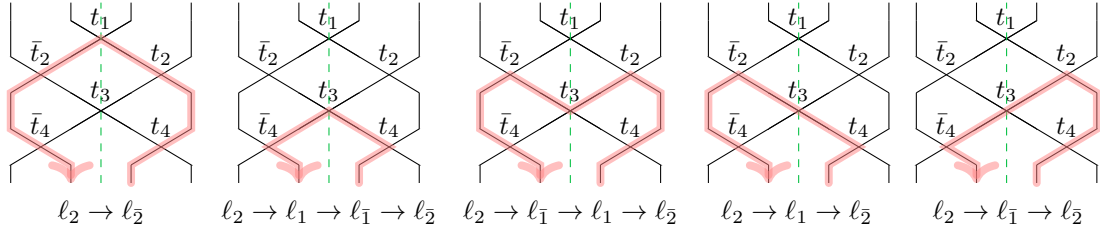
\begin{tabular}{ccccc}
\begin{tikzpicture}[scale = 0.8, yscale=0.6]
\def\nwires{2} 
\def\length{4} 

\foreach \i\y in {1/2,2/1,3/2,4/1}{
	\pgfmathsetmacro\ycoordinate{\length-\i+1}
	\pgfmathsetmacro\yplusone{\y+1}
	\foreach \x in {1,...,\nwires}{
		\ifthenelse{\x = \y}
			{\coordinate (\i LSW) at (-\nwires+\x-0.5, \ycoordinate);
			\coordinate (\i LSE) at (-\nwires+\x +0.5, \ycoordinate);
			\coordinate (\i LNW) at (-\nwires+\x-0.5, \ycoordinate+1);
			\coordinate (\i LNE) at (-\nwires+\x +0.5,\ycoordinate+1);

			\coordinate (\i RSW) at (\nwires - \x -0.5, \ycoordinate);
			\coordinate (\i RSE) at (\nwires - \x + 0.5, \ycoordinate);
			\coordinate (\i RNW) at (\nwires - \x - 0.5, \ycoordinate + 1);
			\coordinate (\i RNE) at (\nwires - \x + 0.5, \ycoordinate + 1);

			\coordinate (\i LC) at (-\nwires + \x, \ycoordinate+ 0.5);
			\coordinate (\i RC) at (\nwires - \x, \ycoordinate+0.5);
		
			\draw (-\nwires+\x-0.5,\ycoordinate) -- (\i LSW)-- (\i LNE) -- (-\nwires+\x+0.5, \ycoordinate+1);
			\draw (-\nwires+\x+0.5,\ycoordinate) --(\i LSE) --(\i LNW) --(-\nwires+\x-0.5, \ycoordinate+1);

			\draw (\nwires - \x -0.5, \ycoordinate) -- (\i RSW)-- (\i RNE) -- (\nwires - \x + 0.5, \ycoordinate + 1);
			\draw (\nwires - \x +0.5, \ycoordinate) -- (\i RSE)-- (\i RNW) -- (\nwires - \x -0.5, \ycoordinate + 1);
			}{}

		\ifthenelse{\x=\y \OR \x = \yplusone}{}{
			\draw
			(-\nwires+\x-0.5,\ycoordinate)--(-\nwires+\x-0.5,\ycoordinate+1);
			\draw
			(\nwires-\x+0.5,\ycoordinate)--(\nwires-\x+0.5,\ycoordinate+1);
		}}

	\ifthenelse{\y=\nwires}{
		\node[above] at (\i LC) {$t_{\i}$};}
	{\node[above] at (\i LC) {$\bar{t}_{\i}$};
	\node[above] at (\i RC) {$t_{\i}$};}

}

\foreach \x in{1,...,\nwires}{
	\draw (\x-0.5,\length+1)--(\x-0.5,\length+1.5);
	\draw (\x-0.5,1)--(\x-0.5,0.5);

	\draw (-\x+0.5,\length+1)--(-\x+0.5,\length+1.5);
	\draw (-\x+0.5,1)--(-\x+0.5,0.5);

}

\draw[teal!50!green, dashed] (0,0.5)--(0,\length+1.5);

\draw[color=red!50, line cap=round, line width=3, opacity=0.5, ->]
(0.5,0.5)--(4RSW)--(4RNE)--(2RSE)--(2RNW)--(1RSE)--(1RC)--(1RSW)--(2LNE)--(2LSW)--(4LNW)--(4LSE)--(-0.5,0.5);

\end{tikzpicture}
&
\begin{tikzpicture}[scale = 0.8, yscale=0.6]
\def\nwires{2} 
\def\length{4} 

\foreach \i\y in {1/2,2/1,3/2,4/1}{
	\pgfmathsetmacro\ycoordinate{\length-\i+1}
	\pgfmathsetmacro\yplusone{\y+1}
	\foreach \x in {1,...,\nwires}{
		\ifthenelse{\x = \y}
			{\coordinate (\i LSW) at (-\nwires+\x-0.5, \ycoordinate);
			\coordinate (\i LSE) at (-\nwires+\x +0.5, \ycoordinate);
			\coordinate (\i LNW) at (-\nwires+\x-0.5, \ycoordinate+1);
			\coordinate (\i LNE) at (-\nwires+\x +0.5,\ycoordinate+1);

			\coordinate (\i RSW) at (\nwires - \x -0.5, \ycoordinate );
			\coordinate (\i RSE) at (\nwires - \x + 0.5, \ycoordinate );
			\coordinate (\i RNW) at (\nwires - \x - 0.5, \ycoordinate +1);
			\coordinate (\i RNE) at (\nwires - \x + 0.5, \ycoordinate +1);

			\coordinate (\i LC) at (-\nwires + \x, \ycoordinate+ 0.5);
			\coordinate (\i RC) at (\nwires - \x, \ycoordinate+0.5);
		
			\draw (-\nwires+\x-0.5,\ycoordinate) -- (\i LSW)-- (\i LNE) -- (-\nwires+\x+0.5, \ycoordinate+1);
			\draw (-\nwires+\x+0.5,\ycoordinate) --(\i LSE) --(\i LNW) --(-\nwires+\x-0.5, \ycoordinate+1);

			\draw (\nwires - \x -0.5, \ycoordinate) -- (\i RSW)-- (\i RNE) -- (\nwires - \x + 0.5, \ycoordinate + 1);
			\draw (\nwires - \x +0.5, \ycoordinate) -- (\i RSE)-- (\i RNW) -- (\nwires - \x -0.5, \ycoordinate + 1);
			}{}

		\ifthenelse{\x=\y \OR \x = \yplusone}{}{
			\draw
			(-\nwires+\x-0.5,\ycoordinate)--(-\nwires+\x-0.5,\ycoordinate+1);
			\draw
			(\nwires-\x+0.5,\ycoordinate)--(\nwires-\x+0.5,\ycoordinate+1);
		}}

	\ifthenelse{\y=\nwires}{
		\node[above] at (\i LC) {$t_{\i}$};}
	{\node[above] at (\i LC) {$\bar{t}_{\i}$};
	\node[above] at (\i RC) {$t_{\i}$};}

}

\foreach \x in{1,...,\nwires}{
	\draw (\x-0.5,\length+1)--(\x-0.5,\length+1.5);
	\draw (\x-0.5,1)--(\x-0.5,0.5);

	\draw (-\x+0.5,\length+1)--(-\x+0.5,\length+1.5);
	\draw (-\x+0.5,1)--(-\x+0.5,0.5);
}

\draw[teal!50!green, dashed] (0,0.5)--(0,\length+1.5);
\draw[color=red!50, line cap=round, line width=3, opacity=0.5, ->]
(0.5,0.5)--(4RSW)--(4RC)--(4RNW)--(3RSE)--(3RC)--(3RSW)--(4LNE)--(4LC)--(4LSE)--(-0.5,0.5);

\end{tikzpicture}
&
\begin{tikzpicture}[scale = 0.8, yscale=0.6]
\def\nwires{2} 
\def\length{4} 

\foreach \i\y in {1/2,2/1,3/2,4/1}{
	\pgfmathsetmacro\ycoordinate{\length-\i+1}
	\pgfmathsetmacro\yplusone{\y+1}
	\foreach \x in {1,...,\nwires}{
		\ifthenelse{\x = \y}
			{\coordinate (\i LSW) at (-\nwires+\x-0.5, \ycoordinate);
			\coordinate (\i LSE) at (-\nwires+\x +0.5, \ycoordinate);
			\coordinate (\i LNW) at (-\nwires+\x-0.5, \ycoordinate+1);
			\coordinate (\i LNE) at (-\nwires+\x +0.5,\ycoordinate+1);

			\coordinate (\i RSW) at (\nwires - \x -0.5, \ycoordinate );
			\coordinate (\i RSE) at (\nwires - \x + 0.5, \ycoordinate );
			\coordinate (\i RNW) at (\nwires - \x - 0.5, \ycoordinate +1);
			\coordinate (\i RNE) at (\nwires - \x + 0.5, \ycoordinate +1);

			\coordinate (\i LC) at (-\nwires + \x, \ycoordinate+ 0.5);
			\coordinate (\i RC) at (\nwires - \x, \ycoordinate+0.5);
		
			\draw (-\nwires+\x-0.5,\ycoordinate) -- (\i LSW)-- (\i LNE) -- (-\nwires+\x+0.5, \ycoordinate+1);
			\draw (-\nwires+\x+0.5,\ycoordinate) --(\i LSE) --(\i LNW) --(-\nwires+\x-0.5, \ycoordinate+1);

			\draw (\nwires - \x -0.5, \ycoordinate) -- (\i RSW)-- (\i RNE) -- (\nwires - \x + 0.5, \ycoordinate + 1);
			\draw (\nwires - \x +0.5, \ycoordinate) -- (\i RSE)-- (\i RNW) -- (\nwires - \x -0.5, \ycoordinate + 1);
			}{}

		\ifthenelse{\x=\y \OR \x = \yplusone}{}{
			\draw
			(-\nwires+\x-0.5,\ycoordinate)--(-\nwires+\x-0.5,\ycoordinate+1);
			\draw
			(\nwires-\x+0.5,\ycoordinate)--(\nwires-\x+0.5,\ycoordinate+1);
		}}

	\ifthenelse{\y=\nwires}{
		\node[above] at (\i LC) {$t_{\i}$};}
	{\node[above] at (\i LC) {$\bar{t}_{\i}$};
	\node[above] at (\i RC) {$t_{\i}$};}

}

\foreach \x in{1,...,\nwires}{
	\draw (\x-0.5,\length+1)--(\x-0.5,\length+1.5);
	\draw (\x-0.5,1)--(\x-0.5,0.5);

	\draw (-\x+0.5,\length+1)--(-\x+0.5,\length+1.5);
	\draw (-\x+0.5,1)--(-\x+0.5,0.5);
}

\draw[teal!50!green, dashed] (0,0.5)--(0,\length+1.5);
\draw[color=red!50, line cap=round, line width=3, opacity=0.5, ->]
(0.5,0.5)--(4RSW)--(4RNE)--(2RSE)--(2RC)--(2RSW)--(3RNE)--(3RC)--(3RNW)--(2LSE)--(2LC)--(2LSW)--(4LNW)--(4LSE)--(-0.5,0.5);
\end{tikzpicture}
&
\begin{tikzpicture}[scale = 0.8, yscale=0.6]
\def\nwires{2} 
\def\length{4} 

\foreach \i\y in {1/2,2/1,3/2,4/1}{
	\pgfmathsetmacro\ycoordinate{\length-\i+1}
	\pgfmathsetmacro\yplusone{\y+1}
	\foreach \x in {1,...,\nwires}{
		\ifthenelse{\x = \y}
			{\coordinate (\i LSW) at (-\nwires+\x-0.5, \ycoordinate);
			\coordinate (\i LSE) at (-\nwires+\x +0.5, \ycoordinate);
			\coordinate (\i LNW) at (-\nwires+\x-0.5, \ycoordinate+1);
			\coordinate (\i LNE) at (-\nwires+\x +0.5,\ycoordinate+1);

			\coordinate (\i RSW) at (\nwires - \x -0.5, \ycoordinate );
			\coordinate (\i RSE) at (\nwires - \x + 0.5, \ycoordinate );
			\coordinate (\i RNW) at (\nwires - \x - 0.5, \ycoordinate +1);
			\coordinate (\i RNE) at (\nwires - \x + 0.5, \ycoordinate +1);

			\coordinate (\i LC) at (-\nwires + \x, \ycoordinate+ 0.5);
			\coordinate (\i RC) at (\nwires - \x, \ycoordinate+0.5);
		
			\draw (-\nwires+\x-0.5,\ycoordinate) -- (\i LSW)-- (\i LNE) -- (-\nwires+\x+0.5, \ycoordinate+1);
			\draw (-\nwires+\x+0.5,\ycoordinate) --(\i LSE) --(\i LNW) --(-\nwires+\x-0.5, \ycoordinate+1);

			\draw (\nwires - \x -0.5, \ycoordinate) -- (\i RSW)-- (\i RNE) -- (\nwires - \x + 0.5, \ycoordinate + 1);
			\draw (\nwires - \x +0.5, \ycoordinate) -- (\i RSE)-- (\i RNW) -- (\nwires - \x -0.5, \ycoordinate + 1);
			}{}

		\ifthenelse{\x=\y \OR \x = \yplusone}{}{
			\draw
			(-\nwires+\x-0.5,\ycoordinate)--(-\nwires+\x-0.5,\ycoordinate+1);
			\draw
			(\nwires-\x+0.5,\ycoordinate)--(\nwires-\x+0.5,\ycoordinate+1);
		}}

	\ifthenelse{\y=\nwires}{
		\node[above] at (\i LC) {$t_{\i}$};}
	{\node[above] at (\i LC) {$\bar{t}_{\i}$};
	\node[above] at (\i RC) {$t_{\i}$};}

}

\foreach \x in{1,...,\nwires}{
	\draw (\x-0.5,\length+1)--(\x-0.5,\length+1.5);
	\draw (\x-0.5,1)--(\x-0.5,0.5);

	\draw (-\x+0.5,\length+1)--(-\x+0.5,\length+1.5);
	\draw (-\x+0.5,1)--(-\x+0.5,0.5);
}

\draw[teal!50!green, dashed] (0,0.5)--(0,\length+1.5);
\draw[color=red!50, line cap=round, line width=3, opacity=0.5, ->]
(0.5,0.5)--(4RSW)--(4RC)--(4RNW)--(3RSE)--(3RC)--(3RNW)--(2LSE)--(2LC)--(2LSW)--(4LNW)--(4LSE)--(-0.5,0.5);
\end{tikzpicture}
&
\begin{tikzpicture}[scale = 0.8, yscale=0.6]
\def\nwires{2} 
\def\length{4} 

\foreach \i\y in {1/2,2/1,3/2,4/1}{
	\pgfmathsetmacro\ycoordinate{\length-\i+1}
	\pgfmathsetmacro\yplusone{\y+1}
	\foreach \x in {1,...,\nwires}{
		\ifthenelse{\x = \y}
			{\coordinate (\i LSW) at (-\nwires+\x-0.5, \ycoordinate);
			\coordinate (\i LSE) at (-\nwires+\x +0.5, \ycoordinate);
			\coordinate (\i LNW) at (-\nwires+\x-0.5, \ycoordinate+1);
			\coordinate (\i LNE) at (-\nwires+\x +0.5,\ycoordinate+1);

			\coordinate (\i RSW) at (\nwires - \x -0.5, \ycoordinate );
			\coordinate (\i RSE) at (\nwires - \x + 0.5, \ycoordinate );
			\coordinate (\i RNW) at (\nwires - \x - 0.5, \ycoordinate +1);
			\coordinate (\i RNE) at (\nwires - \x + 0.5, \ycoordinate +1);

			\coordinate (\i LC) at (-\nwires + \x, \ycoordinate+ 0.5);
			\coordinate (\i RC) at (\nwires - \x, \ycoordinate+0.5);
		
			\draw (-\nwires+\x-0.5,\ycoordinate) -- (\i LSW)-- (\i LNE) -- (-\nwires+\x+0.5, \ycoordinate+1);
			\draw (-\nwires+\x+0.5,\ycoordinate) --(\i LSE) --(\i LNW) --(-\nwires+\x-0.5, \ycoordinate+1);

			\draw (\nwires - \x -0.5, \ycoordinate) -- (\i RSW)-- (\i RNE) -- (\nwires - \x + 0.5, \ycoordinate + 1);
			\draw (\nwires - \x +0.5, \ycoordinate) -- (\i RSE)-- (\i RNW) -- (\nwires - \x -0.5, \ycoordinate + 1);
			}{}

		\ifthenelse{\x=\y \OR \x = \yplusone}{}{
			\draw
			(-\nwires+\x-0.5,\ycoordinate)--(-\nwires+\x-0.5,\ycoordinate+1);
			\draw
			(\nwires-\x+0.5,\ycoordinate)--(\nwires-\x+0.5,\ycoordinate+1);
		}}

	\ifthenelse{\y=\nwires}{
		\node[above] at (\i LC) {$t_{\i}$};}
	{\node[above] at (\i LC) {$\bar{t}_{\i}$};
	\node[above] at (\i RC) {$t_{\i}$};}

}

\foreach \x in{1,...,\nwires}{
	\draw (\x-0.5,\length+1)--(\x-0.5,\length+1.5);
	\draw (\x-0.5,1)--(\x-0.5,0.5);

	\draw (-\x+0.5,\length+1)--(-\x+0.5,\length+1.5);
	\draw (-\x+0.5,1)--(-\x+0.5,0.5);
}

\draw[teal!50!green, dashed] (0,0.5)--(0,\length+1.5);
\draw[color=red!50, line cap=round, line width=3, opacity=0.5, ->]
(0.5,0.5)--(4RSW)--(4RNE)--(2RSE)--(2RC)--(2RSW)--(3RNE)--(3RSW)--(4LNE)--(4LC)--(4LSE)--(-0.5,0.5); 
\end{tikzpicture} \\
$\ell_2 \to \ell_{\bar{2}}$ 
& $\ell_{2} \to \ell_1 \to \ell_{\bar{1}} \to \ell_{\bar{2}}$ 
& $\ell_2 \to \ell_{\bar{1}} \to \ell_1 \to \ell_{\bar{2}}$
& $\ell_2 \to \ell_1 \to \ell_{\bar{2}}$
& $\ell_2 \to \ell_{\bar{1}} \to \ell_{\bar{2}}$
\end{tabular}
The first three paths are symmetric while the last two are not. We note that for the fourth path $P = (\ell_2 \to \ell_1 \to \ell_{\bar{2}})$, its mirror $P^{\vee}$ is the fifth path $(\ell_2 \to \ell_{\bar{1}} \to \ell_{\bar{2}})$.
\end{example}

Recall that $\R^{n(2n-1)}$ and $\R^N$ with $N = n^2$ are vector spaces in which string cones 
$\mathcal{C}_{\bm i}^{(A_{2n-1})}$ and $\mathcal{C}_{\bm i}^{(C_n)}$ live in with the coordinate systems 
${\bm t}$ and ${\bm a}$, respectively. 

\begin{lemma}\label{lemma_hat_F_P_C_and_hat_P}
For ${\bm i} \in \RCn$, define a linear map $\Psi_{\bm i} \colon (\R^{n(2n-1)})^* \to (\R^{N})^*$ by 
\begin{equation}\label{equation_def_Psi}
\Psi_{\bm i}(t_j) = \begin{cases}
a_j & \text{ if } i_j \neq n, \\
2a_j & \text{ if } i_j = n;
\end{cases} 
\quad 
\Psi_{\bm i}(\bar{t}_j) = a_j.
\end{equation}
where $t_j$ and $\bar{t}_j$ are $t_j$-th and $\bar{t}_j$-th coordinate functions of ${\bf t}$ on $\R^{n(2n-1)}$, respectively. 
We similarly regard $a_j$ the $j$-th coorinate function of ${\bf a}$ on $\R^N$.
Then, for $P \in \Paths^{\symp}(\bm i, k)$, it holds that 
\[
\hat{\F}_{\hat{P}}^{(C)}({\bm a})
= \Psi_{\bm i}(\hat{\F}_P({\bm t})).
\]
\end{lemma}
\begin{proof}
Recall from Definition~\ref{def_string_inequalities_BC} that  $\hat{\F}_{\hat{P}}^{(C)}({\bm a}) = \hat{\F}_{\hat{P}}^{(B)}(\Gamma_{\bm i}^{C,B}({\bm a})) = \hat{\F}_{\hat{P}}(\Upsilon_{\bm i}^{B,A}(\Gamma_{\bm i}^{C,B}({\bm a})))$. 
Moreover, the composition $\Upsilon_{\bm i}^{B,A}\circ \Gamma_{\bm i}^{C,B}$ sends ${\bm a} = (a_1,\dots,a_N)$ to 
\[
(\underbrace{m_{i_1}'a_1,\dots,m_{i_1}'a_1}_{m_{i_1}},
\dots,\underbrace{m_{i_N}'a_N,\dots,m_{i_N}'a_N}_{m_{i_N}}). 
\]
Here, $m_{i_j} = 2$ if $i_j \neq n$; $m_{i_j} = 1$ if $i_j = n$, and moreover,  $m_{i_j}' = 1$ if $i_j \neq n$; $m_{i_j} = 2$ if $i_j = n$. 
Indeed, if $i_j \neq n$, then the coordinate $a_j$ sends to $(a_j, a_j)$. If $i_j =n$, then the coordinate $a_j$ sends to $2a_j$. 

On the other hand, when we draw the symplectic wiring diagram $G^{\symp}(\bm i)$, we add two crossings on $i_j$th and $(2n-i_j)$th columns if $i_j \neq n$; otherwise, we add one crossing on $n$th column. Accordingly, by substituting $a_j$ for $t_j$ and $\bar{t}_j$ if $i_j \neq n$; $2a_j$ for $t_j$ if $i_j = n$ in the function $\hat{\F}_P({\bm t})$, we obtain the function $\hat{\F}_{\hat{P}}^{(C)}({\bm a})$. This proves the claim. 
\end{proof}

\begin{example}\label{example_2121}
Let $\bm i = (2,1,2,1)$. Then the map $\Psi_{\bm i}$ is given by 
\[
\Psi_{\bm i}(t_1) = 2a_1, \Psi_{\bm i}(t_2) = \Psi_{\bm i}(\bar{t}_2) = a_2, 
\Psi_{\bm i}(t_3) = 2a_3, \Psi_{\bm i}(t_4) = \Psi_{\bm i}(\bar{t}_4) = a_4.
\]
The five rigorous paths in $\Paths(\bm i, 2)$ in Example~\ref{example_symmetric_paths} provide the functions $\hat{\F}_P({\bm t})$ and $\hat{\F}_{\hat{P}}^{(C)}({\bm a})$ as shown in the third and fourth columns in Table~\ref{table_2121}.
We notice that the first three symmetric paths produce functions $\hat{\F}^{(C)}_P({\bm a})$ divisible by $2$, and the last two nonsymmetric paths provide redundant inequalities because 
$(a_2 - a_3) + (a_3 - a_4) = a_2 - a_4$.
In addition, $\Paths^{\symp}(\bm i, 1)$ consists of only one rigorous path $\ell_1 \to \ell_2$ that provides a non-redundant inequality $a_4 \geq 0$. 
Thus, the number of facets of the string cone $\mathcal{C}_{\bm i}^{(C_2)}$ is $4$. 
\end{example}

\begin{lemma}\label{lemma_lift_of_FP}
Let $\bm i \in R(w_0^{(C_n)})$ and $P\in \Paths^{\symp}(\bm i, k)$. 
Then it holds that $\hat{\F}_{\hat{P}}^{(C)}({\bm a}) = \hat{\F}_{{\hat{P}}^{\vee}}^{(C)}({\bm a})$. 
Moreover, if $P = P^{\vee}$, then the coefficients of $\hat{\F}_{\hat{P}}^{(C)}({\bm a})$ are all even.
\end{lemma}

\begin{proof}
Let $P = (\ell_{r_1} \to \dots \to \ell_{r_{s+1}})$. 
Consider the node $\ell_{r_x} \cap \ell_{r_{x+1}}$ on $P$ for $1 \le x \le s$. Denote by $b_{u_x}$ the node $\ell_{r_x}\cap \ell_{r_{x+1}}$ where $b = t$ or $\bar{t}$. 
Then, for the mirror $P^{\vee} = (\ell_{\overline{r_{s+1}}} \to \dots \to \ell_{\overline{r_1}})$, we have nodes $\{\bar{b}_{u_x} \mid 1\leq x \leq s\}$. Here, we set $\bar{\bar{t}} = t$. 
For $1 \le x \le s$, we have $r_x < r_{x+1}$ if and only if $\overline{r_x} > \overline{r_{x+1}}$ because of 
\[
\overline{r_x} = 2n+1-r_x > 2n+1-r_{x+1} = \overline{r_{x+1}}.
\]
Accordingly, the coefficient of $b_{u_x}$ in $\hat{\F}_P({\bm t})$ is the same as that of $\bar{b}_{u_k}$ in $\hat{\F}_{P^{\vee}}({\bm t})$: 
\begin{equation}\label{eq_FP_and_FPvee}
\hat{\F}_P({\bm t}) 
= \sum_{x=1}^s A_x b_{u_x}\quad \text{if and only if}\quad 
\hat{\F}_{P^{\vee}}({\bm t})
= \sum_{x=1}^s A_x \bar{b}_{u_k}.
\end{equation}
Moreover, if the node $b_{r_k}$ is on the wall, then the node $\bar{b}_{r_k}$ is also on the wall. 
Because of the definition of $\Psi_{\bm i}$ in~\eqref{equation_def_Psi}, we have
\begin{equation}\label{eq_Psi_image}
\begin{cases}
\Psi_{\bm i}(b_{u_x}) = \Psi_{\bm i}(\overline{b}_{u_x}) = a_{u_x} & \text { if }i_{u_x} \neq n; \\
\Psi_{\bm i}(b_{u_x}) = 2 a_{u_x} & \text{ if } i_{u_x} = n. 
\end{cases}
\end{equation}
Therefore, by applying Lemma~\ref{lemma_hat_F_P_C_and_hat_P} on functions in~\eqref{eq_FP_and_FPvee}, we obtain the desired equality $\hat{\F}_{\hat{P}}^{(C)}({\bm a}) = \hat{\F}_{\hat{P}^{\vee}}^{(C)}({\bm a})$. 

If $P = P^{\vee}$, then it follows that $\{b_{u_x} \mid 1 \le x \le s\} = \{\bar{b}_{u_x} \mid 1 \le x \le s\}$. 
This implies that if $i_{u_x} \neq n$, then both $b_{u_x}$ and $\bar{b}_{u_x}$ have the same nonzero coefficient in $\hat{\F}_P(\bm t)$. Accordingly, after applying the map $\Psi_{\bm i}$ to $\hat{\F}_{P}(\bm t)$, we have the function $\hat{\F}_{\hat{P}}^{(C)}(\bm a)$ whose coefficients are all even. This proves the claim. 
\end{proof}

By Lemma~\ref{lemma_lift_of_FP}, we can divide the function $\hat{\F}_{\hat{P}}^{(C)}({\bm a})$ by $2$ if $P = P^{\vee}$. We now define the following. 

\begin{definition}\label{def_FPC}
For $\bm i \in \RCn$ and $P \in \Paths^{\symp}(\bm i, k)$, set 
\[
\F_P^{(C)}(\bm a) \colonequals \begin{cases}
\hat{\F}_{\hat{P}}^{(C)}({\bm a}) & \text{ if } P \neq P^{\vee}, \\
\frac{1}{2} \hat{\F}_{\hat{P}}^{(C)}({\bm a}) & \text{ if }P = P^{\vee}.
\end{cases}
\]
\end{definition}

For instance, the functions $\F_{P}^{(C)}(\bm a)$ for $P$ in $\Paths^{\symp}((2,1,2,1),2)$ are given in the fifth column in Table~\ref{table_2121}. 
\begin{table}[H]
\begin{tabular}{l|l|lll}
\toprule
$P$&symmetric &$\hat{\F}_P({\bm t})$ & $\hat{\F}_{\hat{P}}^{(C)}({\bm a})$ & $\F_P^{(C)}(\bm a)$ \\ 
\midrule
$\ell_2 \to \ell_{\bar{2}}$ & yes& $t_1$ & $2a_1$ & $a_1$\\
$\ell_2 \to \ell_1 \to \ell_{\bar{1}} \to \ell_{\bar{2}}$ & yes&$t_3 - (t_4 + \bar{t}_4)$ & $2a_3 - 2a_4$ & $a_3 - a_4$\\
$\ell_2 \to \ell_{\bar{1}} \to \ell_1 \to \ell_{\bar{2}}$& yes&$(t_2 + \bar{t}_2) - t_3$ & $2a_2 - 2a_3$ & $a_2 - a_3$\\
$\ell_2 \to \ell_1 \to \ell_{\bar{2}}$ & no&$\bar{t}_2 - t_4$ & $a_2 - a_4$ &$a_2 - a_4$ \\
$\ell_2 \to \ell_{\bar{1}} \to \ell_{\bar{2}}$ & no&$t_2 - \bar{t}_4$ & $a_2 - a_4$ & $a_2 - a_4$ \\
\bottomrule
\end{tabular}
\caption{Functions $\hat{\F}_P$, $\hat{\F}_{\hat{P}}^{(C)}$, and $\F_P^{(C)}$ for $P \in \Paths^{\symp}((2,1,2,1),2)$}\label{table_2121}
\end{table}

The following corollary is an immediate consequence of Theorem \ref{thm_string_cones_BC}. 

\begin{corollary}\label{cor_string_cone_inequality_C}
For ${\bm i} \in R(w_0^{(C_n)})$, the following equality holds: 
\begin{equation}\label{eq:descriptions_of_string_cones_in_type_C_n}
\begin{aligned}
&\mathcal{C}_{\bm i}^{(C_n)} = \{{\bm a} \in \R^N \mid \F_P^{(C)}({\bm a}) \geq 0 \quad \text{ for all } 1 \leq k \leq n \text{ and } P \in \Paths^{\symp}(\bm i, k)\}.
\end{aligned}
\end{equation}
\end{corollary}

In the rest of this section, we study non-redundancy of the inequality $\F_P^{(C)}({\bm a}) \geq 0$ in the expression \eqref{eq:descriptions_of_string_cones_in_type_C_n}.

\begin{proposition}\label{prop_non_redundancy_k_not_n_maximal}
Let $\bm i \in R(w_0^{(C_n)})$, and $P \in \Paths^{\symp}(\bm i,k)$. 
If there is no $Q \in \Paths^{\symp}(\bm i,k)$ such that $\chamber(Q) \subseteq \chamber(P) \cup \chamber(P^\vee)$ and such that $\chamber(Q) \nsubseteq \chamber(P)$, then the inequality $\F_P^{(C)}({\bm a}) \geq 0$ is non-redundant in the expression \eqref{eq:descriptions_of_string_cones_in_type_C_n}.
\end{proposition}

\begin{proof}
Assume on the contrary that $\F_P^{(C)}({\bm a}) \geq 0$ is redundant. 
Then there exist positive numbers $q,p \in \Z_{>0}$ (with $q>1$) and $\varepsilon_1,\dots,\varepsilon_q \in \mathbb{Q}_{>0}$ and distinct rigorous paths $Q_1,\dots,Q_q$, which are not $P$, such that 
\begin{equation}\label{equation_redundant}
{\F}_{{P}}^{(C)}({\bm a})  = \varepsilon_1{\F}_{{Q}_{1}}^{(C)}({\bm a}) + \cdots + \varepsilon_p{\F}_{{Q}_p}^{(C)}({\bm a}) + 
\varepsilon_{p+1}{\F}_{{Q}_{p+1}}^{(C)}({\bm a})  + \cdots + \varepsilon_q{\F}_{{Q}_q}^{(C)}({\bm a})
\end{equation}
with $Q_r^\vee \neq Q_r$ for $1 \leq r \leq p$ and $Q_r^\vee = Q_r$ for $p+1 \leq r \leq q$. 
After multiplying $2$ on both sides of~\eqref{equation_redundant} and applying Lemma~\ref{lemma_hat_F_P_C_and_hat_P}, we obtain
\[
\Psi_{\bm i}(
2\hat{\F}_{{P}}({\bm t})
) = 
\Psi_{\bm i}(
2 \varepsilon_1\hat{\F}_{{Q}_1}({\bm t}) + \dots + 2 \varepsilon_p\hat\F_{{Q}_p}({\bm t}) +\varepsilon_{p+1} \hat\F_{{Q}_{p+1}}({\bm t}) + \dots+\varepsilon_{q}\hat{\F}_{{Q}_{q}}({\bm t})
).
\]
Since $\Psi_{\bm i}(\hat\F_{{P}}({\bm t})) = \hat\F_{\hat{P}}^{(C)}({\bm a}) = \hat{\F}_{\hat{P}^{\vee}}^{(C)}({\bm a}) = \Psi_{\bm i}(\hat\F_{P^{\vee}}({\bm t}))$ by Lemma~\ref{lemma_lift_of_FP}, we have
\begin{equation}\label{eq_non_redundant_equality_in_a}
	\begin{split}
\Psi_{\bm i}(\hat{\F}_{{P}}({\bm t}) + \hat{\F}_{{P}^\vee}({\bm t})) &= 
\Psi_{\bm i}(\varepsilon_1(\hat{\F}_{{Q}_1}({\bm t}) + \hat{\F}_{{Q}_1^\vee}({\bm t})) + \dots + \varepsilon_p(\hat\F_{{Q}_{p}}({\bm t}) + \hat\F_{{Q}_p^\vee}({\bm t})) \\
& \qquad \qquad +\varepsilon_{p+1} \hat\F_{{Q}_{p+1}}({\bm t}) + \dots+\varepsilon_q\hat{\F}_{{Q}_{q}}({\bm t})).
\end{split}
\end{equation}
Since $\hat{\F}_{{P}}({\bm t}) + \hat{\F}_{{P}^\vee}({\bm t})$ is a polynomial in $t_j + \bar{t}_j$ with $i_j \neq n$ and in $t_j$ with $i_j = n$, it is obtained from $\Psi_{\bm i}(\hat{\F}_{{P}}({\bm t}) + \hat{\F}_{{P}^\vee}({\bm t}))$ by the following substitutions:
\[a_j \mapsto \begin{cases}
\frac{1}{2} (t_j + \bar{t}_j)\quad(i_j \neq n),\\
\frac{1}{2} t_j \quad(i_j = n).
\end{cases}\]
Similarly, we deduce $\varepsilon_1(\hat{\F}_{{Q}_1}({\bm t}) + \hat{\F}_{{Q}_1^\vee}({\bm t})) + \dots + \varepsilon_p(\hat\F_{{Q}_p}({\bm t}) + \hat\F_{{Q}_p^\vee}({\bm t})) +\varepsilon_{p+1} \hat\F_{{Q}_{p+1}}({\bm t}) + \dots+\varepsilon_q\hat{\F}_{{Q}_{q}}({\bm t})$ from the right hand side of \eqref{eq_non_redundant_equality_in_a} by the same procedure, which implies that 
\[
\hat{\F}_{{P}}({\bm t}) + \hat{\F}_{{P}^\vee}({\bm t}) = 
\varepsilon_1(\hat{\F}_{Q_1}({\bm t}) + \hat{\F}_{{Q}_1^\vee}({\bm t})) + \dots + \varepsilon_p(\hat\F_{{Q}_p}({\bm t}) + \hat\F_{{Q}_p^\vee}({\bm t})) +\varepsilon_{p+1} \hat\F_{{Q}_{p+1}}({\bm t}) + \dots+\varepsilon_q\hat{\F}_{{Q}_{q}}({\bm t}).
\]

Considering the change $\Phi$ of coordinates in Lemma~\ref{lemma_chamber_varialbe}, we obtain 
\begin{equation}\label{eq_lift_of_P_and_chambers}
\chamber_{P} + \chamber_{P^\vee}
= \varepsilon_1(\chamber_{Q_1} + \chamber_{Q_1^\vee}) + \cdots + \varepsilon_p(\chamber_{Q_p} + \chamber_{Q_p^\vee}) + \varepsilon_{p+1}\chamber_{{Q}_{p+1}} + \cdots + \varepsilon_q\chamber_{{Q}_{q}},
\end{equation}
where $\chamber_P \colonequals \sum_{\chamber_j \subseteq \chamber(P)} u_j$.
We notice that every path $Q \in \Paths(\hat{\bm i},k)$ starts at $L_k$ and ends at $L_{k+1}$. This implies that $\chamber(Q)$ contains the chamber having $L_k$ and $L_{k+1}$ if and only if $Q \in \Paths(\hat{\bm i},k)$. 
Since the union $\chamber(P) \cup \chamber(P^\vee)$ does not contain the chamber having $L_b$ and $L_{b+1}$ for $b \neq k, 2n-k$, we see by \eqref{eq_lift_of_P_and_chambers} that 
\[
	\hat{Q}_1, \ldots, \hat{Q}_q \in \Paths(\hat{\bm i},k) \cup \Paths(\hat{\bm i}, 2n-k).
\] 
By re-labeling on $Q_u$'s, if necessary, we may assume that $\hat{Q}_1, \ldots, \hat{Q}_q \in \Paths(\hat{\bm i},k)$. Then $\hat{Q}_1, \ldots, \hat{Q}_q \in \Paths(\hat{\bm i},k)$ contain the chamber having $L_k$ and $L_{k+1}$
and hence we have $\varepsilon_1 + \cdots + \varepsilon_q = 1$ by~\eqref{eq_lift_of_P_and_chambers}. 
Since $\chamber(Q_u) \subseteq \chamber(P)$ by the assumption on $P$ for each $1 \leq u \leq q$, we see by \eqref{eq_lift_of_P_and_chambers} that $\chamber(Q_u) = \chamber(P)$ for all $1 \leq u \leq q$. 
Indeed, if $\chamber(Q_u) \subsetneq \chamber(P)$ for some $1 \leq u \leq q$, then there exists a chamber $\mathscr{D} \subseteq (\chamber(P) \cup \chamber(P^\vee)) \setminus (\chamber(Q_u) \cup \chamber(Q_u^{\vee}))$. 
Comparing the coefficients of the chamber $\mathscr{D}$ in both sides of~\eqref{eq_lift_of_P_and_chambers}, the right hand side of~\eqref{eq_lift_of_P_and_chambers} is strictly less than that in the left hand side of~\eqref{eq_lift_of_P_and_chambers}, which proves that $\chamber(Q_u) \subsetneq \chamber(P)$ cannot occur.
This implies that $Q_u = P$ for all $1 \leq u \leq q$, which contradicts that $Q_u$'s are distinct.
\end{proof}

\begin{remark}\label{rmk_maximal_Q_is_P}
Let $\bm i \in \RCn$, and $P \in \Paths^{\symp}(\bm i,k)$. Suppose that $P$ satisfies the assumption of Proposition~\ref{prop_non_redundancy_k_not_n_maximal}. 
Let $Q \in \Paths^{\symp}(\bm i,k)$ such that $\chamber(Q) \subseteq \chamber(P) \cup \chamber(P^{\vee})$. 
Because of the assumption on $P$, we have $\chamber(Q) \subseteq \chamber(P)$. Accordingly, if $Q$ is maximal in the sense that there does not exist $Q' \in \Paths^{\symp}(\bm i,k)$ such that $\chamber(Q) \subsetneq \chamber(Q') \subseteq \chamber(P) \cup \chamber(P^{\vee})$, then it follows that $Q = P$.
\end{remark}

\begin{figure}[H]
\begin{tikzpicture}[scale = 0.8, yscale=0.6]
\def\nwires{2} 
\def\length{4} 

\foreach \i\y in {1/1,2/2,3/1,4/2}{
	\pgfmathsetmacro\ycoordinate{\length-\i+1}
	\pgfmathsetmacro\yplusone{\y+1}
	\foreach \x in {1,...,\nwires}{
		\ifthenelse{\x = \y}
			{\coordinate (\i LSW) at (-\nwires+\x-0.5, \ycoordinate);
			\coordinate (\i LSE) at (-\nwires+\x +0.5, \ycoordinate);
			\coordinate (\i LNW) at (-\nwires+\x-0.5, \ycoordinate+1);
			\coordinate (\i LNE) at (-\nwires+\x +0.5,\ycoordinate+1);

			\coordinate (\i RSW) at (\nwires - \x -0.5, \ycoordinate);
			\coordinate (\i RSE) at (\nwires - \x + 0.5, \ycoordinate);
			\coordinate (\i RNW) at (\nwires - \x - 0.5, \ycoordinate + 1);
			\coordinate (\i RNE) at (\nwires - \x + 0.5, \ycoordinate + 1);

			\coordinate (\i LC) at (-\nwires + \x, \ycoordinate+ 0.5);
			\coordinate (\i RC) at (\nwires - \x, \ycoordinate+0.5);
		
			\draw (-\nwires+\x-0.5,\ycoordinate) -- (\i LSW)-- (\i LNE) -- (-\nwires+\x+0.5, \ycoordinate+1);
			\draw (-\nwires+\x+0.5,\ycoordinate) --(\i LSE) --(\i LNW) --(-\nwires+\x-0.5, \ycoordinate+1);

			\draw (\nwires - \x -0.5, \ycoordinate) -- (\i RSW)-- (\i RNE) -- (\nwires - \x + 0.5, \ycoordinate + 1);
			\draw (\nwires - \x +0.5, \ycoordinate) -- (\i RSE)-- (\i RNW) -- (\nwires - \x -0.5, \ycoordinate + 1);
			}{}

		\ifthenelse{\x=\y \OR \x = \yplusone}{}{
			\draw
			(-\nwires+\x-0.5,\ycoordinate)--(-\nwires+\x-0.5,\ycoordinate+1);
			\draw
			(\nwires-\x+0.5,\ycoordinate)--(\nwires-\x+0.5,\ycoordinate+1);
		}}

	\ifthenelse{\y=\nwires}{
		\node[above] at (\i LC) {$t_{\i}$};}
	{\node[above] at (\i LC) {$\bar{t}_{\i}$};
	\node[above] at (\i RC) {$t_{\i}$};}

}

\foreach \x in{1,...,\nwires}{
	\draw (\x-0.5,\length+1)--(\x-0.5,\length+1.5);
	\draw (\x-0.5,1)--(\x-0.5,0.5);

	\draw (-\x+0.5,\length+1)--(-\x+0.5,\length+1.5);
	\draw (-\x+0.5,1)--(-\x+0.5,0.5);

}

\draw[teal!50!green, dashed] (0,0.5)--(0,\length+1.5);

\draw[color=red!50, line cap=round, line width=3, opacity=0.5, ->]
(1.5,0.5)--(3RSE)--(3RC)--(4RC)--(4RSE)--(0.5,0.5);

\end{tikzpicture} \quad 
\begin{tikzpicture}[scale = 0.8, yscale=0.6]
\def\nwires{2} 
\def\length{4} 

\foreach \i\y in {1/1,2/2,3/1,4/2}{
	\pgfmathsetmacro\ycoordinate{\length-\i+1}
	\pgfmathsetmacro\yplusone{\y+1}
	\foreach \x in {1,...,\nwires}{
		\ifthenelse{\x = \y}
			{\coordinate (\i LSW) at (-\nwires+\x-0.5, \ycoordinate);
			\coordinate (\i LSE) at (-\nwires+\x +0.5, \ycoordinate);
			\coordinate (\i LNW) at (-\nwires+\x-0.5, \ycoordinate+1);
			\coordinate (\i LNE) at (-\nwires+\x +0.5,\ycoordinate+1);

			\coordinate (\i RSW) at (\nwires - \x -0.5, \ycoordinate);
			\coordinate (\i RSE) at (\nwires - \x + 0.5, \ycoordinate);
			\coordinate (\i RNW) at (\nwires - \x - 0.5, \ycoordinate + 1);
			\coordinate (\i RNE) at (\nwires - \x + 0.5, \ycoordinate + 1);

			\coordinate (\i LC) at (-\nwires + \x, \ycoordinate+ 0.5);
			\coordinate (\i RC) at (\nwires - \x, \ycoordinate+0.5);
		
			\draw (-\nwires+\x-0.5,\ycoordinate) -- (\i LSW)-- (\i LNE) -- (-\nwires+\x+0.5, \ycoordinate+1);
			\draw (-\nwires+\x+0.5,\ycoordinate) --(\i LSE) --(\i LNW) --(-\nwires+\x-0.5, \ycoordinate+1);

			\draw (\nwires - \x -0.5, \ycoordinate) -- (\i RSW)-- (\i RNE) -- (\nwires - \x + 0.5, \ycoordinate + 1);
			\draw (\nwires - \x +0.5, \ycoordinate) -- (\i RSE)-- (\i RNW) -- (\nwires - \x -0.5, \ycoordinate + 1);
			}{}

		\ifthenelse{\x=\y \OR \x = \yplusone}{}{
			\draw
			(-\nwires+\x-0.5,\ycoordinate)--(-\nwires+\x-0.5,\ycoordinate+1);
			\draw
			(\nwires-\x+0.5,\ycoordinate)--(\nwires-\x+0.5,\ycoordinate+1);
		}}

	\ifthenelse{\y=\nwires}{
		\node[above] at (\i LC) {$t_{\i}$};}
	{\node[above] at (\i LC) {$\bar{t}_{\i}$};
	\node[above] at (\i RC) {$t_{\i}$};}

}

\foreach \x in{1,...,\nwires}{
	\draw (\x-0.5,\length+1)--(\x-0.5,\length+1.5);
	\draw (\x-0.5,1)--(\x-0.5,0.5);

	\draw (-\x+0.5,\length+1)--(-\x+0.5,\length+1.5);
	\draw (-\x+0.5,1)--(-\x+0.5,0.5);

}

\draw[teal!50!green, dashed] (0,0.5)--(0,\length+1.5);

\draw[color=red!50, line cap=round, line width=3, opacity=0.5, ->]
(1.5,0.5)--(3RSE)--(2LC)--(3LC)--(4LSE)--(0.5,0.5);

\end{tikzpicture} \quad 
\begin{tikzpicture}[scale = 0.8, yscale=0.6]
\def\nwires{2} 
\def\length{4} 

\foreach \i\y in {1/1,2/2,3/1,4/2}{
	\pgfmathsetmacro\ycoordinate{\length-\i+1}
	\pgfmathsetmacro\yplusone{\y+1}
	\foreach \x in {1,...,\nwires}{
		\ifthenelse{\x = \y}
			{\coordinate (\i LSW) at (-\nwires+\x-0.5, \ycoordinate);
			\coordinate (\i LSE) at (-\nwires+\x +0.5, \ycoordinate);
			\coordinate (\i LNW) at (-\nwires+\x-0.5, \ycoordinate+1);
			\coordinate (\i LNE) at (-\nwires+\x +0.5,\ycoordinate+1);

			\coordinate (\i RSW) at (\nwires - \x -0.5, \ycoordinate);
			\coordinate (\i RSE) at (\nwires - \x + 0.5, \ycoordinate);
			\coordinate (\i RNW) at (\nwires - \x - 0.5, \ycoordinate + 1);
			\coordinate (\i RNE) at (\nwires - \x + 0.5, \ycoordinate + 1);

			\coordinate (\i LC) at (-\nwires + \x, \ycoordinate+ 0.5);
			\coordinate (\i RC) at (\nwires - \x, \ycoordinate+0.5);
		
			\draw (-\nwires+\x-0.5,\ycoordinate) -- (\i LSW)-- (\i LNE) -- (-\nwires+\x+0.5, \ycoordinate+1);
			\draw (-\nwires+\x+0.5,\ycoordinate) --(\i LSE) --(\i LNW) --(-\nwires+\x-0.5, \ycoordinate+1);

			\draw (\nwires - \x -0.5, \ycoordinate) -- (\i RSW)-- (\i RNE) -- (\nwires - \x + 0.5, \ycoordinate + 1);
			\draw (\nwires - \x +0.5, \ycoordinate) -- (\i RSE)-- (\i RNW) -- (\nwires - \x -0.5, \ycoordinate + 1);
			}{}

		\ifthenelse{\x=\y \OR \x = \yplusone}{}{
			\draw
			(-\nwires+\x-0.5,\ycoordinate)--(-\nwires+\x-0.5,\ycoordinate+1);
			\draw
			(\nwires-\x+0.5,\ycoordinate)--(\nwires-\x+0.5,\ycoordinate+1);
		}}

	\ifthenelse{\y=\nwires}{
		\node[above] at (\i LC) {$t_{\i}$};}
	{\node[above] at (\i LC) {$\bar{t}_{\i}$};
	\node[above] at (\i RC) {$t_{\i}$};}

}

\foreach \x in{1,...,\nwires}{
	\draw (\x-0.5,\length+1)--(\x-0.5,\length+1.5);
	\draw (\x-0.5,1)--(\x-0.5,0.5);

	\draw (-\x+0.5,\length+1)--(-\x+0.5,\length+1.5);
	\draw (-\x+0.5,1)--(-\x+0.5,0.5);

}

\draw[teal!50!green, dashed] (0,0.5)--(0,\length+1.5);

\draw[color=red!50, line cap=round, line width=3, opacity=0.5, ->]
(1.5,0.5)--(3RSE)--(1LC)--(1LSW)--(3LNW)--(4LSE)--(0.5,0.5);

\end{tikzpicture} \quad 
\begin{tikzpicture}[scale = 0.8, yscale=0.6]
\def\nwires{2} 
\def\length{4} 

\foreach \i\y in {1/1,2/2,3/1,4/2}{
	\pgfmathsetmacro\ycoordinate{\length-\i+1}
	\pgfmathsetmacro\yplusone{\y+1}
	\foreach \x in {1,...,\nwires}{
		\ifthenelse{\x = \y}
			{\coordinate (\i LSW) at (-\nwires+\x-0.5, \ycoordinate);
			\coordinate (\i LSE) at (-\nwires+\x +0.5, \ycoordinate);
			\coordinate (\i LNW) at (-\nwires+\x-0.5, \ycoordinate+1);
			\coordinate (\i LNE) at (-\nwires+\x +0.5,\ycoordinate+1);

			\coordinate (\i RSW) at (\nwires - \x -0.5, \ycoordinate);
			\coordinate (\i RSE) at (\nwires - \x + 0.5, \ycoordinate);
			\coordinate (\i RNW) at (\nwires - \x - 0.5, \ycoordinate + 1);
			\coordinate (\i RNE) at (\nwires - \x + 0.5, \ycoordinate + 1);

			\coordinate (\i LC) at (-\nwires + \x, \ycoordinate+ 0.5);
			\coordinate (\i RC) at (\nwires - \x, \ycoordinate+0.5);
		
			\draw (-\nwires+\x-0.5,\ycoordinate) -- (\i LSW)-- (\i LNE) -- (-\nwires+\x+0.5, \ycoordinate+1);
			\draw (-\nwires+\x+0.5,\ycoordinate) --(\i LSE) --(\i LNW) --(-\nwires+\x-0.5, \ycoordinate+1);

			\draw (\nwires - \x -0.5, \ycoordinate) -- (\i RSW)-- (\i RNE) -- (\nwires - \x + 0.5, \ycoordinate + 1);
			\draw (\nwires - \x +0.5, \ycoordinate) -- (\i RSE)-- (\i RNW) -- (\nwires - \x -0.5, \ycoordinate + 1);
			}{}

		\ifthenelse{\x=\y \OR \x = \yplusone}{}{
			\draw
			(-\nwires+\x-0.5,\ycoordinate)--(-\nwires+\x-0.5,\ycoordinate+1);
			\draw
			(\nwires-\x+0.5,\ycoordinate)--(\nwires-\x+0.5,\ycoordinate+1);
		}}

	\ifthenelse{\y=\nwires}{
		\node[above] at (\i LC) {$t_{\i}$};}
	{\node[above] at (\i LC) {$\bar{t}_{\i}$};
	\node[above] at (\i RC) {$t_{\i}$};}

}

\foreach \x in{1,...,\nwires}{
	\draw (\x-0.5,\length+1)--(\x-0.5,\length+1.5);
	\draw (\x-0.5,1)--(\x-0.5,0.5);

	\draw (-\x+0.5,\length+1)--(-\x+0.5,\length+1.5);
	\draw (-\x+0.5,1)--(-\x+0.5,0.5);

}

\draw[teal!50!green, dashed] (0,0.5)--(0,\length+1.5);

\draw[color=red!50, line cap=round, line width=3, opacity=0.5, ->]
(0.5,0.5)--(4LSE)--(4LC)--(4LSW)--(-0.5,0.5);

\end{tikzpicture}
\caption{Rigorous paths for ${\bm i} = (1,2,1,2) \in R(w_0^{(C_2)})$.}\label{figure_1212}
\end{figure}
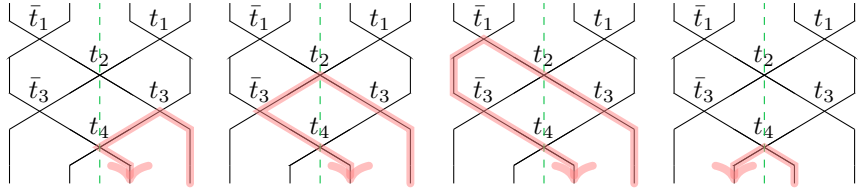

\begin{example}\label{example_1212}
Consider $\bm i = (1,2,1,2) \in \Rcn{2}$. 
Then there are four rigorous paths as shown in Figure~\ref{figure_1212} that all provide non-redundant inequalities. 
Indeed, they satisfy the assumption of Proposition~\ref{prop_non_redundancy_k_not_n_maximal}.
Accordingly, the number of facets of the string cone $\mathcal C_{\bm i}^{(C_2)}$ is $4$.
\end{example}

\begin{proposition}\label{prop_non_redundancy_n}
Let $\bm i \in R(w_0^{(C_n)})$, and $P \in \Paths^{\symp}(\bm i,n)$.
Then the inequality $\F_P^{(C)}({\bm a}) \geq 0$ is non-redundant in the expression \eqref{eq:descriptions_of_string_cones_in_type_C_n} if and only if $P$ is symmetric. 
\end{proposition}

\begin{figure}[t]
\begin{subfigure}[b]{0.35\textwidth}
\begin{tikzpicture}[scale = 0.8, yscale=0.6]
\def\nwires{3} 
\def\length{9} 

\foreach \i\y in {1/1,2/3,3/2,4/1,5/3,6/2,7/1,8/3,9/2}{
	\pgfmathsetmacro\ycoordinate{\length-\i+1}
	\pgfmathsetmacro\yplusone{\y+1}
	\foreach \x in {1,...,\nwires}{
		\ifthenelse{\x = \y}
			{\coordinate (\i LSW) at (-\nwires+\x-0.5, \ycoordinate);
			\coordinate (\i LSE) at (-\nwires+\x +0.5, \ycoordinate);
			\coordinate (\i LNW) at (-\nwires+\x-0.5, \ycoordinate+1);
			\coordinate (\i LNE) at (-\nwires+\x +0.5,\ycoordinate+1);

			\coordinate (\i RSW) at (\nwires - \x -0.5, \ycoordinate );
			\coordinate (\i RSE) at (\nwires - \x + 0.5, \ycoordinate );
			\coordinate (\i RNW) at (\nwires - \x - 0.5, \ycoordinate +1);
			\coordinate (\i RNE) at (\nwires - \x + 0.5, \ycoordinate +1);

			\coordinate (\i LC) at (-\nwires + \x, \ycoordinate+ 0.5);
			\coordinate (\i RC) at (\nwires - \x, \ycoordinate+0.5);
		
			\draw (-\nwires+\x-0.5,\ycoordinate) -- (\i LSW)-- (\i LNE) -- (-\nwires+\x+0.5, \ycoordinate+1);
			\draw (-\nwires+\x+0.5,\ycoordinate) --(\i LSE) --(\i LNW) --(-\nwires+\x-0.5, \ycoordinate+1);

			\draw (\nwires - \x -0.5, \ycoordinate) -- (\i RSW)-- (\i RNE) -- (\nwires - \x + 0.5, \ycoordinate + 1);
			\draw (\nwires - \x +0.5, \ycoordinate) -- (\i RSE)-- (\i RNW) -- (\nwires - \x -0.5, \ycoordinate + 1);
			}{}

		\ifthenelse{\x=\y \OR \x = \yplusone}{}{
			\draw
			(-\nwires+\x-0.5,\ycoordinate)--(-\nwires+\x-0.5,\ycoordinate+1);
			\draw
			(\nwires-\x+0.5,\ycoordinate)--(\nwires-\x+0.5,\ycoordinate+1);
		}}

	\ifthenelse{\y=\nwires}{
		\node[above] at (\i LC) {$t_{\i}$};}
	{\node[above] at (\i LC) {$\bar{t}_{\i}$};
	\node[above] at (\i RC) {$t_{\i}$};}

}

\foreach \x in{1,...,\nwires}{
	\draw (\x-0.5,\length+1)--(\x-0.5,\length+1.5);
	\draw (\x-0.5,1)--(\x-0.5,0.5);

	\draw (-\x+0.5,\length+1)--(-\x+0.5,\length+1.5);
	\draw (-\x+0.5,1)--(-\x+0.5,0.5);
}

\foreach \x in {1,2,3}{
\node[above] at (-3.5+\x,10.5) {$\ell_{\x}$};
\node[above] at (3.5-\x, 10.5){$\ell_{\overline{\x}}$};
}

\draw[teal!50!green, dashed] (0,0.5)--(0,\length+1.5);

\draw[color=red!50, line cap=round, line width=3, opacity=0.5, ->] 
(0.5,0.5)--(9RSW)--(9RNE)--(7RSW)--(7RNE)--(4RSE)--(4RNW)--(3RSE)--(3RC)--(3RSW)--(5RNE)--(5RSW)--(6LNE)--(6LC)--(6LNW)--(4LSE)--(4LC)--(4LSW)--(7LNW)--(7LSE)--(9LNW)--(9LSE)--(-0.5,0.5);

\draw[color=blue, line cap = round, line width = 1.5, dashed, ->] 
(0.5,0.5)--(9RSW)--(9RNE)--(7RSW)--(7RNE)--(4RSE)--(4RC)--(4RSW)--(6RNE)--(6RC)--(6RNW)--(5RSE)--(5RNW)--(3LSE)--(3LC)--(3LSW)--(4LNE)--(4LSW)--(7LNW)--(7LSE)--(9LNW)--(9LSE)--(-0.5,0.5);

\node[circle, draw, red] at (5RC) {};

\end{tikzpicture}
\caption{Non-symmetric paths ${P}$ and ${P}^{\vee}$.}
\label{figure_non_symmetric_paths}
\end{subfigure}
\begin{subfigure}[b]{0.6\textwidth}
\begin{tikzpicture}[scale = 0.8, yscale=0.6]
\def\nwires{3} 
\def\length{9} 

\foreach \i\y in {1/1,2/3,3/2,4/1,5/3,6/2,7/1,8/3,9/2}{
	\pgfmathsetmacro\ycoordinate{\length-\i+1}
	\pgfmathsetmacro\yplusone{\y+1}
	\foreach \x in {1,...,\nwires}{
		\ifthenelse{\x = \y}
			{\coordinate (\i LSW) at (-\nwires+\x-0.5, \ycoordinate);
			\coordinate (\i LSE) at (-\nwires+\x +0.5, \ycoordinate);
			\coordinate (\i LNW) at (-\nwires+\x-0.5, \ycoordinate+1);
			\coordinate (\i LNE) at (-\nwires+\x +0.5,\ycoordinate+1);

			\coordinate (\i RSW) at (\nwires - \x -0.5, \ycoordinate );
			\coordinate (\i RSE) at (\nwires - \x + 0.5, \ycoordinate );
			\coordinate (\i RNW) at (\nwires - \x - 0.5, \ycoordinate +1);
			\coordinate (\i RNE) at (\nwires - \x + 0.5, \ycoordinate +1);

			\coordinate (\i LC) at (-\nwires + \x, \ycoordinate+ 0.5);
			\coordinate (\i RC) at (\nwires - \x, \ycoordinate+0.5);
		
			\draw (-\nwires+\x-0.5,\ycoordinate) -- (\i LSW)-- (\i LNE) -- (-\nwires+\x+0.5, \ycoordinate+1);
			\draw (-\nwires+\x+0.5,\ycoordinate) --(\i LSE) --(\i LNW) --(-\nwires+\x-0.5, \ycoordinate+1);

			\draw (\nwires - \x -0.5, \ycoordinate) -- (\i RSW)-- (\i RNE) -- (\nwires - \x + 0.5, \ycoordinate + 1);
			\draw (\nwires - \x +0.5, \ycoordinate) -- (\i RSE)-- (\i RNW) -- (\nwires - \x -0.5, \ycoordinate + 1);
			}{}

		\ifthenelse{\x=\y \OR \x = \yplusone}{}{
			\draw
			(-\nwires+\x-0.5,\ycoordinate)--(-\nwires+\x-0.5,\ycoordinate+1);
			\draw
			(\nwires-\x+0.5,\ycoordinate)--(\nwires-\x+0.5,\ycoordinate+1);
		}}

	\ifthenelse{\y=\nwires}{
		\node[above] at (\i LC) {$t_{\i}$};}
	{\node[above] at (\i LC) {$\bar{t}_{\i}$};
	\node[above] at (\i RC) {$t_{\i}$};}

}

\foreach \x in{1,...,\nwires}{
	\draw (\x-0.5,\length+1)--(\x-0.5,\length+1.5);
	\draw (\x-0.5,1)--(\x-0.5,0.5);

	\draw (-\x+0.5,\length+1)--(-\x+0.5,\length+1.5);
	\draw (-\x+0.5,1)--(-\x+0.5,0.5);
}

\foreach \x in {1,2,3}{
\node[above] at (-3.5+\x,10.5) {$\ell_{\x}$};
\node[above] at (3.5-\x, 10.5){$\ell_{\overline{\x}}$};
}

\draw[teal!50!green, dashed] (0,0.5)--(0,\length+1.5);

\draw[color = red!70!blue, line cap=round, line width=3, opacity=0.5, ->] 
(0.5,0.5)--(9RSW)--(9RNE)--(7RSW)--(7RNE)--(4RSE)--(4RNW)--(3RSE)--(3RC)--(3RSW)--(5RNE)--(5RC)--(5RNW)--(3LSE)--(3LC)--(3LSW)--(4LNE)--(4LSW)--(7LNW)--(7LSE)--(9LNW)--(9LSE)--(-0.5,0.5);

\end{tikzpicture}%
\begin{tikzpicture}[scale = 0.8, yscale=0.6]
\def\nwires{3} 
\def\length{9} 

\foreach \i\y in {1/1,2/3,3/2,4/1,5/3,6/2,7/1,8/3,9/2}{
	\pgfmathsetmacro\ycoordinate{\length-\i+1}
	\pgfmathsetmacro\yplusone{\y+1}
	\foreach \x in {1,...,\nwires}{
		\ifthenelse{\x = \y}
			{\coordinate (\i LSW) at (-\nwires+\x-0.5, \ycoordinate);
			\coordinate (\i LSE) at (-\nwires+\x +0.5, \ycoordinate);
			\coordinate (\i LNW) at (-\nwires+\x-0.5, \ycoordinate+1);
			\coordinate (\i LNE) at (-\nwires+\x +0.5,\ycoordinate+1);

			\coordinate (\i RSW) at (\nwires - \x -0.5, \ycoordinate );
			\coordinate (\i RSE) at (\nwires - \x + 0.5, \ycoordinate );
			\coordinate (\i RNW) at (\nwires - \x - 0.5, \ycoordinate +1);
			\coordinate (\i RNE) at (\nwires - \x + 0.5, \ycoordinate +1);

			\coordinate (\i LC) at (-\nwires + \x, \ycoordinate+ 0.5);
			\coordinate (\i RC) at (\nwires - \x, \ycoordinate+0.5);
		
			\draw (-\nwires+\x-0.5,\ycoordinate) -- (\i LSW)-- (\i LNE) -- (-\nwires+\x+0.5, \ycoordinate+1);
			\draw (-\nwires+\x+0.5,\ycoordinate) --(\i LSE) --(\i LNW) --(-\nwires+\x-0.5, \ycoordinate+1);

			\draw (\nwires - \x -0.5, \ycoordinate) -- (\i RSW)-- (\i RNE) -- (\nwires - \x + 0.5, \ycoordinate + 1);
			\draw (\nwires - \x +0.5, \ycoordinate) -- (\i RSE)-- (\i RNW) -- (\nwires - \x -0.5, \ycoordinate + 1);
			}{}

		\ifthenelse{\x=\y \OR \x = \yplusone}{}{
			\draw
			(-\nwires+\x-0.5,\ycoordinate)--(-\nwires+\x-0.5,\ycoordinate+1);
			\draw
			(\nwires-\x+0.5,\ycoordinate)--(\nwires-\x+0.5,\ycoordinate+1);
		}}

	\ifthenelse{\y=\nwires}{
		\node[above] at (\i LC) {$t_{\i}$};}
	{\node[above] at (\i LC) {$\bar{t}_{\i}$};
	\node[above] at (\i RC) {$t_{\i}$};}

}

\foreach \x in{1,...,\nwires}{
	\draw (\x-0.5,\length+1)--(\x-0.5,\length+1.5);
	\draw (\x-0.5,1)--(\x-0.5,0.5);

	\draw (-\x+0.5,\length+1)--(-\x+0.5,\length+1.5);
	\draw (-\x+0.5,1)--(-\x+0.5,0.5);
}

\foreach \x in {1,2,3}{
\node[above] at (-3.5+\x,10.5) {$\ell_{\x}$};
\node[above] at (3.5-\x, 10.5){$\ell_{\overline{\x}}$};
}

\draw[teal!50!green, dashed] (0,0.5)--(0,\length+1.5);

\draw[color = red!30!blue, line cap=round, line width=3, opacity=0.4, ->] 
(0.5,0.5)--(9RSW)--(9RNE)--(7RSW)--(7RNE)--(4RSE)--(4RC)--(4RSW)--(6RNE)--(6RC)--(6RNW)--(5RSE)--(5RC)--(5RSW)--(6LNE)--(6LC)--(6LNW)--(4LSE)--(4LC)--(4LSW)--(7LNW)--(7LSE)--(9LNW)--(9LSE)--(-0.5,0.5);

\end{tikzpicture}
\caption{Symmetric paths $Q_1$ and $Q_2$.}
\end{subfigure}
\caption{A non-symmetric path $P = (\ell_3 \to \ell_{\bar{1}} \to \ell_2 \to \ell_{\bar{3}})$ (red highlighted), its mirror ${P}^{\vee} = (\ell_3 \to \ell_{\bar{2}} \to \ell_1 \to \ell_{\bar{3}})$ (blue dashed), and the corresponding symmetric paths $Q_1, Q_2$ for $\bm i = (1,3,2,1,3,2,1,3,2) \in \Rcn{3}$.}
\label{figure_non_symmetric_paths_and_making_symmetric_paths}
\end{figure}
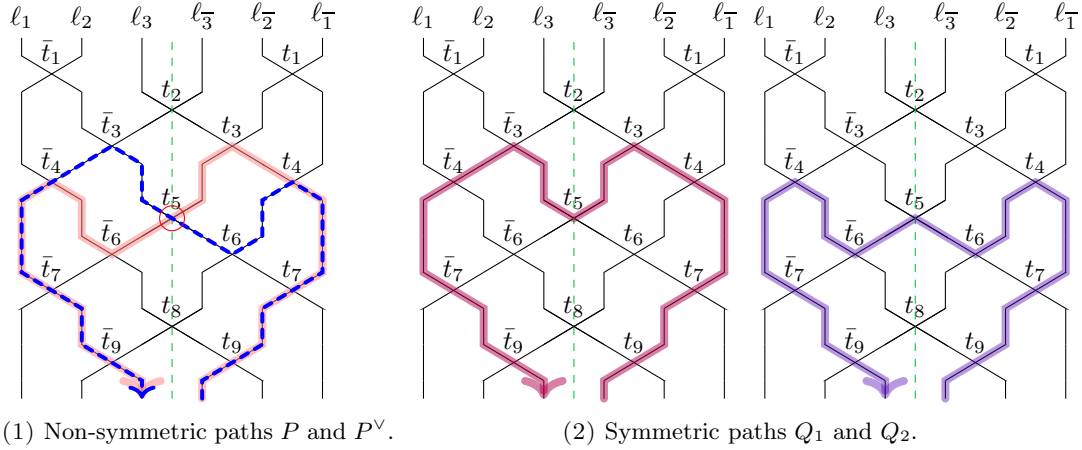

\begin{proof}
If $P$ is symmetric, then the assumption of Proposition \ref{prop_non_redundancy_k_not_n_maximal} is satisfied since $\chamber(P) \cup \chamber(P^\vee) = \chamber(P)$. 
Hence, by Proposition \ref{prop_non_redundancy_k_not_n_maximal}, it suffices to prove that if $P$ is not symmetric, then the inequality $\F_P^{(C)}({\bm a}) \geq 0$ is redundant. 

Assume that $P$ is not symmetric, and write $P = (\ell_{n} \to \ell_{r_1} \to \dots \to \ell_{r_s} \to \ell_{n+1})$. 
Then its mirror is given by ${P}^{\vee} = (\ell_{n} \to \ell_{\overline{r_s}} \to \cdots \to \ell_{\overline{r_1}} \to \ell_{\overline{n}})$. 
Because both paths start at $L_n$ and end at~$L_{\overline{n}}$, each of ${P}$ and ${P}^{\vee}$ meets the wall, which is the central vertical line, once at the same place, say $\ell_{r_p} \cap \ell_{\overline{r_p}}$. Define two \emph{symmetric} paths $Q_1$ and $Q_2$ by
\begin{equation}\label{eq:extended_path_symmetric_case}
\begin{aligned}
Q_1 &\colonequals (\ell_n \to \ell_{r_1} \to \dots \to \ell_{r_p} 
\to \ell_{\overline{r_p}} \to \dots \to \ell_{\overline{r_1}} \to \ell_{\overline{n}}),\\
Q_2 &\colonequals (\ell_n \to \ell_{\overline{r_s}} \to \dots \to \ell_{\overline{r_p}} 
\to \ell_{r_p} \to \dots \to \ell_{r_s}\to \ell_{\overline{n}}).
\end{aligned}
\end{equation}
See Figure~\ref{figure_non_symmetric_paths_and_making_symmetric_paths} and Example~\ref{example_Q1_and_Q2}. 
Then we obtain
\[
\hat{\F}_{{P}}^{(C)}({\bm a}) + \hat{\F}_{{P}^{\vee}}^{(C)}({\bm a}) = \hat{\F}_{Q_1}^{(C)}({\bm a}) + \hat{\F}_{Q_2}^{(C)}({\bm a}),
\]
which implies $\F_{P}^{(C)}({\bm a}) = \F_{Q_1}^{(C)}({\bm a}) + \F_{Q_2}^{(C)}({\bm a})$. 
This proves the proposition.
\end{proof}

\begin{example}\label{example_Q1_and_Q2}
Let $\bm i = (1,3,2,1,3,2,1,3,2) \in R(w_0^{(C_3)})$. Take $P = (\ell_3 \to \ell_{\bar{1}} \to \ell_2 \to \ell_{\bar{3}})$. 
Then its dual ${P}^{\vee}$ is 
$P^{\vee} = (\ell_3 \to \ell_{\bar{2}} \to \ell_1 \to \ell_{\bar{3}})$.
They intersect in the wall at the intersection $\ell_1 \cap \ell_{\bar{1}}$ (see the red circle in Figure~\ref{figure_non_symmetric_paths}). Accordingly, we obtain two symmetric paths $Q_1$ and $Q_2$:
\[
Q_1 = (\ell_3 \to \ell_{\bar{1}} \to \ell_1 \to \ell_{\bar{3}}), \quad 
Q_2 = (\ell_3 \to \ell_{\bar{2}} \to \ell_1 \to \ell_{\bar{1}} \to \ell_2 \to \ell_{\bar{3}}).
\]
Because of $\chamber(\hat{P}) \cup \chamber(\hat{P}^{\vee}) = \chamber(Q_1) \cup \chamber(Q_2)$, we obtain $\hat{\F}_{\hat{P}} + \hat{\F}_{\hat{P}^{\vee}} = \hat{\F}_{\hat{Q}_1}  + \hat{\F}_{\hat{Q}_2}$. 
Indeed, the functions $\hat{\F}_{\mathscr{P}}$, $\hat{\F}_{\mathscr{P}}^{(C)}$, and $\F_{\mathscr{P}}^{(C)}$ for $\mathscr{P} \in \{P,P^{\vee},Q_1,Q_2\}$ are given as follows. 
\begin{center}
\begin{tabular}{l|lll}
\toprule
$\mathscr{P}$ & $\hat{\F}_{\mathscr P}(\bm t)$ & $\hat{\F}_{\hat{\mathscr P}}^{(C)}(\bm a)$ & $\F_{\mathscr P}^{(C)}(\bm a)$  \\
\midrule
$P$ & $t_3 + \bar{t}_4 - \bar{t}_6$ & $a_3 + a_4 - a_6$  
	& $a_3 + a_4 - a_6$ \\
$P^{\vee}$ & $\bar{t}_3 + t_4 - t_6$ & $a_3 + a_4 - a_6$  
	& $a_3 + a_4 - a_6$\\
$Q_1$ & $t_3 + \bar{t}_3 - t_5$& $2a_3 - 2a_5$ 
	& $a_3 - a_5$ \\ 
$Q_2$ & $t_4 + \bar{t}_4 + t_5 - (t_6 + \bar{t}_6)$ & $2a_4 + 2a_5 - 2a_6$ 
	& $a_4 + a_5 - a_6$ \\ 
\bottomrule
\end{tabular}
\end{center}
This provides $\F_{P}^{(C)}({\bm a}) = a_3 + a_4 - a_6 = \F_{Q_1}^{(C)}({\bm a}) + \F_{Q_2}^{(C)}({\bm a})$.
Accordingly, an inequality given by a non-symmetric path can be expressed by the sum of inequalities given by symmetric paths.
\end{example}

\begin{definition}\label{d:extension_k_equal_n}
For $P_1, P_2 \in \Paths^{\symp}(\bm i)$, we say that $P_1$ is an \emph{extension} of $P_2$ if $\chamber(P_2) \subseteq \chamber(P_1)$.
For $P \in \Paths^{\symp}(\bm i,n)$, we define an extension $P_{\rm ex} \in \Paths^{\symp}(\bm i,n)$ as follows.
If $P$ is symmetric, then we set $P_{\rm ex} \coloneqq P$.
Otherwise, $P_{\rm ex}$ is defined to be the symmetric path in $\Paths^{\symp}(\bm i,n)$ such that $\chamber(P_{\rm ex}) = \chamber(P) \cup \chamber(P^\vee)$, which coincides with $Q_1$ or $Q_2$ in \eqref{eq:extended_path_symmetric_case}.
\end{definition}

Since $P_{\rm ex}$ is symmetric, the inequality $\F_{P_{\rm ex}}^{(C)}({\bm a}) \geq 0$ is non-redundant in the expression \eqref{eq:descriptions_of_string_cones_in_type_C_n}.
Let $k \neq n$. 
If $P \in \Paths^{\symp}(\bm i,k)$ does not satisfy the assumption of Proposition \ref{prop_non_redundancy_k_not_n_maximal}, then we can extend it in $\chamber(P) \cup \chamber(P^\vee)$ so that the extended rigorous path $Q$ satisfies the assumption of Proposition~\ref{prop_non_redundancy_k_not_n_maximal}. 
More precisely, we obtain the following. 

\begin{proposition}\label{proposition_extension_exists}
Let $k \neq n$, and $P \in \Paths^{\symp}(\bm i,k)$.
Then there exists a unique extension $P_{\rm ex} \in \Paths^{\symp}(\bm i,k)$ such that $\chamber(P_{\rm ex}) \subseteq \chamber(P) \cup \chamber(P^\vee)$ and such that $P_{\rm ex}$ satisfies the assumption of Proposition \ref{prop_non_redundancy_k_not_n_maximal}.
\end{proposition}
\begin{remark}
	For $k = n$, the extension $P_{\rm ex}$ is symmetric. Hence it provides non-redundant inequality. For $k \neq n$, even though the extension $P_{\rm ex}$ is not symmetric, it provides non-redundant inequality by Proposition \ref{prop_non_redundancy_k_not_n_maximal}.
\end{remark}

\begin{remark}\label{r:maximality_of_extension}
Since $\chamber(P) \subseteq \chamber(P_{\rm ex}) \subseteq \chamber(P) \cup \chamber(P^\vee)$, we have $\chamber(P_{\rm ex}) \cup \chamber(P_{\rm ex}^\vee) = \chamber(P) \cup \chamber(P^\vee)$. 
Hence, for all $Q \in \Paths^{\symp}(\bm i,k)$ such that $\chamber(Q) \subseteq \chamber(P) \cup \chamber(P^\vee)$, it follows that $\chamber(Q) \subseteq \chamber(P_{\rm ex})$ since $P_{\rm ex}$ satisfies the assumption of Proposition \ref{prop_non_redundancy_k_not_n_maximal}. 
\end{remark}

\begin{proof}[Proof of Proposition~\ref{proposition_extension_exists}]
If $P$ satisfies the assumption of Proposition \ref{prop_non_redundancy_k_not_n_maximal}, then we set $P_{\rm ex} \coloneqq P$. 
If $P$ does not satisfy the assumption of Proposition \ref{prop_non_redundancy_k_not_n_maximal}, then there exists $Q \in \Paths^{\symp}(\bm i,k)$ such that $\chamber(Q) \subseteq \chamber(P) \cup \chamber(P^\vee)$ and such that $\chamber(Q) \nsubseteq \chamber(P)$. 
We write all the nodes on the paths $P$ and~$Q$ as 
\[
\begin{split}
P &= (L_k = d_0 \to d_1 \to \cdots \to d_r \to d_{r+1} = L_{k+1}); \\
Q &= (L_k = d_0^\prime \to d_1^\prime \to \cdots \to d_s^\prime \to d_{s+1}^\prime = L_{k+1}),
\end{split}
\]
respectively, where $d_1, \ldots, d_r, d_1^\prime, \ldots, d_s^\prime \in \{t_j \mid 1 \leq j \leq N\} \cup \{\bar{t}_j \mid 1 \leq j \leq N,\ i_j \neq n\}$.
Note that $d_1 = d_1^\prime$ and $d_r = d_s^\prime$ by the shape of the symplectic wiring diagram $G^{\symp}({\bm i}, k)$. 
Let $u_1$ be the minimum of $1 \leq u \leq s-1$ such that the fragment between $d_u^\prime$ and $d_{u+1}^\prime$ is not included in the interior of $\chamber(P)$ and the boundary of $\chamber(P)$. 
Since $\chamber(Q) \nsubseteq \chamber(P)$, such $u$ always exists. 
In addition, denote by $u_2$ the minimum of $u_1 +1 \leq u \leq s$ such that the fragment between $d_u^\prime$ and $d_{u+1}^\prime$ is included in the interior of $\chamber(P)$ and the boundary of $\chamber(P)$. 
Since $d_r = d_s^\prime$, such $u$ always exists. 
Define $1 \leq v_1 < v_2 \leq r$ by $d_{v_1} = d_{u_1}^\prime$ and $d_{v_2} = d_{u_2}^\prime$, respectively.
Now we claim that the following path is a rigorous path in $\Paths^{\symp}(\bm i,k)$: 
\begin{equation}\label{eq:extended_rigorous_path}
\begin{aligned}
(L_k \to d_1 \to \cdots \to d_{v_1} (= d_{u_1}^\prime) \to d_{u_1 +1}^\prime \to \cdots \to d_{u_2}^\prime (= d_{v_2}) \to d_{v_2 +1} \to \cdots \to d_r \to L_{k+1}).
\end{aligned}
\end{equation}
Here, we write all the nodes on the path in~\eqref{eq:extended_rigorous_path}.
It suffices to show that the fragments $(d_{v_1 -1} \to d_{v_1} \to d_{u_1 +1}^\prime)$ and $(d_{u_2 -1}^\prime \to d_{u_2}^\prime \to d_{v_2 +1})$ do not give a forbidden fragment.
We prove only the assertion for $(d_{v_1 -1} \to d_{v_1} \to d_{u_1 +1}^\prime)$ since the proof for $(d_{u_2 -1}^\prime \to d_{u_2}^\prime \to d_{v_2 +1})$ is similar. There are two possibilities: $d_{v_1 -1} = d_{u_1 -1}^\prime$ or $d_{v_1 -1} \neq d_{u_1 -1}^\prime$; we consider two cases separately.
\smallskip 

\noindent \textsf{Case 1. $d_{v_1 -1} = d_{u_1 -1}^\prime$}. If $d_{v_1 -1} = d_{u_1 -1}^\prime$, then the assertion is obvious since the fragment $(d_{v_1 -1} \to d_{v_1} \to d_{u_1 +1}^\prime) = (d_{u_1 -1}^\prime \to d_{u_1}^\prime \to d_{u_1 +1}^\prime)$ is a part of the rigorous path~$Q$. 

\smallskip 

\noindent \textsf{Case 2. $d_{v_1 -1} \neq d_{u_1 -1}^\prime$}.
Suppose that $d_{v_1 -1} \neq d_{u_1 -1}^\prime$. 
Then the fragments $(d_{v_1 -1} \to d_{v_1} \to d_{v_1 +1})$ and $(d_{u_1 -1}^\prime \to d_{u_1}^\prime \to d_{u_1 +1}^\prime)$ are given as in Figure~\ref{fig_fragments}. Indeed, when we write down all the possibilities of the fragments $(d_{v_1 -1} \to d_{v_1} \to d_{v_1 +1})$ and $(d_{u_1 -1}^\prime \to d_{u_1}^\prime \to d_{u_1 +1}^\prime)$, we see that the other possibilities are impossible by the following properties: 
	\begin{itemize}
		\item $P$ and $Q$ do not contain a forbidden fragment,
		\item the fragment between $d_{u_1 -1}^\prime$ and $d_{u_1}^\prime$ is included in the interior of $\chamber(P)$ or in the boundary of $\chamber(P)$ because of the minimality of $u_1$, and
		\item the interior of $\chamber(P)$ is always on the left hand side of the path $P$ because $P$ passes through each node at most once.
	\end{itemize}
\begin{figure}
\begin{tikzpicture}[scale = 0.6]
\draw[-{>[sep=10pt]_}, red] (1,-1) node[below right] {$d_{u_1-1}'$} to (-1,1) node[above left] {$d_{u_1+1}'$};
\draw[-{>[sep=10pt]_}, blue] (1,1) node[above right] {$d_{v_1-1}$} to (-1,-1) node[below left] {$d_{v_1+1}$};
\node[right, blue] at (0,0) {$d_{v_1}$};
\end{tikzpicture} 
\caption{The fragments $(d_{v_1 -1} \to d_{v_1} \to d_{v_1 +1})$ and $(d_{u_1 -1}^\prime \to d_{u_1}^\prime \to d_{u_1 +1}^\prime)$ in the proof of Proposition~\ref{proposition_extension_exists}}\label{fig_fragments}
\end{figure}
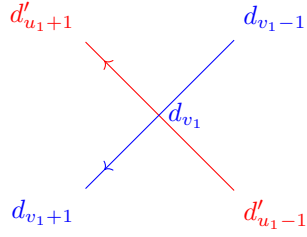
\noindent This implies that the fragment $(d_{v_1 -1} \to d_{v_1} \to d_{u_1 +1}^\prime)$ does not give a forbidden fragment.
Hence our claim follows.

\smallskip 

Repeating this argument by replacing $P$ with the extended path \eqref{eq:extended_rigorous_path}, we deduce an existence of $P_{\rm ex}$.
The uniqueness of $P_{\rm ex}$ follows immediately since $P_{\rm ex}$ satisfies the assumption of Proposition~\ref{prop_non_redundancy_k_not_n_maximal} (see Remarks \ref{rmk_maximal_Q_is_P} and \ref{r:maximality_of_extension}).
\end{proof}

\begin{example}
Let $\bm i = (2,1,3,2,1,3,2,1,3) \in \Rcn{3}$. 
Consider $P = (\ell_2 \to \ell_{\bar{1}} \to \ell_3)$ (see Figure~\ref{fig_extension_example_1}). The path $P$ does not satisfy the assumption of Proposition~\ref{prop_non_redundancy_k_not_n_maximal} because the path $Q = (\ell_2 \to \ell_{\bar{3}} \to \ell_1 \to \ell_{\bar{2}} \to \ell_3)$ satisfies $\chamber(Q) \subseteq \chamber(P) \cup \chamber(P^{\vee})$ and $\chamber(Q) \not\subseteq \chamber(P)$ (see Figure~\ref{fig_extension_example_2}). Following Proposition~\ref{proposition_extension_exists}, there uniquely exists an extension $P_{\rm ex} = (\ell_2 \to \ell_{\bar{1}} \to \ell_1 \to \ell_{\bar{2}} \to \ell_3)$ as depicted in Figure~\ref{fig_extension_example_3}. 
Moreover, we obtain $Q_{\rm ex} = P_{\rm ex}$.
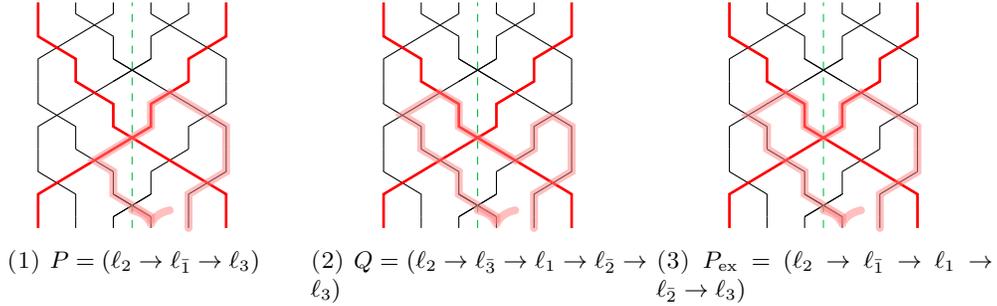
\begin{figure}
\begin{subfigure}[t]{0.3\textwidth}
\centering
\begin{tikzpicture}[scale = 0.5, yscale=0.6]
\def\nwires{3} 
\def\length{9} 

\foreach \i\y in {1/2,2/1,3/3,4/2,5/1,6/3,7/2,8/1,9/3}{
	\pgfmathsetmacro\ycoordinate{\length-\i+1}
	\pgfmathsetmacro\yplusone{\y+1}
	\foreach \x in {1,...,\nwires}{
		\ifthenelse{\x = \y}
			{\coordinate (\i LSW) at (-\nwires+\x-0.5, \ycoordinate);
			\coordinate (\i LSE) at (-\nwires+\x +0.5, \ycoordinate);
			\coordinate (\i LNW) at (-\nwires+\x-0.5, \ycoordinate+1);
			\coordinate (\i LNE) at (-\nwires+\x +0.5,\ycoordinate+1);

			\coordinate (\i RSW) at (\nwires - \x -0.5, \ycoordinate );
			\coordinate (\i RSE) at (\nwires - \x + 0.5, \ycoordinate );
			\coordinate (\i RNW) at (\nwires - \x - 0.5, \ycoordinate +1);
			\coordinate (\i RNE) at (\nwires - \x + 0.5, \ycoordinate +1);

			\coordinate (\i LC) at (-\nwires + \x, \ycoordinate+ 0.5);
			\coordinate (\i RC) at (\nwires - \x, \ycoordinate+0.5);
		
			\draw (-\nwires+\x-0.5,\ycoordinate) -- (\i LSW)-- (\i LNE) -- (-\nwires+\x+0.5, \ycoordinate+1);
			\draw (-\nwires+\x+0.5,\ycoordinate) --(\i LSE) --(\i LNW) --(-\nwires+\x-0.5, \ycoordinate+1);

			\draw (\nwires - \x -0.5, \ycoordinate) -- (\i RSW)-- (\i RNE) -- (\nwires - \x + 0.5, \ycoordinate + 1);
			\draw (\nwires - \x +0.5, \ycoordinate) -- (\i RSE)-- (\i RNW) -- (\nwires - \x -0.5, \ycoordinate + 1);
			}{}

		\ifthenelse{\x=\y \OR \x = \yplusone}{}{
			\draw
			(-\nwires+\x-0.5,\ycoordinate)--(-\nwires+\x-0.5,\ycoordinate+1);
			\draw
			(\nwires-\x+0.5,\ycoordinate)--(\nwires-\x+0.5,\ycoordinate+1);
		}
}


}

\foreach \x in{1,...,\nwires}{
	\draw (\x-0.5,\length+1)--(\x-0.5,\length+1.5);
	\draw (\x-0.5,1)--(\x-0.5,0.5);

	\draw (-\x+0.5,\length+1)--(-\x+0.5,\length+1.5);
	\draw (-\x+0.5,1)--(-\x+0.5,0.5);
}
\draw[color=red, line width = 1] 
(2.5,0.5)--(8RSE)--(6RNW)--(4LSE)--(4LNW)--(2LSE)--(2LNW)--(-2.5,10.5)
(-2.5,0.5)--(8LSW)--(6LNE)--(4RSW)--(4RNE)--(2RSW)--(2RNE)--(2.5,10.5);

\draw[teal!50!green, dashed] (0,0.5)--(0,\length+1.5);
\draw[color=red!50, line cap=round, line width=3, opacity=0.5, ->]
(1.5,0.5)--(8RSW)--(8RNE)--(5RSE)--(4RC)--(4RSW)--(6RNE)--(7LC)--(7LSE)--(9RNW)--(9RSE)--(0.5,0.5);
\end{tikzpicture}
\caption{$P = (\ell_2 \to \ell_{\bar{1}} \to \ell_3)$}\label{fig_extension_example_1}
\end{subfigure}
\begin{subfigure}[t]{0.3\textwidth}
\centering
\begin{tikzpicture}[scale = 0.5, yscale=0.6]
\def\nwires{3} 
\def\length{9} 

\foreach \i\y in {1/2,2/1,3/3,4/2,5/1,6/3,7/2,8/1,9/3}{
	\pgfmathsetmacro\ycoordinate{\length-\i+1}
	\pgfmathsetmacro\yplusone{\y+1}
	\foreach \x in {1,...,\nwires}{
		\ifthenelse{\x = \y}
			{\coordinate (\i LSW) at (-\nwires+\x-0.5, \ycoordinate);
			\coordinate (\i LSE) at (-\nwires+\x +0.5, \ycoordinate);
			\coordinate (\i LNW) at (-\nwires+\x-0.5, \ycoordinate+1);
			\coordinate (\i LNE) at (-\nwires+\x +0.5,\ycoordinate+1);

			\coordinate (\i RSW) at (\nwires - \x -0.5, \ycoordinate );
			\coordinate (\i RSE) at (\nwires - \x + 0.5, \ycoordinate );
			\coordinate (\i RNW) at (\nwires - \x - 0.5, \ycoordinate +1);
			\coordinate (\i RNE) at (\nwires - \x + 0.5, \ycoordinate +1);

			\coordinate (\i LC) at (-\nwires + \x, \ycoordinate+ 0.5);
			\coordinate (\i RC) at (\nwires - \x, \ycoordinate+0.5);
		
			\draw (-\nwires+\x-0.5,\ycoordinate) -- (\i LSW)-- (\i LNE) -- (-\nwires+\x+0.5, \ycoordinate+1);
			\draw (-\nwires+\x+0.5,\ycoordinate) --(\i LSE) --(\i LNW) --(-\nwires+\x-0.5, \ycoordinate+1);

			\draw (\nwires - \x -0.5, \ycoordinate) -- (\i RSW)-- (\i RNE) -- (\nwires - \x + 0.5, \ycoordinate + 1);
			\draw (\nwires - \x +0.5, \ycoordinate) -- (\i RSE)-- (\i RNW) -- (\nwires - \x -0.5, \ycoordinate + 1);
			}{}

		\ifthenelse{\x=\y \OR \x = \yplusone}{}{
			\draw
			(-\nwires+\x-0.5,\ycoordinate)--(-\nwires+\x-0.5,\ycoordinate+1);
			\draw
			(\nwires-\x+0.5,\ycoordinate)--(\nwires-\x+0.5,\ycoordinate+1);
		}
}


}

\foreach \x in{1,...,\nwires}{
	\draw (\x-0.5,\length+1)--(\x-0.5,\length+1.5);
	\draw (\x-0.5,1)--(\x-0.5,0.5);

	\draw (-\x+0.5,\length+1)--(-\x+0.5,\length+1.5);
	\draw (-\x+0.5,1)--(-\x+0.5,0.5);
}
\draw[color=red, line width = 1] 
(2.5,0.5)--(8RSE)--(6RNW)--(4LSE)--(4LNW)--(2LSE)--(2LNW)--(-2.5,10.5)
(-2.5,0.5)--(8LSW)--(6LNE)--(4RSW)--(4RNE)--(2RSW)--(2RNE)--(2.5,10.5);

\draw[teal!50!green, dashed] (0,0.5)--(0,\length+1.5);

\draw[color=red!50, line cap=round, line width=3, opacity=0.5, ->]
(1.5,0.5)--(8RSW)--(8RNE)--(5RSE)--(5RC)--(5RSW)--(7RNE)--(7RC)--(6RNW)--(4LSE)--(4LC)--(5LC)--(5LSE)--(7LNW)--(7LSE)--(9RNW)--(9RSE)--(0.5,0.5);
\end{tikzpicture}
\caption{$Q = (\ell_2 \to \ell_{\bar{3}} \to \ell_1 \to \ell_{\bar{2}} \to \ell_3)$}\label{fig_extension_example_2}
\end{subfigure}
\begin{subfigure}[t]{0.3\textwidth}
\centering
\begin{tikzpicture}[scale = 0.5, yscale=0.6]
\def\nwires{3} 
\def\length{9} 

\foreach \i\y in {1/2,2/1,3/3,4/2,5/1,6/3,7/2,8/1,9/3}{
	\pgfmathsetmacro\ycoordinate{\length-\i+1}
	\pgfmathsetmacro\yplusone{\y+1}
	\foreach \x in {1,...,\nwires}{
		\ifthenelse{\x = \y}
			{\coordinate (\i LSW) at (-\nwires+\x-0.5, \ycoordinate);
			\coordinate (\i LSE) at (-\nwires+\x +0.5, \ycoordinate);
			\coordinate (\i LNW) at (-\nwires+\x-0.5, \ycoordinate+1);
			\coordinate (\i LNE) at (-\nwires+\x +0.5,\ycoordinate+1);

			\coordinate (\i RSW) at (\nwires - \x -0.5, \ycoordinate );
			\coordinate (\i RSE) at (\nwires - \x + 0.5, \ycoordinate );
			\coordinate (\i RNW) at (\nwires - \x - 0.5, \ycoordinate +1);
			\coordinate (\i RNE) at (\nwires - \x + 0.5, \ycoordinate +1);

			\coordinate (\i LC) at (-\nwires + \x, \ycoordinate+ 0.5);
			\coordinate (\i RC) at (\nwires - \x, \ycoordinate+0.5);
		
			\draw (-\nwires+\x-0.5,\ycoordinate) -- (\i LSW)-- (\i LNE) -- (-\nwires+\x+0.5, \ycoordinate+1);
			\draw (-\nwires+\x+0.5,\ycoordinate) --(\i LSE) --(\i LNW) --(-\nwires+\x-0.5, \ycoordinate+1);

			\draw (\nwires - \x -0.5, \ycoordinate) -- (\i RSW)-- (\i RNE) -- (\nwires - \x + 0.5, \ycoordinate + 1);
			\draw (\nwires - \x +0.5, \ycoordinate) -- (\i RSE)-- (\i RNW) -- (\nwires - \x -0.5, \ycoordinate + 1);
			}{}

		\ifthenelse{\x=\y \OR \x = \yplusone}{}{
			\draw
			(-\nwires+\x-0.5,\ycoordinate)--(-\nwires+\x-0.5,\ycoordinate+1);
			\draw
			(\nwires-\x+0.5,\ycoordinate)--(\nwires-\x+0.5,\ycoordinate+1);
		}
}


}

\foreach \x in{1,...,\nwires}{
	\draw (\x-0.5,\length+1)--(\x-0.5,\length+1.5);
	\draw (\x-0.5,1)--(\x-0.5,0.5);

	\draw (-\x+0.5,\length+1)--(-\x+0.5,\length+1.5);
	\draw (-\x+0.5,1)--(-\x+0.5,0.5);
}
\draw[color=red, line width = 1] 
(2.5,0.5)--(8RSE)--(6RNW)--(4LSE)--(4LNW)--(2LSE)--(2LNW)--(-2.5,10.5)
(-2.5,0.5)--(8LSW)--(6LNE)--(4RSW)--(4RNE)--(2RSW)--(2RNE)--(2.5,10.5);

\draw[teal!50!green, dashed] (0,0.5)--(0,\length+1.5);
\draw[color=red!50, line cap=round, line width=3, opacity=0.5, ->]
(1.5,0.5)--(8RSW)--(8RNE)--(5RSE)--(4RC)--(4RSW)--(6RNE)--(6RC)--(6RNW)--(4LSE)--(4LC)--(5LC)--(5LSE)--(7LNW)--(7LSE)--(9LNW)--(9LSE)--(0.5,0.5);
\end{tikzpicture}
\caption{$P_{\rm ex} = (\ell_2 \to \ell_{\bar{1}} \to \ell_1 \to \ell_{\bar{2}} \to \ell_3)$}\label{fig_extension_example_3}
\end{subfigure}
\caption{A rigorous path in $\Paths^{\symp}(2,1,3,2,1,3,2,1,3)$ and its extension.}
\end{figure}
\end{example}

\section{Gelfand--Tsetlin type string polytopes in type $C$}
\label{sec_GCtype}

In this section, we prove Theorems \ref{thm_intro_2} and \ref{thm_intro_3}. 
Indeed, we use an inductive argument on the rank of the Lie algebra. 
In addition, we recall from~\cite{CKLP21_GC} combinatorial properties of rigorous paths in $\Paths(\hat{\bm i})$ for ${\bm i} \in \RCn$. To distinguish rigorous paths in $\Paths^{\symp}(\bm i)$ and that in~$\Paths(\hat{\bm i})$, we call rigorous paths in $\Paths^{\symp}(\bm i)$ \textit{symplectic} rigorous paths in this section.
	
	A {\em $\mcal{D}$-contraction} (resp., an {\em $\mcal{A}$-contraction}) is a process of producing a new reduced word in~$\Ran{n-1}$ 
	from a given element in $\Ran{n}$ (see \cite[Definition~3.6]{CKLP21_GC}). 
	More precisely, a $\mcal{D}$-contraction can be described in terms of a wiring diagram as follows. 
	For a given reduced word ${\bm i} \in \Ran{n}$, we may remove $\ell_{n+1}$ in the wiring diagram $G(\bm {i})$ and so we have a new 
	diagram consisting of $n$ wires $\ell_1, \dots, \ell_n$. 
	The induced diagram becomes a wiring diagram for some reduced decomposition  in $\Ran{n-1}$ since the $n$ wires still meet pairwise exactly once.
	We similarly define an $\mcal{A}$-contraction as a procedure obtaining a new wiring diagram from $G(\bm {i})$ by removing $\ell_1$. 
	We can also define a contraction in the case of type $C_n$ as follows. 

	\begin{definition} 
		For $n \geq 3$ and ${\bm i} \in \RCn$, we define ${\bf cont}({\bm i}) \in \Rcn{n-1}$ to be the reduced word corresponding to the symplectic wiring diagram obtained from 
		$G^{\symp}({\bm i})$ by removing both $\ell_1$ and $\ell_{\bar{1}}$. 
      Then the map 
		\[
			{\bf cont} \colon \RCn \rightarrow \Rcn{n-1},\quad {\bm i} \mapsto {\bf cont}({\bm i}),
		\]
		is called a {\em contraction}. 
	\end{definition}

We notice that for ${\bm i} \in \Rcn{n}$, by applying both the $\mathcal A$-contraction and $\mathcal D$-contraction to the lift $\hat{\bm i} \in \Ran{2n-1}$, we obtain the lift of ${\bf cont}({\bm i}) \in \Rcn{n-1}$.

	\begin{example}
		For $\bm i = (1,3,2,1,3,2,1,3,2) \in \Rcn{3}$, we have 
		\[
			{\bf cont}(\bm i) = (2,1,2,1).
		\]
Moreover, the lift $\hat{\bm i}$ of $\bm i$ is $\hat{\bm i} = (1,5,3,2,4,1,5,3,2,4,1,5,3,2,4) \in \Ran{5}$. 
By applying the $\mathcal A$-contraction and $\mathcal D$-contraction to $\hat{\bm i}$, we obtain $(2,1,3,2,1,3) \in \Ran{3}$, which is the lift of $(2,1,2,1) \in \Rcn{2}$. 
	\end{example}

	For a given reduced word ${\bm i} \in R(w_0^{(A_n)})$, 
	every rigorous path for a ($\mcal{D}$- or $\mcal{A}$-) contraction of ${\bm i}$ uniquely defines a rigorous path for ${\bm i}$ in a natural way (see \cite[Corollary~5.2]{CKLP21_GC}). 
	Similarly, for ${\bm i} \in \RCn$, one can immediately check that every symplectic rigorous path in $G^{\symp}({\bf cont}(\bm i))$ uniquely defines a symplectic rigorous path 
	in $G^{\symp}(\bm i)$, that is, we have a canonical injection:
\[
\iota_{\bm i} \colon \Paths^{\symp}({\bf cont}({\bm i})) \hookrightarrow \Paths^{\symp}({\bm i}).
\]
	
	\begin{definition}
		We say a symplectic rigorous path $P$ in $G^{\symp}(\bm {i})$ is {\em new} if it  is not contained in the image of $\iota_{\bm i}$, that is, it
does not come from a symplectic rigorous path in $G^{\symp}({\bf cont}({\bm i}))$.
	\end{definition} 
	
From now on, we label each node of a symplectic wiring diagram by $t_{i,j} \coloneqq \ell_i \cap \ell_j$ for simplicity. 
	There is a canonical way of constructing a new symplectic rigorous path in $G^{\symp}(\bm {i})$ having a peak $t_{1,i}$ for $i=2,\dots,n,\overline{n}, \ldots, \overline{1}$; or $t_{j,\overline{1}}$ for $j=1,\dots,n,\overline{n}, \ldots, \overline{2}$.
Before providing further explanations, we notice that the wire $\ell_1$ or $\ell_{\bar{1}}$ divides the diagram $G^{\symp}({\bm i})$ into two regions: the \emph{upper} region and the \emph{lower} region. 
For $j =1, \overline{1}$, we say that a path $P$ is \emph{below} $\ell_j$ if $P$ travels only the lower region given by $\ell_j$. 
Similarly, a node $t$ is said to be \emph{below} (resp., \emph{above}) $\ell_j$ if $t$ lies in the lower region (resp., the upper region) given by $\ell_j$. 
We also use similar notations and terminologies for rigorous paths in $G(\hat{\bm i})$. 
	Recall from \cite[Propositions 5.7, 5.8]{CKLP21_GC} that for each $2 \leq i \leq 2n$, there is a unique new rigorous path $\hat{P}_{1,i}$ in $G(\hat{\bm i})$, called a {\em canonical rigorous path}, such that 
	\begin{itemize}
		\item it has a unique peak $t_{1,i}$, that is, $\Lambda(\hat{P}_{1,i}) = \{t_{1,i}\}$;
		\item it is below $\ell_1$;
		\item its wire-expression is of the form: 
			\[
				\ell_{s_1} \rightarrow \cdots \rightarrow \ell_{s_q} \rightarrow \ell_1 \rightarrow \ell_{i} \rightarrow \ell_{r_1} \rightarrow \cdots \rightarrow \ell_{r_p},
			\]
			where the sequences $r_1, \dots, r_p$ and $s_1, \dots, s_q$ are decreasing;
		\item there is no wire $\ell_y$ with $y > s_1$ such that $(\ell_{s_1} \rightarrow \cdots \rightarrow \ell_{s_q} \rightarrow \ell_1)$ crosses $\ell_y$ before $t_{1,i}$.
	\end{itemize}
	Similarly, there is a unique new rigorous path $\hat{P}_{j,2n}$ in $G(\hat{\bm i})$ for $1 \leq j \leq 2n-1$, called a {\em canonical rigorous path}, such that 
	\begin{itemize}
		\item it has a unique peak $t_{j,2n}$, that is, $\Lambda(\hat{P}_{j,2n}) = \{t_{j,2n}\}$;
		\item it is below $\ell_{2n}$;
		\item its wire-expression is of the form: 
			\[
				\ell_{r_p} \rightarrow \cdots \rightarrow \ell_{r_1} \rightarrow \ell_j \rightarrow \ell_{2n} \rightarrow \ell_{s_q} \rightarrow \cdots \rightarrow \ell_{s_1},
			\]
			where the sequences $r_1, \dots, r_p$ and $s_1, \dots, s_q$ are increasing;
		\item there is no wire $\ell_y$ with $y > r_p$ such that $(\ell_{r_p} \rightarrow \cdots \rightarrow \ell_{r_1} \rightarrow \ell_j)$ crosses $\ell_y$ before $t_{j,2n}$.
	\end{itemize}
Let $P_{1,i}$ and $P_{j,2n}$ be symplectic rigorous paths in $G^{\symp}({\bm i})$ whose lifts are $\hat{P}_{1,i}$ and $\hat{P}_{j,2n}$, respectively. 
Then it is easy to see that $P_{1,i}^\vee = P_{2n+1-i,2n}$ for all $2 \leq i \leq 2n$. 
In the case of type $C_n$, we define a {\em canonical} symplectic rigorous path as follows.
	
\begin{definition}
Let $P \in \Paths^{\symp}({\bm i})$ be a symplectic rigorous path whose lift $\hat{P}$ is a canonical rigorous path in $G(\hat{\bm i},k)$.
Then we define the corresponding canonical symplectic rigorous path in the following way. 
\begin{itemize}
\item[Step 1.] If $k \leq n$, then set $\tilde{P} = P$. Otherwise, set $\tilde{P} = P^{\vee}$. 
\item[Step 2.] The extension $\tilde{P}_{\rm ex}$ of $\tilde{P}$ given in Definition \ref{d:extension_k_equal_n} or Proposition \ref{proposition_extension_exists} is called the {\em canonical symplectic rigorous path} corresponding to $P$.
\end{itemize}
A symplectic rigorous path is called a {\em canonical symplectic rigorous path} if it is obtained in this way. 
\end{definition}
		
It is straightforward by definition and Proposition \ref{prop_non_redundancy_k_not_n_maximal} that a canonical symplectic rigorous path $P$ gives a non-redundant inequality $\F_P^{(C)}({\bm a}) \geq 0$ in the expression \eqref{eq:descriptions_of_string_cones_in_type_C_n}. 
In addition, every canonical symplectic rigorous path in $\Paths^{\symp}(\bm i,n)$ is symmetric. 
For canonical rigorous path $\hat{P}_{1,i}$ and $\hat{P}_{j,2n}$, let $Q_{1,i}$ and $Q_{j,2n}$ denote the corresponding canonical symplectic rigorous paths, respectively. 
If $i \geq n+1$, then we also write $P_{1, \overline{2n+1-i}} \coloneqq P_{1, i}$ and $Q_{1, \overline{2n+1-i}} \coloneqq Q_{1, i}$. 
We use similar notations for $P_{j,2n}$ and $Q_{j,2n}$. 
Since $P_{1,i}^\vee = P_{2n+1-i,2n}$, we have $Q_{1,i} = Q_{2n+1-i,2n}$.
Hence it suffices to consider $\{Q_{1,2n}, Q_{2,2n}, \ldots, Q_{2n-1,2n}\}$ or $\{Q_{1,2}, Q_{1,3}, \ldots, Q_{1,n}, Q_{1,2n}, Q_{2,2n}, \ldots, Q_{n,2n}\}$.

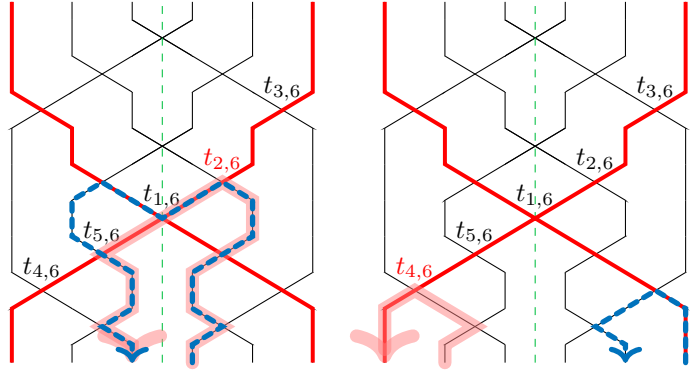
\begin{figure}
\begin{tikzpicture}[baseline =-0.5ex, scale = 0.8, yscale=0.6]
\def\nwires{3} 
\def\length{9} 

\foreach \i\y in {1/3, 2/2, 3/1, 4/3, 5/2, 6/3, 7/2, 8/1, 9/2}{
	\pgfmathsetmacro\ycoordinate{\length-\i+1}
	\pgfmathsetmacro\yplusone{\y+1}
	\foreach \x in {1,...,\nwires}{
		\ifthenelse{\x = \y}
			{\coordinate (\i LSW) at (-\nwires+\x-0.5, \ycoordinate);
			\coordinate (\i LSE) at (-\nwires+\x +0.5, \ycoordinate);
			\coordinate (\i LNW) at (-\nwires+\x-0.5, \ycoordinate+1);
			\coordinate (\i LNE) at (-\nwires+\x +0.5,\ycoordinate+1);

			\coordinate (\i RSW) at (\nwires - \x -0.5, \ycoordinate );
			\coordinate (\i RSE) at (\nwires - \x + 0.5, \ycoordinate );
			\coordinate (\i RNW) at (\nwires - \x - 0.5, \ycoordinate +1);
			\coordinate (\i RNE) at (\nwires - \x + 0.5, \ycoordinate +1);

			\coordinate (\i LC) at (-\nwires + \x, \ycoordinate+ 0.5);
			\coordinate (\i RC) at (\nwires - \x, \ycoordinate+0.5);
		
			\draw (-\nwires+\x-0.5,\ycoordinate) -- (\i LSW)-- (\i LNE) -- (-\nwires+\x+0.5, \ycoordinate+1);
			\draw (-\nwires+\x+0.5,\ycoordinate) --(\i LSE) --(\i LNW) --(-\nwires+\x-0.5, \ycoordinate+1);

			\draw (\nwires - \x -0.5, \ycoordinate) -- (\i RSW)-- (\i RNE) -- (\nwires - \x + 0.5, \ycoordinate + 1);
			\draw (\nwires - \x +0.5, \ycoordinate) -- (\i RSE)-- (\i RNW) -- (\nwires - \x -0.5, \ycoordinate + 1);
			}{}

		\ifthenelse{\x=\y \OR \x = \yplusone}{}{
			\draw
			(-\nwires+\x-0.5,\ycoordinate)--(-\nwires+\x-0.5,\ycoordinate+1);
			\draw
			(\nwires-\x+0.5,\ycoordinate)--(\nwires-\x+0.5,\ycoordinate+1);
		}
}


}

\foreach \x in{1,...,\nwires}{
	\draw (\x-0.5,\length+1)--(\x-0.5,\length+1.5);
	\draw (\x-0.5,1)--(\x-0.5,0.5);

	\draw (-\x+0.5,\length+1)--(-\x+0.5,\length+1.5);
	\draw (-\x+0.5,1)--(-\x+0.5,0.5);
}

\draw[teal!50!green, dashed] (0,0.5)--(0,\length+1.5);

\draw[color=red, line width = 1.5] 
(-2.5,0.5)--(8LSW)--(5RNE)--(3RSW)--(3RNE)--(2.5,10.5)
(2.5,0.5)--(8RSE)--(5LNW)--(3LSE)--(3LNW)--(-2.5,10.5);

\node[above] at (3RC) {$t_{3,6}$};
\node[above, red] at (5RC) {$t_{2,6}$};
\node[above] at (6RC) {$t_{1,6}$};
\node[above] at (7LC) {$t_{5,6}$};
\node[above] at (8LC) {$t_{4,6}$};

\draw[color=red!50, line cap=round, line width=5, opacity=0.5, ->]
(0.5,0.5)--(9RSW)--(9RC)--(9RNW)--(7RSW)--(7RNE)--(5RSE)--(5RC)--(7LC)--(7LSE)--(9LNE)--(9LC)--(9LSE)--(-0.5,0.5);

\draw[color=cyan!50!blue, line cap = round, line width = 2, dashed, ->] 
(0.5,0.5)--(9RSW)--(9RC)--(9RNW)--(7RSW)--(7RNE)--(5RSE)--(5RC)--(6RC)--(5LC)--(5LSW)--(7LNW)--(7LSE)--(9LNE)--(9LC)--(9LSE)--(-0.5,0.5);
\end{tikzpicture} \quad
\begin{tikzpicture}[baseline =-0.5ex, scale = 0.8, yscale=0.6]
\def\nwires{3} 
\def\length{9} 

\foreach \i\y in {1/3, 2/2, 3/1, 4/3, 5/2, 6/3, 7/2, 8/1, 9/2}{
	\pgfmathsetmacro\ycoordinate{\length-\i+1}
	\pgfmathsetmacro\yplusone{\y+1}
	\foreach \x in {1,...,\nwires}{
		\ifthenelse{\x = \y}
			{\coordinate (\i LSW) at (-\nwires+\x-0.5, \ycoordinate);
			\coordinate (\i LSE) at (-\nwires+\x +0.5, \ycoordinate);
			\coordinate (\i LNW) at (-\nwires+\x-0.5, \ycoordinate+1);
			\coordinate (\i LNE) at (-\nwires+\x +0.5,\ycoordinate+1);

			\coordinate (\i RSW) at (\nwires - \x -0.5, \ycoordinate );
			\coordinate (\i RSE) at (\nwires - \x + 0.5, \ycoordinate );
			\coordinate (\i RNW) at (\nwires - \x - 0.5, \ycoordinate +1);
			\coordinate (\i RNE) at (\nwires - \x + 0.5, \ycoordinate +1);

			\coordinate (\i LC) at (-\nwires + \x, \ycoordinate+ 0.5);
			\coordinate (\i RC) at (\nwires - \x, \ycoordinate+0.5);
		
			\draw (-\nwires+\x-0.5,\ycoordinate) -- (\i LSW)-- (\i LNE) -- (-\nwires+\x+0.5, \ycoordinate+1);
			\draw (-\nwires+\x+0.5,\ycoordinate) --(\i LSE) --(\i LNW) --(-\nwires+\x-0.5, \ycoordinate+1);

			\draw (\nwires - \x -0.5, \ycoordinate) -- (\i RSW)-- (\i RNE) -- (\nwires - \x + 0.5, \ycoordinate + 1);
			\draw (\nwires - \x +0.5, \ycoordinate) -- (\i RSE)-- (\i RNW) -- (\nwires - \x -0.5, \ycoordinate + 1);
			}{}

		\ifthenelse{\x=\y \OR \x = \yplusone}{}{
			\draw
			(-\nwires+\x-0.5,\ycoordinate)--(-\nwires+\x-0.5,\ycoordinate+1);
			\draw
			(\nwires-\x+0.5,\ycoordinate)--(\nwires-\x+0.5,\ycoordinate+1);
		}
}


}

\foreach \x in{1,...,\nwires}{
	\draw (\x-0.5,\length+1)--(\x-0.5,\length+1.5);
	\draw (\x-0.5,1)--(\x-0.5,0.5);

	\draw (-\x+0.5,\length+1)--(-\x+0.5,\length+1.5);
	\draw (-\x+0.5,1)--(-\x+0.5,0.5);
}

\draw[teal!50!green, dashed] (0,0.5)--(0,\length+1.5);

\draw[color=red, line width = 1.5] 
(-2.5,0.5)--(8LSW)--(5RNE)--(3RSW)--(3RNE)--(2.5,10.5)
(2.5,0.5)--(8RSE)--(5LNW)--(3LSE)--(3LNW)--(-2.5,10.5);

\node[above] at (3RC) {$t_{3,6}$};
\node[above] at (5RC) {$t_{2,6}$};
\node[above] at (6RC) {$t_{1,6}$};
\node[above] at (7LC) {$t_{5,6}$};
\node[above, red] at (8LC) {$t_{4,6}$};

\draw[color=red!50, line cap=round, line width=5, opacity=0.5, ->]
(-1.5,0.5)--(9LSW)--(9LC)--(8LC)--(8LSW)--(-2.5,0.5);

\draw[color=cyan!50!blue, line cap = round, line width = 2, dashed, ->] 
(2.5,0.5)--(8RSE)--(8RC)--(9RC)--(9RSE)--(1.5,0.5);
\end{tikzpicture} 
\caption{Canonical rigorous paths (red highlighted) and the corresponding canonical symplectic rigorous paths (blue dashed) for nodes $t_{2,6}$ and $t_{4,6}$, where $\bm i = (3,2,1,3,2,3,2,1,2)$.}\label{figure_canonical_symplectic_GP_paths}
\end{figure}

\begin{example}
Let $\bm i = (3,2,1,3,2,3,2,1,2) \in \Rcn{3}$. 
For each node $t_{j, 2n}$ on $G(\hat{\bm i})$, the canonical rigorous path $\hat{P}_{j, 2n}$ and the corresponding canonical symplectic rigorous path $Q_{j, 2n}$ are given as follows.
\begin{center}
\begin{tabular}{c|ll}
\toprule 
node & canonical rigorous path & canonical symplectic rigorous path \\
\midrule 
$t_{3,6}$ & $\ell_3 \to \ell_6 \to \ell_5 \to \ell_4$ 
& $\ell_3 \to \ell_{\bar{1}} \to \ell_{1} \to \ell_{\bar{3}}$ \\
$t_{2,6}$ & $\ell_3 \to \ell_2 \to \ell_6 \to \ell_5 \to \ell_4$ & $\ell_3 \to \ell_2 \to \ell_{\bar{1}} \to \ell_1 \to \ell_{\bar{2}} \to \ell_{\bar{3}}$ \\
$t_{1,6}$ & $\ell_3 \to \ell_2 \to \ell_1 \to \ell_6 \to \ell_5 \to \ell_4$ &  $\ell_3 \to \ell_2 \to \ell_1 \to \ell_{\bar{1}} \to \ell_{\bar{2}} \to \ell_{\bar{3}}$ \\
$t_{5,6}$ & $\ell_5 \to \ell_6$ & $\ell_{1} \to \ell_2$ \\
$t_{4,6}$ & $\ell_{5} \to \ell_{4} \to \ell_6$ & $\ell_1 \to \ell_3 \to \ell_2$ \\
\bottomrule
\end{tabular}
\end{center}
We depict paths for nodes $t_{2,6}$ and $t_{4,6}$ in Figure~\ref{figure_canonical_symplectic_GP_paths}. 
\end{example}

	\begin{proposition}\label{proposition_number_of_canonical_paths}
		For ${\bm i} \in \RCn$ with $n \geq 3$, there are exactly $(2n-1)$ number of canonical symplectic rigorous paths in $G^{\symp}({\bm i})$ each of which is indexed by the $(2n-1)$ nodes $t_{1,\bar{1}}, \dots, t_{n,\bar{1}}, t_{\bar{n},\bar{1}}, \ldots, t_{\bar{2},\bar{1}}$ on $\ell_{\bar{1}}$.
	\end{proposition}
	
	\begin{proof}
It is enough to show that $Q_{i,\bar{1}} \neq Q_{j,\bar{1}}$ if $i \neq j$. 
If $Q_{i,\bar{1}} = Q_{j,\bar{1}}$, then we have
		\[
			 \chamber(P_{i,\bar{1}}) \cup \chamber(P_{i,\bar{1}}^\vee) = \chamber(Q_{i,\bar{1}}) \cup \chamber(Q_{i,\bar{1}}^\vee)
			 = \chamber(Q_{j,\bar{1}}) \cup \chamber(Q_{j,\bar{1}}^\vee) = \chamber(P_{j,\bar{1}}) \cup \chamber(P_{j,\bar{1}}^\vee),
		\]
which gives a contradiction since $\chamber(P_{i,\bar{1}}) \cup \chamber(P_{i,\bar{1}}^\vee)$ has two peaks $t_{i,\bar{1}}$ and $t_{1,\bar{i}}$ whereas the peaks of
		$\chamber(P_{j,\bar{1}}) \cup \chamber(P_{j,\bar{1}}^\vee)$ are $t_{j,\bar{1}}$ and $t_{1,\bar{j}}$. 
		This completes the proof. 
	\end{proof}

Let $\lVert {\bm i} \rVert$ denote the number of facets of the string cone $\mathcal{C}_{\bm i}^{(C_n)}$.
By Proposition~\ref{proposition_number_of_canonical_paths}, we obtain the following corollary.

\begin{corollary}\label{cor_ineq_conti}
For ${\bm i} \in \RCn$ with $n \geq 3$, it holds that
\begin{equation}\label{equation_increasing}
			\lVert {\bm i} \rVert \geq \lVert {\bf cont}(\bm i) \rVert + (2n-1). 
\end{equation}
 In particular, we get $\lVert {\bm i} \rVert \geq 1 + 3 + \dots + (2n-1) = n^2$.
\end{corollary}

	\begin{lemma}\label{lemma_new_canonical}
		Take ${\bm i} \in \RCn$ with $n \geq 3$. 
Suppose that there is no node below $\ell_1$ or $\ell_{\bar{1}}$ except for nodes on $\ell_1$ and $\ell_{\bar{1}}$. Then every new symplectic rigorous path defining a non-redundant inequality is canonical. 
	\end{lemma}
	
	\begin{proof}
By the assumption, the symplectic wiring diagram $G^{\symp}({\bm i})$ is given as in Figure~\ref{figure_no_node_under_l1}.
		
\begin{figure}[H]
\centering
\begin{tikzpicture}[scale = 0.8, yscale=0.6]
\def\nwires{3} 
\def\length{9} 

\foreach \i\y in {1/2,2/3,3/2,4/3,5/1,6/2,7/3,8/2,9/1}{
	\pgfmathsetmacro\ycoordinate{\length-\i+1}
	\pgfmathsetmacro\yplusone{\y+1}
	\foreach \x in {1,...,\nwires}{
		\ifthenelse{\x = \y}
			{\coordinate (\i LSW) at (-\nwires+\x-0.5, \ycoordinate);
			\coordinate (\i LSE) at (-\nwires+\x +0.5, \ycoordinate);
			\coordinate (\i LNW) at (-\nwires+\x-0.5, \ycoordinate+1);
			\coordinate (\i LNE) at (-\nwires+\x +0.5,\ycoordinate+1);

			\coordinate (\i RSW) at (\nwires - \x -0.5, \ycoordinate );
			\coordinate (\i RSE) at (\nwires - \x + 0.5, \ycoordinate );
			\coordinate (\i RNW) at (\nwires - \x - 0.5, \ycoordinate +1);
			\coordinate (\i RNE) at (\nwires - \x + 0.5, \ycoordinate +1);

			\coordinate (\i LC) at (-\nwires + \x, \ycoordinate+ 0.5);
			\coordinate (\i RC) at (\nwires - \x, \ycoordinate+0.5);
		
			\draw (-\nwires+\x-0.5,\ycoordinate) -- (\i LSW)-- (\i LNE) -- (-\nwires+\x+0.5, \ycoordinate+1);
			\draw (-\nwires+\x+0.5,\ycoordinate) --(\i LSE) --(\i LNW) --(-\nwires+\x-0.5, \ycoordinate+1);

			\draw (\nwires - \x -0.5, \ycoordinate) -- (\i RSW)-- (\i RNE) -- (\nwires - \x + 0.5, \ycoordinate + 1);
			\draw (\nwires - \x +0.5, \ycoordinate) -- (\i RSE)-- (\i RNW) -- (\nwires - \x -0.5, \ycoordinate + 1);
			}{}

		\ifthenelse{\x=\y \OR \x = \yplusone}{}{
			\draw
			(-\nwires+\x-0.5,\ycoordinate)--(-\nwires+\x-0.5,\ycoordinate+1);
			\draw
			(\nwires-\x+0.5,\ycoordinate)--(\nwires-\x+0.5,\ycoordinate+1);
		}
}


}

\foreach \x in{1,...,\nwires}{
	\draw (\x-0.5,\length+1)--(\x-0.5,\length+1.5);
	\draw (\x-0.5,1)--(\x-0.5,0.5);

	\draw (-\x+0.5,\length+1)--(-\x+0.5,\length+1.5);
	\draw (-\x+0.5,1)--(-\x+0.5,0.5);
}
\draw[color=red, line width = 1.5] 
(-2.5,0.5)--(9LSW)--(9LNE)--(8LSW)--(8LNE)--(7LSW)--(7LNE)--(6RSW)--(6RNE)--(5RSW)--(5RNE)--(2.5,10.5);
\draw[color=red, line width = 1.5] 
(2.5,0.5)--(9RSE)--(9RNW)--(8RSE)--(8RNW)--(7RSE)--(7RNW)--(6LSE)--(6LNW)--(5LSE)--(5LNW)--(-2.5,10.5);

\node[above] at (9RC) {$t_{1,2}$};
\node[above] at (8RC) {$t_{1,3}$};
\node[above] at (5RC) {$t_{2,\bar{1}}$};
\node[above] at (6RC) {$t_{3,\bar{1}}$};

\end{tikzpicture}
\caption{No node below $\ell_1$ or $\ell_{\bar{1}}$ except for nodes on $\ell_1$ and $\ell_{\bar{1}}$.}\label{figure_no_node_under_l1}
\end{figure}
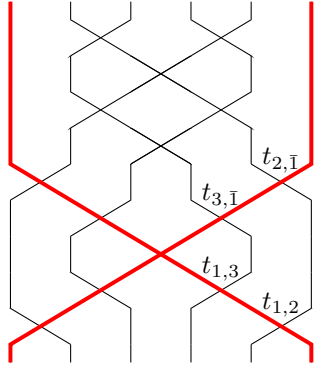

		Using Figure~\ref{figure_no_node_under_l1}, we describe all canonical symplectic rigorous paths explicitly as follows. 
		For $t_{1,i}$ with $2 \leq i \leq n$, the canonical symplectic rigorous path $Q_{1,i}$ having a peak $t_{1,i}$ is given by 
		\[
			\begin{array}{cccccl}
				\ell_{i-1} & \rightarrow & \ell_1 & \rightarrow & \ell_{i} & \quad \quad \text{if $i > 2$,} \\ 
				\ell_1 & \rightarrow & \ell_2 & & & \quad \quad \text{if $i = 2$.} \\ 
			\end{array}
		\] 
		Similarly, for $t_{j,\bar{1}}$ with $2 \leq j \leq n$, the canonical symplectic rigorous path $Q_{j,\bar{1}}$ having a peak $t_{j,\bar{1}}$ is given by 
		\[
			\begin{array}{cccccccl}
				\ell_j & \rightarrow & \ell_{\bar{1}} & \rightarrow & \ell_{j+1} & & & \quad \quad \text{if $2 \leq j < n$,} \\ 
				\ell_n & \rightarrow & \ell_{\bar{1}} & \rightarrow & \ell_{1} & \rightarrow & \ell_{\bar{n}} & \quad \quad \text{if $j = n$.} \\ 
			\end{array}
		\] 
		For $t_{1,\bar{1}}$, the canonical symplectic rigorous path $Q_{1,\bar{1}}$ having a peak $t_{1,\bar{1}}$ is given by 
		\[
				\ell_n \rightarrow \ell_1 \rightarrow \ell_{\bar{1}} \rightarrow \ell_{\bar{n}}.
		\]	
		It is straightforward that there is no other new symplectic rigorous path defining a non-redundant inequality (cf.\ Figure~\ref{figure_no_node_under_l1}). 
		We notice that the path $(\ell_{n}\to \ell_{\overline{1}}\to \ell_{\overline{n}})$ and its mirror $(\ell_{{n}}\to \ell_{{1}}\to \ell_{\overline{n}})$ are new symplectic rigorous paths but they define redundant inequalities by Proposition~\ref{prop_non_redundancy_n}.
		This completes the proof.
	\end{proof}

We denote by 
\[
			\begin{split}
				&{\bm i}_{C}^{(n)} \coloneqq (n,n-1,n,n-1,n-2,n-1,n,n-1,n-2,\dots,2,\dots,n\dots,2,1,2,\dots,n,\dots,2,1),\\
				&{\bm j}_{C}^{(n)} \coloneqq (n-1,n,n-1,n,n-2,n-1,n,n-1,n-2,\dots,2,\dots,n\dots,2,1,2,\dots,n,\dots,2,1).
			\end{split}
		\]

	\begin{corollary}\label{corollary_simplicial}
		If $\bm i$ is either ${\bm i}_{C}^{(n)}$ or ${\bm j}_{C}^{(n)}$, then the string cone $\mcal{C}_{\bm i}^{(C_n)}$ is simplicial, that is, $\lVert \bm i \rVert = n^2$.
	\end{corollary}
	
	\begin{proof}
		We proceed by induction on $n$. 
If $n = 2$, then we see by Examples~\ref{example_2121} and~\ref{example_1212} that 
\[
\lVert (2,1,2,1) \rVert = \lVert (1,2,1,2) \rVert = 4.
\]
For $n \geq 3$, each reduced expression ends with $(1,2,\dots,n,\dots,2,1)$, and this implies that there is no node below $\ell_1$ or $\ell_{\bar{1}}$ except for nodes on $\ell_1$ and $\ell_{\bar{1}}$.  
		By Lemma \ref{lemma_new_canonical}, every new symplectic rigorous path defining a non-redundant inequality is canonical. Hence there are exactly $(2n-1)$ new ones, and each other path can be obtained from a symplectic rigorous path for ${\bf cont}({\bm i})$, that is, we have 
\[
\lVert {\bm i} \rVert = \lVert {\bf cont}({\bm i}) \rVert + (2n-1).
\] 
This implies by the induction hypothesis that 
\[
\lVert {\bm i} \rVert = (n-1)^2 + (2n-1) = n^2,
\]
which coincides with the dimension of $\mcal{C}_{\bm i}^{(C_n)}$. 
Thus, the result follows. 
	\end{proof}

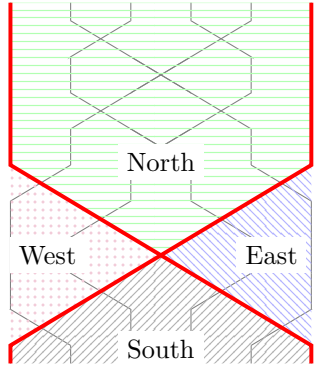
\begin{figure}[H]
\begin{tikzpicture}[scale = 0.8, yscale=0.6]
\def\nwires{3} 
\def\length{9} 

\foreach \i\y in {1/2,2/3,3/2,4/3,5/1,6/2,7/3,8/2,9/1}{
	\pgfmathsetmacro\ycoordinate{\length-\i+1}
	\pgfmathsetmacro\yplusone{\y+1}
	\foreach \x in {1,...,\nwires}{
		\ifthenelse{\x = \y}
			{\coordinate (\i LSW) at (-\nwires+\x-0.5, \ycoordinate);
			\coordinate (\i LSE) at (-\nwires+\x +0.5, \ycoordinate);
			\coordinate (\i LNW) at (-\nwires+\x-0.5, \ycoordinate+1);
			\coordinate (\i LNE) at (-\nwires+\x +0.5,\ycoordinate+1);

			\coordinate (\i RSW) at (\nwires - \x -0.5, \ycoordinate );
			\coordinate (\i RSE) at (\nwires - \x + 0.5, \ycoordinate );
			\coordinate (\i RNW) at (\nwires - \x - 0.5, \ycoordinate +1);
			\coordinate (\i RNE) at (\nwires - \x + 0.5, \ycoordinate +1);

			\coordinate (\i LC) at (-\nwires + \x, \ycoordinate+ 0.5);
			\coordinate (\i RC) at (\nwires - \x, \ycoordinate+0.5);
		
			\draw[gray] (-\nwires+\x-0.5,\ycoordinate) -- (\i LSW)-- (\i LNE) -- (-\nwires+\x+0.5, \ycoordinate+1);
			\draw[gray] (-\nwires+\x+0.5,\ycoordinate) --(\i LSE) --(\i LNW) --(-\nwires+\x-0.5, \ycoordinate+1);

			\draw[gray] (\nwires - \x -0.5, \ycoordinate) -- (\i RSW)-- (\i RNE) -- (\nwires - \x + 0.5, \ycoordinate + 1);
			\draw[gray] (\nwires - \x +0.5, \ycoordinate) -- (\i RSE)-- (\i RNW) -- (\nwires - \x -0.5, \ycoordinate + 1);
			}{}

		\ifthenelse{\x=\y \OR \x = \yplusone}{}{
			\draw[gray]
			(-\nwires+\x-0.5,\ycoordinate)--(-\nwires+\x-0.5,\ycoordinate+1);
			\draw[gray]
			(\nwires-\x+0.5,\ycoordinate)--(\nwires-\x+0.5,\ycoordinate+1);
		}
}


}

\foreach \x in{1,...,\nwires}{
	\draw[gray] (\x-0.5,\length+1)--(\x-0.5,\length+1.5);
	\draw[gray] (\x-0.5,1)--(\x-0.5,0.5);

	\draw[gray] (-\x+0.5,\length+1)--(-\x+0.5,\length+1.5);
	\draw[gray] (-\x+0.5,1)--(-\x+0.5,0.5);
}

\filldraw[draw=none,pattern=north east lines,pattern color=black!30]
(-2.5,0.5)--(2.5,0.5)--(9RSE)--(9RC)--(7RC)--(9LC)--(9LSW)--cycle; 
\node[fill=white, below=of 7RC] {South};

\filldraw[draw=none, pattern=north west lines, pattern color = blue!30]
(9RSE)--(7RC)--(5RNE)--cycle;
\node[fill=white, right=of 7RC] {East};

\filldraw[draw=none, pattern= dots, pattern color = purple!20]
(9LSW)--(7RC)--(5LNW)--cycle;
\node[fill=white, left=of 7RC] {West};

\filldraw[draw=none, pattern=horizontal lines, pattern color = green!30] 
(5LNW)--(7RC)--(5RNE)--(2.5, 10.5)--(-2.5, 10.5)--cycle;
\node[fill=white, above=of 7RC] {North};

\draw[color=red, line width = 1.5] 
(-2.5,0.5)--(9LSW)--(9LNE)--(8LSW)--(8LNE)--(7LSW)--(7LNE)--(6RSW)--(6RNE)--(5RSW)--(5RNE)--(2.5,10.5);
\draw[color=red, line width = 1.5] 
(2.5,0.5)--(9RSE)--(9RNW)--(8RSE)--(8RNW)--(7RSE)--(7RNW)--(6LSE)--(6LNW)--(5LSE)--(5LNW)--(-2.5,10.5);


\end{tikzpicture}
\caption{The north, east, west, and south sectors.}\label{fig_NEWS_sectors}
\end{figure}

We notice that the union of wires $\ell_1$ and $\ell_{\bar{1}}$ divides the diagram $G^{\symp}(\bm i)$ into four regions: the \emph{north sector}, the \emph{east sector}, the \emph{west sector}, and the \emph{south sector}.
See Figure~\ref{fig_NEWS_sectors}. 
	\begin{lemma}\label{lemma_wall}
Take ${\bm i} \in \RCn$ with $n \geq 3$. 
If there is a node on the wall below $\ell_1$ \textup{(}as well as $\ell_{\bar{1}}$\textup{)} on $G^{\symp}(\bm i)$ except for $t_{1, \bar{1}}$, then there is a non-canonical new symplectic rigorous path $P$ such that the inequality $\F_P^{(C)}({\bm a}) \geq 0$ is non-redundant in the expression \eqref{eq:descriptions_of_string_cones_in_type_C_n}.
	\end{lemma}
	
	\begin{proof}
We divide the proof into three cases. 

		\noindent
		{\bf Case 1:} there exists $1 \leq k \leq n-1$ such that $t_{k,\overline{k}}$ is above $\ell_1$ and such that $t_{k+1,\overline{k+1}}$ is below $\ell_1$ (see Figure~\ref{fig_forbidden_to_nonforbidden}). 
In this case, we further take cases: 

		\noindent
		{\bf Case 1-a:} $\ell_k$ and $\ell_{\overline{k+1}}$ do not meet in the east sector (see Figure~\ref{fig_forbidden_to_nonforbidden}). 
Then we consider the following two new symplectic rigorous paths: 
		\[
			\breve{P}_{k,\bar{1}} \colonequals (\ell_k \rightarrow \ell_{\bar{1}} \rightarrow \ell_{\overline{k+1}} \rightarrow \ell_{1} \rightarrow \ell_{k+1}) \quad \text{and} \quad 
			\widetilde{P}_{k,\bar{1}} \colonequals (\ell_k \rightarrow \ell_{\bar{1}} \rightarrow \ell_{\overline{k+1}} \rightarrow \ell_{k+1})
		\]
		which are the red highlighted path and the blue dashed path in Figure~\ref{fig_forbidden_to_nonforbidden}, respectively. 
If $\breve{P}_{k,\bar{1}}$ or $\widetilde{P}_{k,\bar{1}}$ contains a forbidden pattern, we modify it into a non-forbidden path as in Figure~\ref{fig_forbidden_to_nonforbidden}, where $\breve{P}_{k,\bar{1}}$ is thought of as a modified path obtained from 
\[(\ell_k \rightarrow \ell_{\bar{1}} \rightarrow \ell_{1} \rightarrow \ell_{k+1}).\]
It is easy to see that such modification is always possible. 
Since $\widetilde{P}_{k,\bar{1}}$ does not cross the wall, it is obvious that $\widetilde{P}_{k,\bar{1}}$ satisfies the assumption of Proposition \ref{prop_non_redundancy_k_not_n_maximal}, and hence that $\widetilde{P}_{k,\bar{1}} = Q_{k,\bar{1}}$.
Since $\chamber(\widetilde{P}_{k,\bar{1}}) \subsetneq \chamber(\breve{P}_{k,\bar{1}}) \subseteq \chamber((\breve{P}_{k,\bar{1}})_{\rm ex})$, we have $(\breve{P}_{k,\bar{1}})_{\rm ex} \neq Q_{k,\bar{1}}$. 
Since $\ell_1, \ldots, \ell_k$ are oriented upward and $\ell_{k+1}, \ldots, \ell_n, \ell_{\bar{n}}, \ldots \ell_{\bar{1}}$ are oriented downward, it is easy to see that $(\breve{P}_{k,\bar{1}})_{\rm ex}$ has a peak $t_{k,\bar{1}}$, which implies that $(\breve{P}_{k,\bar{1}})_{\rm ex}$ is a new non-canonical symplectic rigorous path.\vs{0.1cm}

\begin{figure}[H]
\begin{tikzpicture}[scale = 0.55, yscale=0.6]
\def\nwires{5} 
\def\length{15} 

\foreach \i\y in {1/1, 2/2, 3/3, 4/4, 5/5, 6/3,
7/4, 8/2, 9/3, 10/2, 11/4, 12/1, 13/5, 14/2, 15/4}{
	\pgfmathsetmacro\ycoordinate{\length-\i+1}
	\pgfmathsetmacro\yplusone{\y+1}
	\foreach \x in {1,...,\nwires}{
		\ifthenelse{\x = \y}
			{\coordinate (\i LSW) at (-\nwires+\x-0.5, \ycoordinate);
			\coordinate (\i LSE) at (-\nwires+\x +0.5, \ycoordinate);
			\coordinate (\i LNW) at (-\nwires+\x-0.5, \ycoordinate+1);
			\coordinate (\i LNE) at (-\nwires+\x +0.5,\ycoordinate+1);

			\coordinate (\i RSW) at (\nwires - \x -0.5, \ycoordinate );
			\coordinate (\i RSE) at (\nwires - \x + 0.5, \ycoordinate );
			\coordinate (\i RNW) at (\nwires - \x - 0.5, \ycoordinate +1);
			\coordinate (\i RNE) at (\nwires - \x + 0.5, \ycoordinate +1);

			\coordinate (\i LC) at (-\nwires + \x, \ycoordinate+ 0.5);
			\coordinate (\i RC) at (\nwires - \x, \ycoordinate+0.5);
		
			}{}

}


}

\node[circle, draw, red, pin={[pin distance = 1cm]0:forbidden}] at (2RC) {};

\draw[teal!50!green, dashed] (0,0.5)--(0,\length+1.5);

\draw[color=red, line width=1.5]
(-4.5,0.5)--(12LSW)--(12LNE)--(10LSW)--(10LNE)--(9LSW)--(9LNE)
--(7LSW)--(7LNE)--(5LSW)--(5LNE)--(4RSW)--(4RNE)--(3RSW)--(3RNE)
--(2RSW)--(2RNE)--(1RSW)--(1RNE)--(4.5,16.5)
(4.5,0.5)--(12RSE)--(12RNW)--(10RSE)--(10RNW)--(9RSE)--(9RNW)
--(7RSE)--(7RNW)--(5RSE)--(5RNW)--(4LSE)--(4LNW)--(3LSE)--(3LNW)
--(2LSE)--(2LNW)--(1LSE)--(1LNW)--(-4.5,16.5);

\draw 
(2.5,0.5)--(14RSW)--(14RNE)--(12RSW)--(12RNE)--(1RSE)--(1RNW)--(3.5,16.5)
(1.5,0.5)--(15RSE)--(15RNW)--(13RSE)--(13RNW)--(11LSE)--(11LNW)--
(9LSE)--(9LNW)--(8LSE)--(8LNW)--(2LSW)--(2LNE)--(-2.5,16.5);

\draw[color=teal]
(15LNE)--(13LSW)--(13LNE)--(11RSW)--(11RNE)--(9RSW)--(9RNE)
--(8RSW)--(8RNE)--(2RSE)--(2RNW)--(2.5,16);

\draw [color=red!50, line cap=round, line width=5, opacity=0.5, ->] 
	(2.5,0)--(14RSW)--(14RNE)--(12RSW)--(12RNE)--(1RSE)--(1RC)--(2RC)--(2RSE)--(8RNE)--(9RC)--(9RNW)--(7RSE)--(7RNW)--(5RSE)--(2LC)--(2LSW)--(8LNW)--(9LSE)--(11LNW)--(11LSE)--(13LNW)--(13LSE)--(15RNW)--(15RSE)--(1.5,0);

\draw[color=cyan!50!blue, line cap = round, line width = 2, dashed, ->] 
	(2.5,0)--(14RSW)--(14RNE)--(12RSW)--(12RNE)--(1RSE)--(1RC)--(2RC)--(2RSE)--(8RNE)--(9RSW)--(11RNE)--(11RSW)--(13LNE)--(13LC)--(13RSE)--(15RNW)--(15RSE)--(1.5,0);

\node[below] at (2.5,-0.5) {$L_k$};
\node[below] at (1.5,-0.5) {$L_{k+1}$};
\end{tikzpicture}
\caption{$t_{k+1,\overline{k+1}}$ on the wall below $\ell_1$ and $\ell_{\bar{1}}$; $\ell_k$ and $\ell_{\overline{k+1}}$ do not meet on the east sector.}\label{fig_forbidden_to_nonforbidden}
\end{figure}

		\noindent
		{\bf Case 1-b:} $\ell_k$ and $\ell_{\overline{k+1}}$ meet in the east sector (see Figure~\ref{fig_node_on_the_wall}). 
Then we consider the following two symplectic rigorous paths: 
		\[
			\breve{P}_{k,\bar{1}} \colonequals (\ell_k \rightarrow \ell_{\bar{1}} \rightarrow \ell_1 \rightarrow \ell_{\bar{k}} \rightarrow \ell_{k+1}) \quad \text{and} \quad 
			\widetilde{P}_{k,\bar{1}} \colonequals (\ell_k \rightarrow \ell_{\bar{1}} \rightarrow \ell_1 \rightarrow \ell_{k+1})
		\]
		which are the red highlighted path and the blue dashed path in Figure~\ref{fig_node_on_the_wall}, respectively. 
If $\breve{P}_{k,\bar{1}}$ or $\widetilde{P}_{k,\bar{1}}$ contains a forbidden pattern, we modify it into a non-forbidden path as in Case 1-a. 
Since $\ell_1, \ldots, \ell_k$ are oriented upward and $\ell_{k+1}, \ldots, \ell_n, \ell_{\bar{n}}, \ldots \ell_{\bar{1}}$ are oriented downward, it is easy to see that $\breve{P}_{k,\bar{1}}$ and $\widetilde{P}_{k,\bar{1}}$ satisfy the assumption of Proposition \ref{prop_non_redundancy_k_not_n_maximal}. 
In particular, these paths give non-redundant inequalities. 
Since 
\[
{P}_{k,\bar{1}}  \colonequals (\ell_k \rightarrow \ell_{\bar{1}} \rightarrow \ell_{k+1}),
\] 
it follows that
\[
\chamber({P}_{k,\bar{1}}) \subseteq \chamber(\breve{P}_{k,\bar{1}}) \subseteq  \chamber(\breve{P}_{k,\bar{1}}) \cup \chamber(\breve{P}_{k,\bar{1}}^\vee) = \chamber(P_{k,\bar{1}}) \cup \chamber(P_{k,\bar{1}}^\vee). 
\]
In particular, we have $\breve{P}_{k,\bar{1}} = Q_{k,\bar{1}}$ is a canonical symplectic rigorous path, which implies that $\widetilde{P}_{k,\bar{1}}$ is a non-canonical new symplectic rigorous path. \vs{0.1cm}

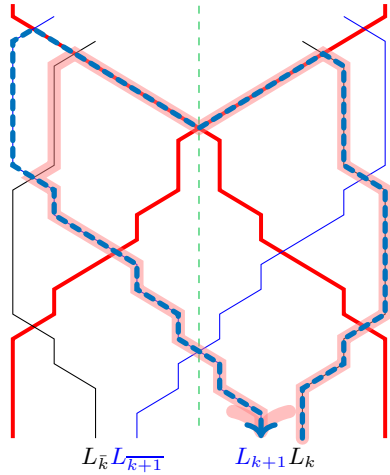
\begin{figure}[H]
\begin{tikzpicture}[scale = 0.55, yscale=0.6]
\def\nwires{5} 
\def\length{16} 

\foreach \i\y in {1/1, 2/2, 3/3, 4/4, 5/5, 6/3,
7/1, 
8/4, 9/2, 10/3, 11/2, 12/4, 13/1, 14/5, 15/2, 16/4}{
	\pgfmathsetmacro\ycoordinate{\length-\i+1}
	\pgfmathsetmacro\yplusone{\y+1}
	\foreach \x in {1,...,\nwires}{
		\ifthenelse{\x = \y}
			{\coordinate (\i LSW) at (-\nwires+\x-0.5, \ycoordinate);
			\coordinate (\i LSE) at (-\nwires+\x +0.5, \ycoordinate);
			\coordinate (\i LNW) at (-\nwires+\x-0.5, \ycoordinate+1);
			\coordinate (\i LNE) at (-\nwires+\x +0.5,\ycoordinate+1);

			\coordinate (\i RSW) at (\nwires - \x -0.5, \ycoordinate );
			\coordinate (\i RSE) at (\nwires - \x + 0.5, \ycoordinate );
			\coordinate (\i RNW) at (\nwires - \x - 0.5, \ycoordinate +1);
			\coordinate (\i RNE) at (\nwires - \x + 0.5, \ycoordinate +1);

			\coordinate (\i LC) at (-\nwires + \x, \ycoordinate+ 0.5);
			\coordinate (\i RC) at (\nwires - \x, \ycoordinate+0.5);
		
			}{}

}

}

\draw[teal!50!green, dashed] (0,0.5)--(0,\length+1.5);
\draw[color=red, line width=1.5]
(-4.5,0)--(13LSW)--(13LNE)--(11LSW)--(11LNE)--(10LNE)--(8LSW)--(8LNE)--(5LSW)--(5LNE)--(1RNE)--(4.5, 17.5)
(4.5,0)--(13RSE)--(13RNW)--(11RSE)--(11RNW)--(10RNW)--(8RSE)--(8RNW)--(5RSE)--(5RNW)--(1LNW)--(-4.5,17.5);

\draw 
(2.5,0) node[below] {\small $L_k$}--(15RSW)--(15RNE)--(13RSW)--(13RNE)--(7RSE)--(7RNW)--(2RSE)--(2RNW)
(-2.5,0) node[below] {\small $L_{\bar{k}}$}--(15LSE)--(15LNW)--(13LSE)--(13LNW)--(7LSW)--(7LNE)--(2LSW)--(2LNE);

\draw[blue]
(1.5,0) node[below] {\small $L_{k+1}$}--(16RSE)--(16RNW)--(14RSE)--(14RNW)--(12LSE)--(12LNW)--(10LSE)--(9LNW)--(7LSE)--(7LNW)--(1LSW)--(1LNE)
(-1.5,0) node[below] {\small $L_{\overline{k+1}}$}--(16LSW)--(16LNE)--(14LSW)--(14LNE)--(12RSW)--(12RNE)--(10RSW)--(9RNE)--(7RSW)--(7RNE)--(1RSE)--(1RNW);

\draw [color=red!50, line cap=round, line width=5, opacity=0.5, ->] 
(2.5,0)--(15RSW)--(15RNE)--(13RSW)--(13RNE)--(7RSE)--(7RNW)--(2RSE)--(2RC)--(5RC)--(2LC)--(2LSW)--(7LNE)--(7LC)--(7LSE)--(9LNW)--(10LSE)--(12LNW)--(12LSE)--(14LNW)--(14LSE)--(16RNW)--(16RSE)--(1.5,0);

\draw[color=cyan!50!blue, line cap = round, line width = 2, dashed, ->] 
(2.5,0)--(15RSW)--(15RNE)--(13RSW)--(13RNE)--(7RSE)--(7RNW)--(2RSE)--(2RC)--(5RC)--(1LC)--(1LSW)--(7LNW)--(7LSE)--(9LNW)--(10LSE)--(12LNW)--(12LSE)--(14LNW)--(14LSE)--(16RNW)--(16RSE)--(1.5,0);

\end{tikzpicture}
\caption{$t_{k+1,\overline{k+1}}$ on the wall below $\ell_1$ and $\ell_{\bar{1}}$; $\ell_k$ and $\ell_{\overline{k+1}}$ meet on the east sector.}\label{fig_node_on_the_wall}
\end{figure}

		\noindent
		{\bf Case 2:} $t_{2,\overline{2}}$ is below $\ell_1$ and $t_{n,\overline{n}}$ is above $\ell_1$ (see Figure~\ref{fig_case_2_paths_modification}). 
In this case, there are two new symplectic rigorous paths:
\[
\breve{P}_{1,\bar{1}} \colonequals (\ell_n \to \ell_1 \to \ell_{\bar{1}} \to \ell_{\bar{n}}) \quad \text{ and } \quad 
\widetilde{P}_{1,\bar{1}} \colonequals (\ell_1 \to \ell_{\bar{1}} \to \ell_{\bar{n}} \to \ell_2)
\]
which are the red highlighted path and the blue dashed path in Figure~\ref{fig_case_2_paths_modification}, respectively.
If $\breve{P}_{1,\bar{1}}$ or $\widetilde{P}_{1,\bar{1}}$ contains a forbidden pattern, we modify it into a non-forbidden path as in Figure~\ref{fig_case_2_paths_modification}. 
Since $\breve{P}_{1,\bar{1}}$ is symmetric, we have $\breve{P}_{1,\bar{1}} = Q_{1,\bar{1}}$, which implies that $(\widetilde{P}_{1,\bar{1}})_{\rm ex} \neq Q_{1,\bar{1}}$,
Since $(\widetilde{P}_{1,\bar{1}})_{\rm ex}$ has a peak $t_{1,\bar{1}}$, we see that $(\widetilde{P}_{1,\bar{1}})_{\rm ex}$ is a non-canonical new symplectic rigorous path. \vs{0.1cm}

\begin{figure}[H]
\centering
\begin{tikzpicture}[scale = 0.55, yscale=0.6]
\def\nwires{4} 
\def\length{15} 

\foreach \i\y in {
1/1, 2/2, 3/4, 4/3, 5/4,
6/2, 7/3, 8/2, 9/1, 10/3, 
11/2, 12/3, 13/4, 14/3, 15/2
}{
	\pgfmathsetmacro\ycoordinate{\length-\i+1}
	\pgfmathsetmacro\yplusone{\y+1}
	\foreach \x in {1,...,\nwires}{
		\ifthenelse{\x = \y}
			{\coordinate (\i LSW) at (-\nwires+\x-0.5, \ycoordinate);
			\coordinate (\i LSE) at (-\nwires+\x +0.5, \ycoordinate);
			\coordinate (\i LNW) at (-\nwires+\x-0.5, \ycoordinate+1);
			\coordinate (\i LNE) at (-\nwires+\x +0.5,\ycoordinate+1);

			\coordinate (\i RSW) at (\nwires - \x -0.5, \ycoordinate );
			\coordinate (\i RSE) at (\nwires - \x + 0.5, \ycoordinate );
			\coordinate (\i RNW) at (\nwires - \x - 0.5, \ycoordinate +1);
			\coordinate (\i RNE) at (\nwires - \x + 0.5, \ycoordinate +1);

			\coordinate (\i LC) at (-\nwires + \x, \ycoordinate+ 0.5);
			\coordinate (\i RC) at (\nwires - \x, \ycoordinate+0.5);
		

			}{}


			}{}

}

\node[below] at (2.5,0) {$L_2$};
\node[below] at (-2.5,0) {$L_{\bar{2}}$};
\node[below] at (0.5,0) {$L_n$};
\node[below] at (-0.5,0) {$L_{\bar{n}}$};
\node[right] at (7RNE) {$\ell_u$};
\node[left] at (7LNW) {$\ell_{\bar{u}}$};
\draw[color=red, line width=1.5] 
(-3.5, 0.5)--(9LSW)--(7LNE)--(5LSW)--(4RNE)--(2RSW)--(1RNE)
(3.5, 0.5)--(9RSE)--(7RNW)--(5RSE)--(4LNW)--(2LSE)--(1LNW)
;
\draw[color=blue]
(-2.5,0.5)--(15LSW)--(11RNE)--(9RSW)--(9RNE)
(2.5,0.5)--(15RSE)--(11LNW)--(9LSE)--(9LNW)
;
%
\draw 
(0.5,0.5)--(14RSW)--(14RNE)--(12RSE)--(12RNW)--(10RSW)--(10RNE)--(8RSW)--(8RNE)--(6RSE)--(6RNW)--(4RSE)--(3RNW)
(-0.5,0.5)--(14LSE)--(14LNW)--(12LSW)--(12LNE)--(10LSE)--(10LNW)--(8LSE)--(8LNW)--(6LSW)--(6LNE)--(4LSW)--(3LNE)
;
\draw[color=teal]
(7RNE)--(7RSW)--(10RNW)--(10RSE)
(7LNW)--(7LSE)--(10LNE)--(10LSW)
;
\node[circle, draw, red, pin={[pin distance = 1.5cm]180:forbidden}] at (7LC) {};
\node[circle, draw, red, pin={[pin distance = 1.5cm]0:forbidden}] at (7RC) {};
\draw[color=red!50, line cap=round, line width=5, opacity=0.5, ->] 
(0.5,0.5)--(14RSW)--(14RNE)--(12RSE)--(12RNW)--(10RSW)--(10RC)--(10RNW)--(7RSW)--(7RC)--(7RNW)--(5RSE)--(5RC)
--(5LSW)--(7LNE)--(7LC)--(7LSE)--(10LNE)--(10LC)--(10LSE)--(12LNE)--(12LSW)--(14LNW)--(14LSE)--(-0.5,0.5); 

\draw[color=cyan!50!blue, line cap = round, line width = 2, dashed, ->] 
(3.5,0.5)--(9RSE)--(7RNW)--(5RSE)--(5RC)--(5RSW)--(7LNE)--(7LC)--(7LSE)--(10LNE)--(10LC)--(10LSE)--(12LNE)--(12LC)--(15RSE)--(2.5,0.5);

\end{tikzpicture}
\caption{In Case 2, $(\ell_n \to \ell_1 \to \ell_{\bar{1}} \to \ell_{\bar{n}})$ is modified to $(\ell_n \to \ell_u \to \ell_1 \to \ell_{\bar{1}} \to \ell_{\bar{u}} \to \ell_{\bar{n}})$;
$(\ell_1 \to \ell_{\bar{1}} \to \ell_{\bar{n}} \to \ell_2)$ is modified to $(\ell_1 \to \ell_{\bar{1}} \to \ell_{\bar{u}} \to \ell_{\bar{n}} \to \ell_2)$.
}\label{fig_case_2_paths_modification}
\end{figure}
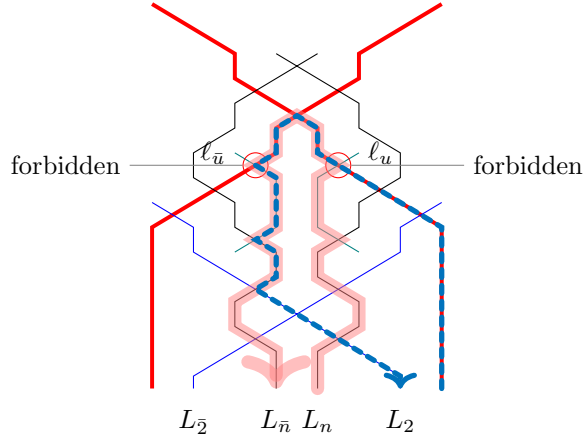

		\noindent
		{\bf Case 3:} $t_{k,\overline{k}}$ is below $\ell_1$ for all $2 \leq k \leq n$ (see Figure~\ref{fig_case_3_paths}). 
In this case, we label each node on $\ell_1$ above $\ell_{\bar{1}}$ by $t_2, \dots, t_n$ from right to left (see Figure~\ref{fig_nodes_on_ell_bar1}).

\begin{figure}[H]
\begin{tikzpicture}[scale = 0.55, yscale=0.6]
\def\nwires{5} 
\def\length{9} 

\foreach \i\y in {1/1, 2/2, 3/3, 4/4, 5/5, 6/4, 7/3, 8/2, 9/1}{
	\pgfmathsetmacro\ycoordinate{\length-\i+1}
	\pgfmathsetmacro\yplusone{\y+1}
	\foreach \x in {1,...,\nwires}{
		\ifthenelse{\x = \y}
			{\coordinate (\i LSW) at (-\nwires+\x-0.5, \ycoordinate+0);
			\coordinate (\i LSE) at (-\nwires+\x +0.5, \ycoordinate+0);
			\coordinate (\i LNW) at (-\nwires+\x-0.5, \ycoordinate+1);
			\coordinate (\i LNE) at (-\nwires+\x +0.5,\ycoordinate+1);

			\coordinate (\i RSW) at (\nwires - \x -0.5, \ycoordinate+0);
			\coordinate (\i RSE) at (\nwires - \x + 0.5, \ycoordinate + 0);
			\coordinate (\i RNW) at (\nwires - \x - 0.5, \ycoordinate + 1);
			\coordinate (\i RNE) at (\nwires - \x + 0.5, \ycoordinate + 1);

			\coordinate (\i LC) at (-\nwires + \x, \ycoordinate+ 0.5);
			\coordinate (\i RC) at (\nwires - \x, \ycoordinate+0.5);

}


}

}

\draw
(-4.5,11)--(1LNW)--(9RSE)--(4.5,0);

\foreach \x in {9,8,6,5}{
\draw (\x RNE)--(\x RSW);
}

\node[above] at (9RC) {$t_2$};
\node[above] at (8RC) {$t_3$};
\node[above] at (6RC) {$t_{n}$};

\node[above] at (7RC) {\rotatebox{150}{$\cdots$}};
%

\draw[teal!50!green, dashed] (0,1)--(0,10);

\node[below] at (4.5,-0.5) {$L_{1}$};

\end{tikzpicture}
\caption{Nodes on $\ell_{1}$ below $\ell_{\bar{1}}$.}\label{fig_nodes_on_ell_bar1}
\end{figure}
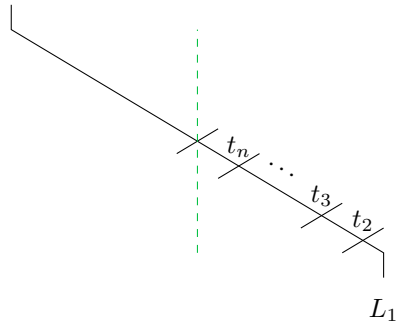

By the assumption of Case 3, we have $t_n = t_{1,\overline{k}}$ for some $2 \leq k \leq n$. 
It is easy to see that there exists a sequence $(k = r_1 > \cdots > r_q = 2)$ with $q \geq 1$ such that the following path on $G^{\symp}({\bm i},1)$ does not contain a forbidden pattern: 
\[\breve{P}_{1,\overline{k}} = (\ell_1 \to \ell_{\overline{k}} \to \ell_{r_1} \to \cdots \to \ell_{r_q}).\]
We also consider the following new symplectic rigorous path on $G^{\symp}({\bm i},k)$:
\[
\widetilde{P}_{1,\overline{k}} \colonequals (\ell_k \to \ell_{\bar{1}} \to \ell_{k+1}). 
\]
$\breve{P}_{1,\overline{k}}$ and $\widetilde{P}_{1,\overline{k}}$ are the red highlighted path and the blue dashed path in Figure~\ref{fig_case_3_paths}, respectively.
If $\widetilde{P}_{1,\overline{k}}$ contains a forbidden pattern, we modify it into a non-forbidden path as in Case 1.
Since $\breve{P}_{1,\overline{k}}$ does not cross the wall, it is obvious that $\breve{P}_{1,\overline{k}}$ satisfies the assumption of Proposition \ref{prop_non_redundancy_k_not_n_maximal}. 
Since $\chamber(\breve{P}_{1,\overline{k}}) \not\subseteq \chamber(\widetilde{P}_{1,\overline{k}}) \cup \chamber(\widetilde{P}_{1,\overline{k}}^\vee)$, it follows that $\breve{P}_{1,\overline{k}} \neq (\widetilde{P}_{1,\overline{k}})_{\rm ex} = Q_{k,\overline{1}}$. 
Since $\breve{P}_{1,\overline{k}}$ has a peak $t_{1,\overline{k}}$, we see that $\breve{P}_{1,\overline{k}}$ is a non-canonical new symplectic rigorous path. 
This completes the proof. 
	\end{proof}

\begin{figure}[H]
\centering
\begin{tikzpicture}[scale = 0.8, yscale=0.5]
\def\nwires{4} 
\def\length{14} 

\foreach \i\y in {
1/1, 2/2, 3/3, 4/4, 5/3,
6/4, 7/2, 8/3, 9/4, 10/1, 
11/2, 12/3, 13/2, 14/4
}{
	\pgfmathsetmacro\ycoordinate{\length-\i+1}
	\pgfmathsetmacro\yplusone{\y+1}
	\foreach \x in {1,...,\nwires}{
		\ifthenelse{\x = \y}
			{\coordinate (\i LSW) at (-\nwires+\x-0.5, \ycoordinate);
			\coordinate (\i LSE) at (-\nwires+\x +0.5, \ycoordinate);
			\coordinate (\i LNW) at (-\nwires+\x-0.5, \ycoordinate+1);
			\coordinate (\i LNE) at (-\nwires+\x +0.5,\ycoordinate+1);

			\coordinate (\i RSW) at (\nwires - \x -0.5, \ycoordinate );
			\coordinate (\i RSE) at (\nwires - \x + 0.5, \ycoordinate );
			\coordinate (\i RNW) at (\nwires - \x - 0.5, \ycoordinate +1);
			\coordinate (\i RNE) at (\nwires - \x + 0.5, \ycoordinate +1);

			\coordinate (\i LC) at (-\nwires + \x, \ycoordinate+ 0.5);
			\coordinate (\i RC) at (\nwires - \x, \ycoordinate+0.5);
		
			\draw (-\nwires+\x-0.5,\ycoordinate) -- (\i LSW)-- (\i LNE) -- (-\nwires+\x+0.5, \ycoordinate+1);
			\draw (-\nwires+\x+0.5,\ycoordinate) --(\i LSE) --(\i LNW) --(-\nwires+\x-0.5, \ycoordinate+1);

			\draw (\nwires - \x -0.5, \ycoordinate) -- (\i RSW)-- (\i RNE) -- (\nwires - \x + 0.5, \ycoordinate + 1);
			\draw (\nwires - \x +0.5, \ycoordinate) -- (\i RSE)-- (\i RNW) -- (\nwires - \x -0.5, \ycoordinate + 1);
			}{}

		\ifthenelse{\x=\y \OR \x = \yplusone}{}{
			\draw
			(-\nwires+\x-0.5,\ycoordinate)--(-\nwires+\x-0.5,\ycoordinate+1);
			\draw
			(\nwires-\x+0.5,\ycoordinate)--(\nwires-\x+0.5,\ycoordinate+1);
		}

			}{}


}

\node[below] at (2.5,0.5) {\small $L_2$};
\node[below] at (-2.5,0.5) {\small $L_{\bar{2}}$};
\node[below] at (0.5,0.5) {\small $L_{k+1}$};
\node[below] at (-0.5,0.5) {\small $L_{\overline{k+1}}$};
\node[below] at (1.5,0.5) {\small $L_k$};
\node[below] at (-1.5,0.5) {\small $L_{\bar{k}}$};
\draw[color=red, line width=1.5] 
(-3.5, 1)--(10LSW)--(10LNE)--(7LSW)--(7LNE)--(5LSW)--(1RNE)
(3.5, 1)--(10RSE)--(10RNW)--(7RSE)--(7RNW)--(5RSE)--(1LNW)
;
\draw[color=blue]
(-2.5, 1)--(13LSW)--(12LNE)--(9LSW)--(7RNE)--(2RSE)--(2RNW)
(2.5, 1)--(13RSE)--(12RNW)--(9RSE)--(7LNW)--(2LSW)--(2LNE)
;

\draw[color=red!50, line cap=round, line width=5, opacity=0.5, ->] 
(3.5,1)--(10RSE)--(10RNW)--(7RSE)--(7RNW)--(5RSE)--(5RC)--(6RC)--(6RSE)--(8RNW)--(8RSE)--(11RNW)--(11RSE)--(13RNE)--(13RC)--(13RSE)--(2.5,1);

\draw[color=cyan!50!blue, line cap = round, line width = 2, dashed, ->] 
(1.5,1)--(13RSW)--(13RNE)--(11RSE)--(11RNW)--(8RSE)--(8RNW)--(6RSE)--(5LC)--(5LSW)--(7LNE)--(7LSW)--(10LNE)--(10LC)--(12LSE)--(14RNW)--(14RSE);

\end{tikzpicture}
\caption{Case 3.}
\label{fig_case_3_paths}
\end{figure}

	\begin{lemma}\label{lemma_south}
		Suppose that there is a crossing on the south or the east sector of~$G^{\symp}({\bm i})$. 
		Then there is a non-canonical new symplectic rigorous path $P$ that gives a non-redundant inequality $\F_P^{(C)}({\bm a}) \geq 0$ in the expression \eqref{eq:descriptions_of_string_cones_in_type_C_n}.
	\end{lemma}
	
	\begin{proof}
		By Lemma \ref{lemma_wall}, we may assume that there is no crossing on the wall below $\ell_1$ and $\ell_{\bar{1}}$.
		We label each node on $\ell_{\bar{1}}$ above $\ell_1$ by $t_2, \dots, t_n$ from right to left (see Figure~\ref{fig_nodes_on_ell_1}). 

\begin{figure}[H]
\begin{tikzpicture}[scale = 0.55, yscale=0.6]
\def\nwires{5} 
\def\length{9} 

\foreach \i\y in {1/1, 2/2, 3/3, 4/4, 5/5, 6/4, 7/3, 8/2, 9/1}{
	\pgfmathsetmacro\ycoordinate{\length-\i+1}
	\pgfmathsetmacro\yplusone{\y+1}
	\foreach \x in {1,...,\nwires}{
		\ifthenelse{\x = \y}
			{\coordinate (\i LSW) at (-\nwires+\x-0.5, \ycoordinate+0);
			\coordinate (\i LSE) at (-\nwires+\x +0.5, \ycoordinate+0);
			\coordinate (\i LNW) at (-\nwires+\x-0.5, \ycoordinate+1);
			\coordinate (\i LNE) at (-\nwires+\x +0.5,\ycoordinate+1);

			\coordinate (\i RSW) at (\nwires - \x -0.5, \ycoordinate+0);
			\coordinate (\i RSE) at (\nwires - \x + 0.5, \ycoordinate + 0);
			\coordinate (\i RNW) at (\nwires - \x - 0.5, \ycoordinate + 1);
			\coordinate (\i RNE) at (\nwires - \x + 0.5, \ycoordinate + 1);

			\coordinate (\i LC) at (-\nwires + \x, \ycoordinate+ 0.5);
			\coordinate (\i RC) at (\nwires - \x, \ycoordinate+0.5);

}


}

}

\draw (-4.5,0)--(9LSW)--(9LNE)--(8LSW)--(8LNE)--(7LSW)--(7LNE)--(6LSW)--(6LNE)--(5LSW)--(5LNE)--(4RSW)--(4RNE)--(3RSW)--(3RNE)--(2RSW)--(2RNE)--(1RSW)--(1RNE)--(4.5,11);

\foreach \x in {1,2, 4,5}{
\draw (\x RNW)--(\x RSE);
}

\node[above] at (1RC) {$t_2$};
\node[above] at (2RC) {$t_3$};
\node[above] at (4RC) {$t_{n}$};

\node[above] at (3RC) {\rotatebox{30}{$\cdots$}};
%

\draw[teal!50!green, dashed] (0,1)--(0,10);

\node[below] at (-4.5,-0.5) {$L_{\bar{1}}$};

\end{tikzpicture}
\caption{Nodes on $\ell_{\bar{1}}$ above $\ell_1$.}\label{fig_nodes_on_ell_1}
\end{figure}
		Let $2 \leq j \leq n$ be the largest integer such that $t_j \neq \ell_j \cap \ell_{\bar{1}}$. Such an index $j$ exists by our assumption that 
		there exists a crossing on the south or the east sector. 
Write $t_j = \ell_k \cap \ell_{\bar{1}}$. We divide the proof into three cases. \vs{0.2cm}
		
		\noindent
		{\bf Case 1: $j=n$.} In this case, the symmetric symplectic rigorous path  
		\[
			\ell_n \rightarrow \ell_{\bar{1}} \rightarrow \ell_k \rightarrow \ell_{\bar{k}} \rightarrow \ell_1 \rightarrow \ell_{\bar{n}}
		\]
		is our desired path. See the blue dashed path in Figure~\ref{fig_case_1}.
Here, we note that the canonical symplectic rigorous path for $t_{n, \bar{1}}$ is $\ell_n \to \ell_{\bar{1}} \to \ell_1 \to \ell_{\bar{n}}$, which is highlighted in red in Figure~\ref{fig_case_1}.
\begin{figure}[H]
\begin{tikzpicture}[scale = 0.55, yscale=0.6]
\def\nwires{6} 
\def\length{17} 

\foreach \i\y in {
1/1, 2/2, 3/3, 
4/4, 5/6,
6/5, 
7/3, 
8/2, 9/6, 10/5,
11/4, 12/3, 13/2, 14/1, 15/3, 16/4, 17/5}{
	\pgfmathsetmacro\ycoordinate{\length-\i+1}
	\pgfmathsetmacro\yplusone{\y+1}
	\foreach \x in {1,...,\nwires}{
		\ifthenelse{\x = \y}
			{\coordinate (\i LSW) at (-\nwires+\x-0.5, \ycoordinate);
			\coordinate (\i LSE) at (-\nwires+\x +0.5, \ycoordinate);
			\coordinate (\i LNW) at (-\nwires+\x-0.5, \ycoordinate+1);
			\coordinate (\i LNE) at (-\nwires+\x +0.5,\ycoordinate+1);

			\coordinate (\i RSW) at (\nwires - \x -0.5, \ycoordinate );
			\coordinate (\i RSE) at (\nwires - \x + 0.5, \ycoordinate );
			\coordinate (\i RNW) at (\nwires - \x - 0.5, \ycoordinate +1);
			\coordinate (\i RNE) at (\nwires - \x + 0.5, \ycoordinate +1);

			\coordinate (\i LC) at (-\nwires + \x, \ycoordinate+ 0.5);
			\coordinate (\i RC) at (\nwires - \x, \ycoordinate+0.5);

			}{}


}
}

\draw[teal!50!green, dashed] (0,0.5)--(0,\length+1.5);
\draw[red, line width= 1.5] (5.5,\length+1.5)--(1RNE)--(1RSW)--(4RSW)--(6RNE)--(6RSW)--(9RNE)--(9RSW)--(10LNE)--(10LSW)--(11LNE)--(11LSW)--(12LNE)--(12LSW)--(13LNE)--(13LSW)--(14LNE)--(14LSW)--(-5.5,0)

(-5.5,\length+1.5)--(1LNW)--(1LSE)--(4LSE)--(6LNW)--(6LSE)--(9LNW)--(9LSE)--(10RNW)--(10RSE)--(11RNW)--(11RSE)--(12RNW)--(12RSE)--(13RNW)--(13RSE)--(14RNW)--(14RSE)--(5.5,0);


\draw 
(5RNE)--(5RSW)
(5LNW)--(5LSE);
\node[above] at (5RNW) {$\ell_k$};
\node[above] at (5LNE) {$\ell_{\bar{k}}$};

%

\node[below] at (0.5,-1) {$L_{n}$};
\node[below] at (-0.5,-1) {$L_{\bar{n}}$};

\draw[color=red!50, line cap=round, line width=5, opacity=0.5, ->] 
(0.5,0)--(17RSW)--(15RNE)--(13RSW)--(13RNE)--(8RSE)--(7RNW)--(4RSE)--(4RC)--(4RSW)--(6RNE)--(6RSW)--(9RNE)--(9RC)--(9RNW)--(6LSE)--(6LNW)--(4LSE)--(4LC)--(4LSW)--(7LNE)--(8LSW)--(13LNW)--(13LSE)--(15LNW)--(17LSE)--(-0.5,0);
\draw[color=cyan!50!blue, line cap = round, line width = 2, dashed, ->] 
(0.5,0)--(17RSW)--(15RNE)--(13RSW)--(13RNE)--(8RSE)--(7RNW)--(4RSE)--(4RC)--(4RSW)--(6RNE)--(6RC)--(5RC)--(6LC)--(6LNW)--(4LSE)--(4LC)--(4LSW)--(7LNE)--(8LSW)--(13LNW)--(13LSE)--(15LNW)--(17LSE)--(-0.5,0);

\node[above] at (6RC) {$t_{k,\bar{1}}$};
\end{tikzpicture}
\caption{Case 1: $j=n$.}\label{fig_case_1}
\end{figure}
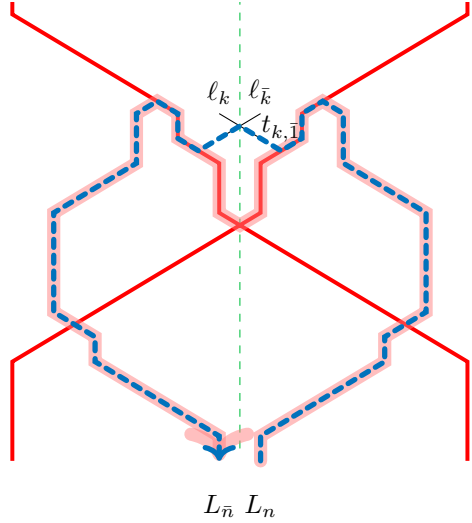

		\noindent
		{{\bf{Case 2:}} $j<n$ and $\ell_j \cap \ell_k$ is in the east sector.} In this case, 
		we have two new symplectic rigorous paths having the same peak $t_{j,\bar{1}}$: 
		\[
			\ell_j \rightarrow \ell_{\bar{1}} \rightarrow \ell_{j+1} \quad \text{and} \quad \ell_j \rightarrow \ell_{\bar{1}} \rightarrow \ell_k \rightarrow \ell_{j+1},
		\]		
	which satisfy the assumption of Proposition~\ref{prop_non_redundancy_k_not_n_maximal}.
The first path is the canonical symplectic rigorous path for $t_{j, \bar{1}}$, and the second path is not canonical. The canonical path is highlighted in red and the second path is blue dashed in Figure~\ref{fig_east}.
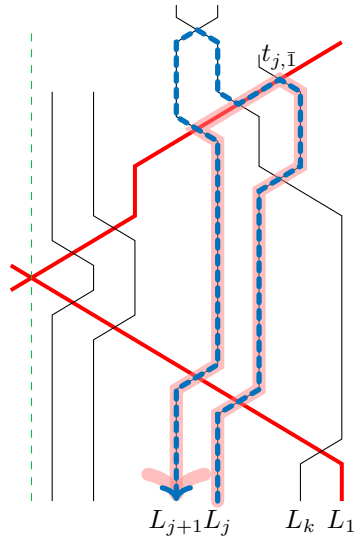
\begin{figure}[H]
\begin{tikzpicture}[scale = 0.55, yscale=0.6]
\def\nwires{8} 
\def\length{18} 

\foreach \i\y in {0/4, 1/1,2/2, 3/3, 4/4, 5/5, 6/2, 7/1, 8/6, 9/7, 10/8,
11/7, 12/6, 13/5, 14/4, 15/3, 16/2, 17/1}{
	\pgfmathsetmacro\ycoordinate{\length-\i+1}
	\pgfmathsetmacro\yplusone{\y+1}
	\foreach \x in {1,...,\nwires}{
		\ifthenelse{\x = \y}{

			\coordinate (\i RSW) at (\nwires - \x -0.5, \ycoordinate );
			\coordinate (\i RSE) at (\nwires - \x + 0.5, \ycoordinate );
			\coordinate (\i RNW) at (\nwires - \x - 0.5, \ycoordinate +1);
			\coordinate (\i RNE) at (\nwires - \x + 0.5, \ycoordinate +1);

			\coordinate (\i LC) at (-\nwires + \x, \ycoordinate+ 0.5);
			\coordinate (\i RC) at (\nwires - \x, \ycoordinate+0.5);
		
			}{}

}

}

\draw[red, line width = 1.5] (1RNE)--(1RSW)--(2RNE)--(2RSW)--(3RNE)--(3RSW)--(4RNE)--(4RSW)--(5RNE)--(5RSW)--(8RNE)--(8RSW)--(9RNE)--(9RSW)--(10RNE)--(10RSW);

\draw[red, line width = 1.5] (10RNW)--(10RSE)--(11RNW)--(11RSE)--(12RNW)--(12RSE)--(13RNW)--(13RSE)--(14RNW)--(14RSE)--(15RNW)--(15RSE)--(16RNW)--(16RSE)--(17RNW)--(17RSE)--(7.5,0.5);

\draw (0.5, 17)--(9RNW)--(9RSE)--(11RNE)--(11RSW)--(0.5,0.5);
\draw (1.5, 17)--(8RNW)--(8RSE)--(12RNE)--(12RSW)--(1.5,0.5);
\draw (4.5, 20.5)--(0RNE)--(0RSW)--(4RNW)--(4RSE)--(14RNE)--(14RSW)--(3.5,0.5);
\draw (3.5, 20.5)--(0RNW)--(0RSE)--(3RNW)--(3RSE)--(6RNW)--(6RSE)--(7RNW)--(7RSE)--(17RNE)--(17RSW)--(6.5,0.5);
\draw (5.5,18.5)--(2RNW)--(2RSE)--(6RNE)--(6RSW)--(15RNE)--(15RSW)--(4.5,0.5);

\node[below] at (4.5,0.5) {$L_j$};
\node[below] at (3.5,0.5) {$L_{j+1}$};
\node[below] at (6.5,0.5) {$L_k$};
\node[below] at (7.5,0.5) {$L_1$};

\draw[teal!50!green, dashed] (0,0.5)--(0,\length+1.5);

\draw[color=red!50, line cap=round, line width=5, opacity=0.5, ->] 
	(4.5,0.5)--(15RSW)--(15RNE)--(6RSW)--(6RNE)--(2RSE)--(2RC)--(2RSW)--(4RC)--(4RSE)--(14RNE)--(14RSW)--(3.5,0.5);

\draw[color=cyan!50!blue, line cap = round, line width = 2, dashed, ->] 
	(4.5,0.5)--(15RSW)--(15RNE)--(6RSW)--(6RNE)--(2RSE)--(2RC)--(3RC)--(3RNW)--(0RSE)--(0RC)--(0RSW)--(4RNW)--(4RSE)--(14RNE)--(14RSW)--(3.5,0.5);

\node[above] at (2RC) {$t_{j,\bar{1}}$};
\end{tikzpicture}
\caption{$j < n$ and $\ell_j \cap \ell_k$ is in the east sector.}\label{fig_east}
\end{figure}
		
		\noindent
		{\bf Case 3:} $j<n$ and $\ell_j \cap \ell_k$ is in the south sector. In this case, we similarly obtain 		
		two new symplectic rigorous paths having the same peak $t_{j,\bar{1}}$: 
		\[
			\ell_j \rightarrow \ell_{\bar{1}} \rightarrow \ell_{j+1} \quad \text{and} \quad \ell_j \rightarrow \ell_{\bar{1}} \rightarrow \ell_k \rightarrow \ell_{j+1},
		\]		
which satisfy the assumption of Proposition~\ref{prop_non_redundancy_k_not_n_maximal}.
The first path is the canonical symplectic rigorous path for $t_{j, \bar{1}}$, and the second path is not canonical. The canonical path is highlighted in red and the second path is blue dashed in Figure~\ref{fig_south}.

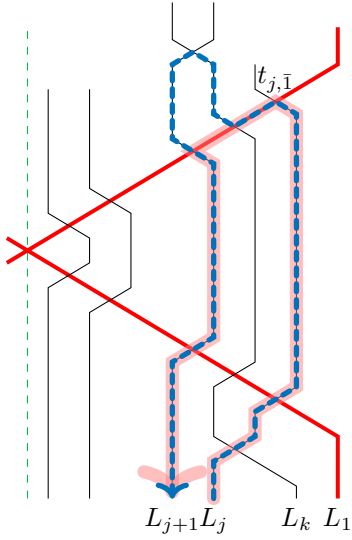
\begin{figure}[H]
\begin{tikzpicture}[scale = 0.55, yscale=0.6]
\def\nwires{8} 
\def\length{18} 

\foreach \i\y in {1/4, 2/1, 3/2, 4/3, 5/4, 6/5, 7/6, 8/7, 9/8, 10/7, 11/6, 12/5, 13/4, 14/3, 15/2, 16/1, 17/3, 18/2}{
	\pgfmathsetmacro\ycoordinate{\length-\i+1}
	\pgfmathsetmacro\yplusone{\y+1}
	\foreach \x in {1,...,\nwires}{
		\ifthenelse{\x = \y}{

			\coordinate (\i RSW) at (\nwires - \x -0.5, \ycoordinate );
			\coordinate (\i RSE) at (\nwires - \x + 0.5, \ycoordinate );
			\coordinate (\i RNW) at (\nwires - \x - 0.5, \ycoordinate +1);
			\coordinate (\i RNE) at (\nwires - \x + 0.5, \ycoordinate +1);

			\coordinate (\i LC) at (-\nwires + \x, \ycoordinate+ 0.5);
			\coordinate (\i RC) at (\nwires - \x, \ycoordinate+0.5);

			}{}

}

	\ifthenelse{\y=\nwires}{}{
}
}

\draw[teal!50!green, dashed] (0,0.5)--(0,\length+1.5);

\draw[red, line width = 1.5]
	(7.5, \length+1.5)--(2RNE)--(9RSW)
	(9RNW)--(16RSE)--(7.5,0.5);

\draw (0.5, 0.5)--(10RSW)--(10RNE)--(8RSE)--(8RNW)--(0.5, 17);
\draw (1.5,0.5)--(11RSW)--(11RNE)--(7RSE)--(7RNW)--(1.5,17);

\draw (3.5,0.5)--(13RSW)--(13RNE)--(5RSE)--(5RNW)--(1RSW)--(1RNE)--(4.5,20.5);
\draw (4.5,0.5)--(17RSW)--(17RNE)--(15RSW)--(15RNE)--(3RSE)--(3RNW)--(5.5,18);
\draw (6.5,0.5)--(18RSE)--(17RNW)--(14RSW)--(14RNE)--(4RSE)--(4RNW)--(1RSE)--(1RNW)--(3.5,20.5);

\node[below] at (4.5,0.5) {$L_j$};
\node[below] at (3.5,0.5) {$L_{j+1}$};
\node[below] at (6.5,0.5) {$L_k$};
\node[below] at (7.5,0.5) {$L_1$};

\draw[color=red!50, line cap=round, line width=5, opacity=0.5, ->] 
	(4.5,0.5)--(17RSW)--(17RNE)--(15RSW)--(15RNE)--(3RSE)--(3RC)--(5RC)--(5RSE)--(13RNE)--(13RSW)--(3.5,0.5);

\node[above] at (3RC) {$t_{j,\bar{1}}$};

\draw[color=cyan!50!blue, line cap = round, line width = 2, dashed, ->] 
	(4.5,0.5)--(17RSW)--(17RNE)--(15RSW)--(15RNE)--(3RSE)--(3RC)--(4RC)--(4RNW)--(1RSE)--(1RC)--(1RSW)--(5RNW)--(5RSE)--(13RNE)--(13RSW)--(3.5,0.5); 
\end{tikzpicture}
\caption{$j < n$ and $\ell_j \cap \ell_k$ is in the south sector.}\label{fig_south}
\end{figure}		
Hence the result follows.
	\end{proof}

Combining Lemmas \ref{lemma_new_canonical}, \ref{lemma_wall}, and \ref{lemma_south},
	we obtain the following corollary. 

	\begin{corollary}\label{corollary_strictly_increasing}
The equality holds in~\eqref{equation_increasing} if and only if $\bm i$ is of the form 
\[(\ldots,\underbrace{1,2,\dots,n,\dots,2,1}_{2n-1}).\]
	\end{corollary}

For $\bm i \in \RCn$ and a dominant integral weight $\lambda$, the string polytope $\Delta^{(C)}_{\bm i}(\lambda)$ in type $C_n$ is defined from $\mathcal{C}_{\bm i}^{(C_n)}$ in a way similar to Definition \ref{def_string_polytope}.
Let $P_{++}^{(C)}$ denote the set of regular dominant integral weights. 
Summarizing the above arguments, we know the following.

\begin{theorem}\label{thm_simplicial_string_cones_BC}
For $n \geq 2$ and $\bm i \in \RCn$, the following are equivalent.
\begin{enumerate}
\item The number of facets of $\Delta_{\bm i}(\lambda)$ is $2N$ for every $\lambda \in P_{++}^{(C)}$.
\item The string cone $\mathcal{C}_{\bm i}^{(C_n)}$ is simplicial.
\item $\bm i$ is either ${\bm i}_{C}^{(n)}$ or ${\bm j}_{C}^{(n)}$. 
\end{enumerate}
\end{theorem}

\begin{proof}
Fix $\lambda \in P_{++}^{(C)}$.
By the definition of the string polytope $\Delta^{(C)}_{\bm i}(\lambda)$, we see that the number of facets of $\Delta_{\bm i}^{(C)}(\lambda)$ is the number of facets of the string cone $\mathcal{C}_{\bm i}^{(C_n)}$ plus $N$. 
This provides the equivalence between (1) and (2). 
The equivalence between (2) and (3) comes from Corollaries \ref{corollary_simplicial} and \ref{corollary_strictly_increasing}. 
Hence the result follows.
\end{proof}

\begin{remark}\label{r:second_main_type_B}
By Theorem \ref{t:similarity_between_B_and_C}, we see that Theorem \ref{thm_simplicial_string_cones_BC} is naturally extended to the case of type $B_n$.
\end{remark}

Using Theorem~\ref{thm_simplicial_string_cones_BC}, we provide a classification of Gelfand--Tsetlin type string polytopes in type~$C_n$. 
Let $\lambda = \sum_{i=1}^n \lambda_i \varpi_i$ be a dominant integral weight, where $\varpi_1,\dots,\varpi_n$ are fundamental weights. 
The Gelfand--Tsetlin polytope $GT_{C_n}(\lambda)$ of type $C_n$ is defined to be the set of points
\[
(a^{(1)}_1, \underbrace{b^{(2)}_1, a^{(1)}_2, a^{(2)}_1}_3, 
\underbrace{b^{(3)}_1, b^{(2)}_2, a^{(1)}_3, a^{(2)}_2, a^{(3)}_1}_5,
\dots,
\underbrace{b^{(n)}_1,\dots,b^{(2)}_{n-1},a^{(1)}_n,\dots,a^{(n)}_1}_{2n-1}) \in \R^{N}
\]
satisfying the following inequalities. 
\[
\begin{tikzcd}[column sep = 0cm, row sep = 0cm]
\lambda_{\geq 1} \arrow[rd, "\textcolor{black}{{\geq}}" description,sloped,	color=white]   
	&& \lambda_{\geq 2} \arrow[rd, "\textcolor{black}{{\geq}}" description,sloped,	color=white] \arrow[ld, "\textcolor{black}{{\geq}}" description,sloped,	color=white] 
	&& \cdots 
	& \lambda_{\geq n} \arrow[rd, "\textcolor{black}{{\geq}}" description,sloped,	color=white] 
	&& 0 \arrow[ld, "\textcolor{black}{{\geq}}" description,sloped,	color=white] \\
& a^{(1)}_1 \arrow[rd, "\textcolor{black}{{\geq}}" description,sloped,	color=white] 
	&& a^{(1)}_2 \arrow[ld, "\textcolor{black}{{\geq}}" description,sloped,	color=white] && 
	& a^{(1)}_n \arrow[rd, "\textcolor{black}{{\geq}}" description,sloped,	color=white] \arrow[ld, "\textcolor{black}{{\geq}}" description,sloped,	color=white] \\
&& b^{(2)}_1 \arrow[rd, "\textcolor{black}{{\geq}}" description,sloped,	color=white] && \cdots 
	& b^{(2)}_{n-1} \arrow[rd, "\textcolor{black}{{\geq}}" description,sloped,	color=white] 
	&& 0 \arrow[ld, "\textcolor{black}{{\geq}}" description,sloped,	color=white] \\
&&& a^{(2)}_1 \arrow[rd, "\textcolor{black}{{\geq}}" description,sloped,	color=white] &&& a^{(2)}_{n-1} \\
&&&& \text{\rotatebox{-45}{$\cdots$}} \arrow[rd, "\textcolor{black}{{\geq}}" description,sloped,	color=white]  &&& \vdots \\
&&&&& b^{(n)}_1 \arrow[rd, "\textcolor{black}{{\geq}}" description,sloped,	color=white] 
	&& 0 \arrow[ld, "\textcolor{black}{{\geq}}" description,sloped,	color=white] \\
&&&&&& a^{(n)}_1
\end{tikzcd}
\]
Here, $\lambda_{\geq k} \coloneqq \lambda_k + \lambda_{k+1} + \cdots + \lambda_n$ for $1\leq k \leq n$.

\begin{theorem}[{\cite[Corollary~7]{Littelmann98}}]\label{t:Littelmann_type_C_GT}
For $n \geq 2$ and a dominant integral weight $\lambda$, the string polytope $\Delta_{{\bm i}_C^{(n)}}^{(C)}(\lambda)$ is unimodularly equivalent to the Gelfand--Tsetlin polytope $GT_{C_n}(\lambda)$ in type $C_n$. 
\end{theorem}

Let $\rho \coloneqq \varpi_1 + \varpi_2 + \cdots + \varpi_n \in P_{++}^{(C)}$ be the sum of fundamental weights.
As an application of Theorem \ref{thm_simplicial_string_cones_BC}, we can prove that the converse of Theorem \ref{t:Littelmann_type_C_GT} is also the case. 
More precisely, we obtain the following. 

\begin{theorem}\label{thm_GT_string_polytopes_C}
For $n \geq 2$ and $\bm i \in \RCn$, 
the string polytope $\Delta_{\bm i}^{(C)}(\rho)$ is unimodularly equivalent to the Gelfand--Tsetlin polytope $GT_{C_n}(\rho)$ in type $C_n$ if and only if ${\bm i} = {\bm i}_{C}^{(n)}$.
\end{theorem}

To prove Theorem~\ref{thm_GT_string_polytopes_C}, we prepare one lemma.

\begin{lemma}\label{lemma_non_integral_vertices}
Let ${\bm j}_n \coloneqq {\bm j}_{C}^{(n)}$.
Then the point $(0, \frac{3}{2}, 3,1,0,\dots,0) \in \R^N$ is a vertex of the string polytope~$\Delta_{{\bm j}_n}^{(C)}(\rho)$. Indeed, $\Delta_{{\bm j}_n}^{(C)}(\rho)$ is \emph{non}-integral.
\end{lemma}

\begin{figure}
\begin{subfigure}[t]{0.45\textwidth}
\centering
\begin{tikzpicture}[scale = 0.8, yscale=0.6]
\def\nwires{2} 
\def\length{4} 

\foreach \i\y in {1/1,2/2,3/1,4/2}{
	\pgfmathsetmacro\ycoordinate{\length-\i+1}
	\pgfmathsetmacro\yplusone{\y+1}
	\foreach \x in {1,...,\nwires}{
		\ifthenelse{\x = \y}
			{\coordinate (\i LSW) at (-\nwires+\x-0.5, \ycoordinate);
			\coordinate (\i LSE) at (-\nwires+\x +0.5, \ycoordinate);
			\coordinate (\i LNW) at (-\nwires+\x-0.5, \ycoordinate+1);
			\coordinate (\i LNE) at (-\nwires+\x +0.5,\ycoordinate+1);

			\coordinate (\i RSW) at (\nwires - \x -0.5, \ycoordinate );
			\coordinate (\i RSE) at (\nwires - \x + 0.5, \ycoordinate );
			\coordinate (\i RNW) at (\nwires - \x - 0.5, \ycoordinate +1);
			\coordinate (\i RNE) at (\nwires - \x + 0.5, \ycoordinate +1);

			\coordinate (\i LC) at (-\nwires + \x, \ycoordinate+ 0.5);
			\coordinate (\i RC) at (\nwires - \x, \ycoordinate+0.5);
		
			\draw (-\nwires+\x-0.5,\ycoordinate) -- (\i LSW)-- (\i LNE) -- (-\nwires+\x+0.5, \ycoordinate+1);
			\draw (-\nwires+\x+0.5,\ycoordinate) --(\i LSE) --(\i LNW) --(-\nwires+\x-0.5, \ycoordinate+1);

			\draw (\nwires - \x -0.5, \ycoordinate) -- (\i RSW)-- (\i RNE) -- (\nwires - \x + 0.5, \ycoordinate + 1);
			\draw (\nwires - \x +0.5, \ycoordinate) -- (\i RSE)-- (\i RNW) -- (\nwires - \x -0.5, \ycoordinate + 1);
			}{}

		\ifthenelse{\x=\y \OR \x = \yplusone}{}{
			\draw
			(-\nwires+\x-0.5,\ycoordinate)--(-\nwires+\x-0.5,\ycoordinate+1);
			\draw
			(\nwires-\x+0.5,\ycoordinate)--(\nwires-\x+0.5,\ycoordinate+1);
		}}

	\ifthenelse{\y=\nwires}{
		\node[above] at (\i LC) {$t_{\i}$};}
	{\node[above] at (\i LC) {$\bar{t}_{\i}$};
	\node[above] at (\i RC) {$t_{\i}$};}

}

\foreach \x in{1,...,\nwires}{
	\draw (\x-0.5,\length+1)--(\x-0.5,\length+1.5);
	\draw (\x-0.5,1)--(\x-0.5,0.5);

	\draw (-\x+0.5,\length+1)--(-\x+0.5,\length+1.5);
	\draw (-\x+0.5,1)--(-\x+0.5,0.5);
}

\end{tikzpicture}
\caption{$G^{\symp}({\bm j}_2)$}
\end{subfigure}
\begin{subfigure}[t]{0.45\textwidth}
\centering
\begin{tikzpicture}[scale = 0.8, yscale=0.6]
\def\nwires{3} 
\def\length{9} 

\foreach \i\y in {1/2,2/3,3/2,4/3,5/1,6/2,7/3,8/2,9/1}{
	\pgfmathsetmacro\ycoordinate{\length-\i+1}
	\pgfmathsetmacro\yplusone{\y+1}
	\foreach \x in {1,...,\nwires}{
		\ifthenelse{\x = \y}
			{\coordinate (\i LSW) at (-\nwires+\x-0.5, \ycoordinate);
			\coordinate (\i LSE) at (-\nwires+\x +0.5, \ycoordinate);
			\coordinate (\i LNW) at (-\nwires+\x-0.5, \ycoordinate+1);
			\coordinate (\i LNE) at (-\nwires+\x +0.5,\ycoordinate+1);

			\coordinate (\i RSW) at (\nwires - \x -0.5, \ycoordinate );
			\coordinate (\i RSE) at (\nwires - \x + 0.5, \ycoordinate );
			\coordinate (\i RNW) at (\nwires - \x - 0.5, \ycoordinate +1);
			\coordinate (\i RNE) at (\nwires - \x + 0.5, \ycoordinate +1);

			\coordinate (\i LC) at (-\nwires + \x, \ycoordinate+ 0.5);
			\coordinate (\i RC) at (\nwires - \x, \ycoordinate+0.5);
		
			\draw (-\nwires+\x-0.5,\ycoordinate) -- (\i LSW)-- (\i LNE) -- (-\nwires+\x+0.5, \ycoordinate+1);
			\draw (-\nwires+\x+0.5,\ycoordinate) --(\i LSE) --(\i LNW) --(-\nwires+\x-0.5, \ycoordinate+1);

			\draw (\nwires - \x -0.5, \ycoordinate) -- (\i RSW)-- (\i RNE) -- (\nwires - \x + 0.5, \ycoordinate + 1);
			\draw (\nwires - \x +0.5, \ycoordinate) -- (\i RSE)-- (\i RNW) -- (\nwires - \x -0.5, \ycoordinate + 1);
			}{}

		\ifthenelse{\x=\y \OR \x = \yplusone}{}{
			\draw
			(-\nwires+\x-0.5,\ycoordinate)--(-\nwires+\x-0.5,\ycoordinate+1);
			\draw
			(\nwires-\x+0.5,\ycoordinate)--(\nwires-\x+0.5,\ycoordinate+1);
		}}

	\ifthenelse{\y=\nwires}{
		\node[above] at (\i LC) {$t_{\i}$};}
	{\node[above] at (\i LC) {$\bar{t}_{\i}$};
	\node[above] at (\i RC) {$t_{\i}$};}

}

\foreach \x in{1,...,\nwires}{
	\draw (\x-0.5,\length+1)--(\x-0.5,\length+1.5);
	\draw (\x-0.5,1)--(\x-0.5,0.5);

	\draw (-\x+0.5,\length+1)--(-\x+0.5,\length+1.5);
	\draw (-\x+0.5,1)--(-\x+0.5,0.5);
}
\draw[color=red, line width = 1.5] 
(-2.5,0.5)--(9LSW)--(9LNE)--(8LSW)--(8LNE)--(7LSW)--(7LNE)--(6RSW)--(6RNE)--(5RSW)--(5RNE)--(2.5,10.5);
\draw[color=red, line width = 1.5] 
(2.5,0.5)--(9RSE)--(9RNW)--(8RSE)--(8RNW)--(7RSE)--(7RNW)--(6LSE)--(6LNW)--(5LSE)--(5LNW)--(-2.5,10.5);

\end{tikzpicture}
\caption{$G^{\symp}({\bm j}_3)$}\label{fig_G_symp_i3}
\end{subfigure}
\caption{Symplectic wiring diagrams for ${\bm j}_2 = (1,2,1,2)$ and ${\bm j}_3 = (2,3,2,3,1,2,3,2,1)$.}\label{fig_wiring_diagrams_i2_i3}
\end{figure}
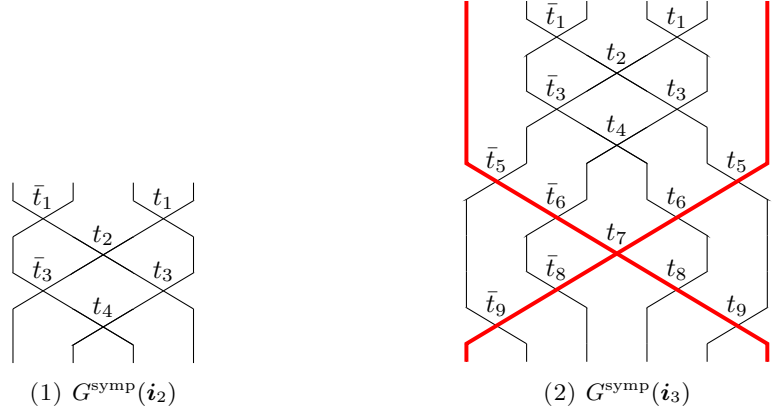

\begin{proof}
For $n = 2$, the string cone $\mathcal C_{\bm j_2}^{(C_n)}$ is given by the following inequalities: 
\begin{equation}\label{eq_C_i2}
a_1 \geq 0, \quad 
2a_2 \geq a_3 \geq 2 a_4 \geq 0. 
\end{equation}
By putting two wires on the bottom of the symplectic wiring diagram for ${\bm j}_2$, we obtain the symplectic wiring diagram for ${\bm j}_3$ as depicted in Figure~\ref{fig_wiring_diagrams_i2_i3}. 
These two added wires are shown in red in Figure~\ref{fig_G_symp_i3}. 
Considering symplectic rigorous paths, one sees that the string cone $\mathcal C_{{\bm j}_3}^{(C_n)}$ is given by the inequalities in~\eqref{eq_C_i2} and the inequalities $a_5 \geq a_6 \geq a_7 \geq a_8 \geq a_9 \geq 0$. 
Repeating this argument, we see that the string cone $\mathcal C_{{\bm j}_n}^{(C_n)}$ is given by the following inequalities:
\[
\begin{split}
&a_1 \geq 0, \\
&2a_2 \geq a_3 \geq 2 a_4 \geq 0, \\
&a_5  \geq  a_6 \geq a_7 \geq a_8 \geq a_9 \geq 0, \\
&\vdots  \\
&a_{N-2n+2} \geq a_{N-2n+3} \geq \cdots \geq  a_N \geq 0.
\end{split}
\]

In addition, the string polytope $\Delta_{{{\bm j}}_n}^{(C)}(\rho)$ is obtained from $\mathcal C_{{\bm j}_n}^{(C_n)}$ by the following inequalities:
\[
\begin{split}
a_1 & \leq 1 + 2 a_2 - 2a_3 + 2a_4 - \sum_{k \geq 5} \langle \alpha_{i_k}, h_{n-1}\rangle a_k,  \\
a_2 &\leq 1 + a_3 - 2a_4 - \sum_{k \geq 5} \langle \alpha_{i_k}, h_{n}\rangle a_k, \\
a_3 &\leq 1 + 2 a_4 - \sum_{k \geq 5} \langle \alpha_{i_k}, h_{n-1}\rangle a_k, \\
a_4 &\leq 1 - \sum_{k \geq 5} \langle \alpha_{i_k}, h_{n}\rangle a_k, \\
\vdots & \\
a_N & \leq 1.
\end{split}
\]
The point $(0, \frac{3}{2}, 3,1,0,\dots,0) \in \R^N$ becomes a vertex of the polytope $\Delta_{{\bm j}_n}^{(C)}(\rho)$ because it satisfies the following $N$ equalities and $N$ inequalities:

\[
\begin{split}
&\textcolor{red}{a_1 = 0}, \\
&\textcolor{red}{2a_2 = a_3}\, \geq 2 a_4 \geq 0, \\
&\textcolor{red}{a_5  =  a_6 = a_7 = a_8 = a_9 = 0}, \\
&\vdots  \\
&\textcolor{red}{a_{N-2n+2} = a_{N-2n+3} = \cdots =  a_N = 0},
\end{split} \qquad 
\begin{split}
a_1 & \leq 1 + 2 a_2 - 2a_3 + 2a_4 - \sum_{k \geq 5} \langle \alpha_{i_k}, h_{n-1}\rangle a_k,  \\
a_2 &\leq 1 + a_3 - 2a_4 - \sum_{k \geq 5} \langle \alpha_{i_k}, h_{n}\rangle a_k, \\
\textcolor{red}{a_3}  & \,\textcolor{red}{= 1 + 2 a_4 - \sum_{k \geq 5} \langle \alpha_{i_k}, h_{n-1}\rangle a_k}, \\
\textcolor{red}{a_4} &\,\textcolor{red}{= 1 - \sum_{k \geq 5} \langle \alpha_{i_k}, h_{n}\rangle a_k}, \\
\vdots & \\
a_N & \leq 1. 
\end{split}
\]
Hence the result follows.
\end{proof}

\begin{proof}[Proof of Theorem~\ref{thm_GT_string_polytopes_C}]
By Theorem~\ref{thm_simplicial_string_cones_BC}, there are only two possibilities of reduced words ${\bm i}$ such that the string polytope $\Delta_{\bm i}^{(C)}(\rho)$ has the same number of facets as the Gelfand--Tsetlin polytope~$GT_{C_n}(\rho)$. 
Because of Lemma~\ref{lemma_non_integral_vertices}, the string polytope $\Delta_{{\bm j}_C^{(n)}}^{(C)}(\rho)$ has a \emph{non-integral} vertex. Since the Gelfand--Tsetlin polytope $GT_{C_n}(\rho)$ is integral, this proves the statement.
\end{proof}

\begin{remark}
For $n = 3$ and ${\bm i} = {\bm j}_C^{(3)} = (2,3,2,3,1,2,3,2,1)$, the string polytope $\Delta_{\bm i}^{(C)}(\rho)$ is not combinatorially equivalent to the Gelfand--Tsetlin polytope $GT_{C_n}(\rho)$ since the $f$-vector of the Gelfand--Tsetlin polytope $GT_{C_n}(\rho)$ is $(1, 176, 936, 2244, 3126, 2760, 1590, 594, 138, 18, 1)$ while that of the string polytope $\Delta_{\bm i}^{(C)}(\rho)$ is $(1, 175, 933, 2241, 3125, 2760, 1590, 594, 138, 18, 1)$. Indeed, their numbers of vertices are different.
\end{remark}

\providecommand{\bysame}{\leavevmode\hbox to3em{\hrulefill}\thinspace}
\providecommand{\MR}{\relax\ifhmode\unskip\space\fi MR }
\providecommand{\MRhref}[2]{%
  \href{http://www.ams.org/mathscinet-getitem?mr=#1}{#2}
}
\providecommand{\href}[2]{#2}

\end{document}